\documentclass[letterpaper,10pt]{amsart}

\calclayout
\pdfoutput=1
\usepackage[english]{babel}

\usepackage[final]{microtype}
\usepackage{csquotes}
\usepackage[T1]{fontenc}
\usepackage{lmodern}

\usepackage{accents}
\usepackage{stackengine}
\usepackage{amsmath}
\usepackage{amsfonts}
\usepackage{amsthm}
\usepackage{mathtools}
\usepackage{graphicx}
\usepackage{mleftright}
\usepackage{amssymb}
\usepackage{dsfont}
\usepackage{xcolor}
\usepackage{bbm}
\usepackage[colorlinks=true, allcolors=blue]{hyperref}
\usepackage{subfiles}
\usepackage{cleveref}
\usepackage{tikz-cd}
\usepackage{enumitem}
\usepackage{nicematrix}
\usepackage{stmaryrd}
\usepackage{mathrsfs}

\usepackage{hyperref}

\begingroup
\catcode`\ =9
\catcode`\^^I=9
\endgroup

\newcommand{\SLipdf}{SL∞(Fq)}
\newcommand{\GLipdf}{GL∞(Fq)}

\newcommand{\redpdf}{red}

\newtheorem{theorem}{Theorem}[section]
\newtheorem{proposition}[theorem]{Proposition}
\newtheorem{corollary}[theorem]{Corollary}
\newtheorem{claim}[theorem]{Claim}
\newtheorem{lemma}[theorem]{Lemma}
\numberwithin{equation}{section}

\theoremstyle{definition}
\newtheorem{definition}[theorem]{Definition}
\newtheorem{remark}[theorem]{Remark}
\newtheorem*{example}{Example}
\newtheorem*{claim*}{Claim}

\newcommand{\p}{\varphi}

\renewcommand{\l}{\ell}
\newcommand{\CC}{\mathbb{C}}

\newcommand{\NN}{\mathbb{N}}
\newcommand{\ZZ}{\mathbb{Z}}
\newcommand{\RR}{\mathbb{R}}
\newcommand{\FF}{\mathbb{F}}
\newcommand{\one}{\mathbbm{1}}
\newcommand{\x}{\times}

\newcommand{\inter}{\cap}
\newcommand{\biginter}{\bigcap}
\newcommand{\union}{\cup}
\newcommand{\bigunion}{\bigcup}

\newcommand{\comp}[1]{{#1}^\mathsf{c}}

\newcommand{\angles}[1]{\left\langle #1 \right\rangle} 
\newcommand{\paren}[1]{\left( #1 \right)}
\newcommand{\abs}[1]{\left| #1 \right|}
\newcommand{\setr}[1]{\left\llbracket #1 \right\rrbracket}

\newcommand{\Fix}{\operatorname{Fix}}
\newcommand{\Gr}{\operatorname{Gr}}
\newcommand{\codim}{\operatorname{codim}}
\newcommand{\rk}{\operatorname{rk}}

\newcommand{\spn}{\operatorname{span}}
\newcommand{\diag}{\operatorname{diag}}
\newcommand{\id}{\operatorname{id}}
\newcommand{\sgn}{\operatorname{sgn}}
\newcommand{\Sym}{\operatorname{Sym}}
\newcommand{\Alt}{\operatorname{Alt}}
\newcommand{\proj}{\operatorname{proj}}

\newcommand{\Gal}{\operatorname{Gal}}

\newcommand*\diff{\mathop{}\!\mathrm{d}}

\newcommand{\set}[1]{{\left\{ \, #1 \, \right\}}}
\newcommand{\setcomp}[2]{\set{#1 \, : \, #2}}
\newcommand{\sset}[1]{{\{  #1  \}}}
\newcommand{\ssetcomp}[2]{\sset{#1  :  #2}}

\newcommand{\PP}{\mathbb{P}}
\newcommand{\EE}[1]{\mathbb{E}\mleft[#1\mright]}

\newcommand{\cA}{\mathcal{A}}
\newcommand{\cD}{\mathcal{D}}
\newcommand{\cE}{\mathcal{E}}
\newcommand{\cM}{\mathcal{M}}
\newcommand{\cN}{\mathcal{N}}
\newcommand{\cO}{\mathcal{O}}
\newcommand{\cP}{\mathcal{P}}

\newcommand{\cT}{\mathcal{T}}

\newcommand{\ubar}[1]{\stackunder[1.2pt]{$#1$}{\rule{.8ex}{.075ex}}}

\newcommand{\Sub}[1]{\operatorname{Sub}(#1)}
\newcommand{\IRS}{\operatorname{IRS}}
\newcommand{\Extr}{\operatorname{Extr}}
\newcommand{\Tr}{\operatorname{Tr}}
\newcommand{\tr}{\operatorname{tr}}
\newcommand{\Ch}{\operatorname{Ch}}
\newcommand{\Stab}{\operatorname{Stab}}

\newcommand{\Haar}{\operatorname{Haar}}
\newcommand{\triv}{{\set{I}}}

\newcommand{\WSub}[1]{\widehat{\operatorname{Sub}}(#1)}
\newcommand{\hmu}{\hat{\mu}}
\newcommand{\FP}{\operatorname{FP}}

\newcommand{\SL}{\operatorname{SL}}
\newcommand{\PSL}{\operatorname{PSL}}
\newcommand{\GL}{\operatorname{GL}}
\newcommand{\M}{\operatorname{M}}

\newcommand{\U}{\operatorname{U}}
\newcommand{\SU}{\operatorname{SU}}

\newcommand{\Sp}{\operatorname{Sp}}
\newcommand{\SO}{\operatorname{SO}}

\newcommand{\cB}{\mathcal{B}}
\newcommand{\FFq}{\FF_q}
\newcommand{\rkt}{\widetilde{\operatorname{rk}}}
\newcommand{\red}{\mathsf{red}}
\newcommand{\Col}{\operatorname{Col}}
\newcommand{\col}{\operatorname{color}}

\newcommand{\Ai}{{A_\infty}}

\newcommand{\GLj}[1]{{\GL_{#1}(\FFq)}}
\newcommand{\GLn}{\GLj{n}}
\newcommand{\GLm}{\GLj{m}}
\newcommand{\GLno}{\GLj{n-1}}
\newcommand{\GLi}{\GLj{\infty}}

\newcommand{\SLj}[1]{{\SL_{#1}(\FFq)}}
\newcommand{\PSLj}[1]{{\PSL_{#1}(\FFq)}}
\newcommand{\SLn}{{\SLj{n}}}

\newcommand{\SLi}{{\SL_\infty(\FFq)}}
\newcommand{\SLz}{{\SL_{\ZZ}(\FFq)}}

\newcommand{\Mj}[2]{\M_{#1 \x #2}(\FFq)}

\newcommand{\GLB}{{\operatorname{GLB}(\FFq)}}

\newcommand{\hV}{\widehat{V}}

\newcommand{\uv}{{\ubar v}}
\newcommand{\ux}{{\ubar x}}
\newcommand{\uw}{{\ubar w}}

\newcommand{\up}{{\ubar \p}}
\newcommand{\ua}{{\ubar a}}

\title[Invariant random subgroups of $\SLi$]{Invariant random subgroups of the finitary special linear group over a finite field}
\author{Noam Atar \and Arie Levit}

\begin{document}

\begin{abstract}
    An invariant random subgroup of a discrete group is a conjugation-invariant probability measure on the space of its subgroups.
We confirm a conjecture of Simon Thomas and classify the invariant random subgroups of the finitary  special and general linear groups over a finite field.
These are shown to be stabilizers of random  tuples of vectors and functionals. The classification uses the character theory of the finitary  special linear group established by Skudlarek and the theory of the invariant random subgroups of the finitary  permutation group established by Vershik. In addition, we classify the invariant probability measures for several natural actions of these finitary groups.

\end{abstract}

\maketitle

\section{Introduction}
\label{section:introduction}
Let $q$ be a prime power and $\FF_q$ the finite field of order $q$. Consider the  finitary special linear group $\SLi = \bigcup_n \SLn$. The group $\SLi$ consists of all infinite matrices $A$ over the field $\FFq$ satisfying $A_{ij}=\delta_{ij}$ for all but finitely many pairs of indices $i,j \in \NN$ as well as the requirement $\det A = 1$. 

Our main goal is to classify the invariant random subgroups of the group $\SLi$. This means the following. To any given discrete group $G$ we  associate its space of subgroups $\Sub{G}$. This is a compact space endowed with the Chabauty topology. The group $G$ is acting on $\Sub{G}$ by homeomorphisms via conjugation. An 
\textit{invariant random subgroup} of $G$ is a probability measure on $\Sub{G}$ preserved under this action. One example is a Dirac point mass $\delta_N$ at some normal subgroup $N \lhd G$. The most general example is the stabilizer distribution of some probability measure preserving action of $G$. 
We denote the convex set of all invariant random subgroups of $G$ by $\IRS(G)$. An invariant random subgroup of $G$ is called \emph{ergodic} if it is an extreme point of the convex compact set $\IRS(G)$. Invariant random subgroups were introduced in \cite{AGV2014}. See e.g.
\cite{gelander2015lectureinvariantrandomsubgroups} for a survey.

We show that the invariant random subgroups of  $\SLi$ can be described by means of certain natural probability measure preserving actions. Consider the the infinite-dimensional vector space $V = \prod_{n \in \NN} \FFq$ over the field $\FFq$. The group $\SLi$ admits the standard action on $V$, namely $A\cdot v = Av$.
In addition, let $\hV$ be the same underlying vector space as $V$, equipped  with the dual action of the group $\SLi$, namely $A \cdot \p = (A^{-1})^t \p$.

For each  $k, \l \in \NN_0= \NN \union \{0\}$  consider the vector space  $V^{k, \l} = V^k \x \hV^\l$ equipped with the natural product $\SLi$-action. 
Note that $V^{k, \l}$ is a compact abelian group and let $\Haar^{k, \l}$ be its Haar probability measure. 
The stabilizer map  $\Stab_{\SLi} : V^{k, \l} \to \Sub{\SLi}$ induces an invariant random subgroup of $\SLi$ denoted $\mu_{k, \l}$ as the pushforward of the measure $\Haar^{k, \l}$. Notice that $\mu_{0, 0} = \delta_{\SLi}$.
In addition, let $\mu_\infty = \delta_\triv$ denote the Dirac measure at the trivial subgroup. Our main result is the following.

\begin{theorem} \label{thm:main_theorem}
    The ergodic invariant random subgroups of the group $\SLi$ are  
    \[\IRS_\mathrm{ergodic}(\SLi) = \setcomp{\mu_{k, \l}}{k, \l \in \NN_0} \union \set{\mu_\infty}.\]
    These are pairwise distinct.
\end{theorem}

\Cref{thm:main_theorem} has been conjectured by Simon Thomas \cite{thomasbook}. One reading of \Cref{thm:main_theorem} is as a far-reaching generalization of the fact that the group $\SLi$ is simple and admits no finite-index subgroups. The weak-$*$ topology on the space $\IRS_\mathrm{ergodic}(\SLi)$ is such that $\mu_\infty$ is the only accumulation point of the points $\mu_{k,\l}$ (see \Cref{prop:extreme_points_topology}). In particular $\IRS(\SLi)$ is a Bauer simplex, c.f. \cite{glasner1997kazhdan}.

There is an alternative interpretation of each invariant random subgroup $\mu_{k,\l}$ as the stabilizer of a uniformly random point in the direct product of two infinite Grassmannian spaces $\Gr_k(V) \x \Gr_\l(\hV)$ (see \Cref{subsub:Grassmannians} for details). 

It is interesting to note a difference between the topological and measurable dynamics of the $\SLi$-action  on its Chabauty space of subgroups. 
Thomas \cite{thomas2018URS} has classified the minimal closed orbits (i.e. \textit{uniformly recurrent subgroups}) for this the action. The only such orbits are the trivial ones, namely the singleton $\set{N}$ where $N = \triv$ or $N = \SLi$.

\begin{example}
It is a consequence of Theorem \ref{thm:main_theorem} that the only meaningful information that can be \enquote{extracted} given a random tuple $\ux=(v^1, \ldots, v^k, \p^1, \ldots, \p^\l) \in V^{k,\l}$ of $k$ vectors and $\l$ functionals with the purpose of defining an invariant random subgroup is its pointwise stabilizer. At first glance, this might seem surprising. For instance, one might be tempted to consider the stabilizer of a flag defined by $\ux$ or of its projectivization. It turns out that all such ideas eventually lead to the same class of invariant random subgroups.
For example, let $v \in V^{1,0}$ be a generic vector. It is easy to see that $\Haar^{1,0}$-almost surely one has
\[ \setcomp{A \in \SLi}{ Av=v} = \setcomp{A \in \SLi}{ \text{$Av = \lambda v$ for some $\lambda \in \FFq$}}.\]
\end{example}

We also consider the finitary general linear group $\GLi = \bigcup_n \GLn$. It is defined analogously to the group $\SLi$ with the requirement that the determinant is non-zero. Certainly $\SLi$ is a normal subgroup in $\GLi$ of index $q-1$. 
Given a fixed
 multiplicative subgroup $K \le \FFq^\x$, a fixed pair  $k,\l \in \NN_0$ and a tuple $\ux \in V^{k,\l}$ there is an associated subgroup
\[
    \Stab_\GLi^K(\ux) = \setcomp{g \in \GLi}{\text{$g\ux = \ux$ \, and \, $ \det g \in K$}}.
\]
Pushing forward the measure $\Haar^{k,\l}$ via the map $\Stab_\GLi^K$ to the  space $\Sub{\GLi}$ gives rise to an  invariant random subgroup denoted $\mu_{k, \l}^K$. In addition, we denote $\mu_\infty = \delta_\triv$.

\begin{theorem} \label{thm:gl_inf_theorem}
    The ergodic invariant random subgroups of the group $\GLi$ are
    \[
    \setcomp{\mu_{k, \l}^K}{\text{$k, \l \in \NN_0$ and $K \leq \FFq^\x$ is a subgroup}} \union \set{\mu_\infty}.
    \]
    These are pairwise distinct.
\end{theorem}

\subsection{Classification of invariant measures on \texorpdfstring{$V^{k, \l}$}{V\textasciicircum k,l}}
As part of the proof of \Cref{thm:main_theorem} we are led to classify the $\SLi$-invariant probability measures on the vector spaces $V^{k, \l}$. We find this ergodic-theoretic result to be of independent interest. 
Here is the most essential case of this classification.

\begin{theorem} \label{thm:xkl_measures}
    Let $k, \l \in \NN_0$. Let $\mu$ be an $\SLi$-invariant probability measure on $V^{k, \l}$. If for $\mu$-almost every tuple $\ux=(v^1, \ldots, v^k, \p^1, \ldots, \p^\l) \in V^{k, \l}$ both collections $\set{v^1, \ldots, v^k}$ and $\set{\p^1, \ldots, \p^\l}$ are linearly independent then $\mu = \Haar^{k, \l}$.
\end{theorem}

The proof of the above \Cref{thm:xkl_measures} depends on the pointwise ergodic theorem for actions of amenable groups and on various linear-algebraic counting problems over the field $\FFq$.
See  \Cref{theorem:invariant_measures_on_v} for the full classification of all $\SLi$-invariant measures on $V^{k, \l}$. Furthermore, in the related \Cref{prop:grassmannian_measures} we show unique ergodicity for the product of the infinite Grassmannian spaces $\Gr_k(V) \x \Gr_\l(\hV)$. 

\begin{remark}
In the case where $q$ is a prime and either $k=0$ or $\l=0$, our Theorem \ref{thm:xkl_measures} follows from \cite[Theorem 4.1]{bowen2017characteristic} with a different proof. While it is likely that the proof method of \cite{bowen2017characteristic} can be extended also when $q$ is a prime power and possibly also for general $k$ and $\l$, our methods give the stronger  \Cref{prop:stabilizer_dual_uniquely_ergodic} and are somewhat more effective.
\end{remark}

\subsection{Thoma characters}

Traces and characters in the sense of Thoma play an important role in our work. A \emph{trace} on a given discrete group $G$ is a conjugation-invariant positive-definite function $\varphi : G \to \CC$ normalized so that $\varphi(e)=1$. If $\pi$ is a finite-dimensional representation of $G$  then the function $\frac{\operatorname{tr} \pi}{\dim \pi} $ is a Thoma trace.
A general Thoma trace can be viewed as the analogue of such a map when $\pi$ is allowed to be infinite-dimensional. We denote by $\Tr(G)$ the space of all Thoma traces on $G$. It is a compact convex set with the pointwise convergence topology. These are of independent representation theoretic interest, see  e.g.  \cite{bekka_la_harpe2020} for more details.

The extreme points of the convex set $\Tr(G)$ are called the \textit{Thoma characters} of $G$. The set of characters is denoted $\Ch(G)$.
Skudlarek has classified the Thoma characters of the group $\SLi$.
For each $r \in \NN_0$ define the function $\chi_r : \SLi \to [0, 1]$ by 
\[
\chi_k(g) = q^{-r \cdot\rk(g-I)} \quad \forall g \in \SLi.
\]
Note that the rank of the matrix $g-I$ is always finite, provided $g \in \GLi$.
In addition, we denote $\chi_\infty = \delta_I$.

\begin{theorem}[Skudlarek \cite{Skudlarek1976}] \label{thm:sl_infty_characters}
$\Ch(\SLi) = \setcomp{\chi_r}{r \in \NN_0} \union \set{\chi_\infty}$.
\end{theorem}

We point out that our proof of Theorem \ref{thm:main_theorem} relies crucially and in several ways on the a priori knowledge of Skudlarek's character classification for $\SLi$.

To an invariant random subgroup $\mu$ of a given discrete group $G$ we associate its \textit{Vershik trace} $\chi_\mu$ given by
\[\chi_\mu(g) = \mu(\setcomp{H \in \Sub{G}}{g \in H}) \quad \forall g \in G.\]
The Vershik trace $\chi_\mu$ is a trace in the sense of Thoma. 
A direct calculation shows that the Vershik trace of the invariant random subgroup $\mu_{k, \l}$ of the group $\SLi$ is the character $\chi_{r}$  where $r = k + \l $ (\Cref{prop:trace_calculation}).

The Vershik trace map $\IRS(G) \to \Tr(G)$ is an affine map of compact convex sets. It sets up a very general connection between invariant random subgroups and the  representation theory of the group $G$. One may be interested in studying the properties of this map, i.e.\ when is it surjective,  injective, or extreme point preserving?

The finitary symmetric and alternating groups $S_\infty = \Sym_\text{fin}(\NN)$ and $\Ai = \Alt_\text{fin}(\NN)$  provide a negative answer to all three of the above questions  (see \Cref{subsection:irs_a_infty_description} for a detailed discussion).
Even so, the Vershik trace map for these two groups is \enquote{close enough} to be a bijection and carries a great deal of information.
We will make essential use of this behavior in the proof of \Cref{thm:main_theorem}, as the alternating group $\Ai$ is embedded in the group $\SLi$ via permutation matrices.

In the case of the group $\SLi$,  the Vershik trace map takes ergodic invariant random subgroups of $\SLi$ to characters, and is surjective and finite-to-one on these.
This gives an interesting example where the Vershik trace map is not injective, but is surjective and extreme point preserving.

\subsection{Other groups and future directions}
The two groups $\Ai$ and $S_\infty$ can be seen as the $q=1$ analogues of the groups $\SLi$ and $\GLi$ respectively.   
The classification of their Thoma characters and invariant random subgroups has been known for a relatively long time. The characters were classified by Thoma in \cite{thoma1964unzerlegbaren}.
The invariant random subgroups were classified by Vershik in \cite{Vershik2012}.
The calculation of the Vershik traces and the description of these as a weighted average over the set of Thoma characters are done in \cite{thomas2022_sinfty_irs_and_characters}. 
These results will be described precisely in \Cref{subsection:irs_a_infty_description}.

For the permutation groups $S_\infty$ and $\Ai$, both the ergodic invariant random subgroups and Thoma characters form uncountable sets.
In addition, although the Vershik trace of an ergodic invariant random subgroup may not be a character, every character has a canonically associated ergodic invariant random subgroup (see \cite{thomas2022_sinfty_irs_and_characters}).
This turns out to be  very far from the situation for the group $\SLi$,  the starkest difference being that the number of  ergodic invariant random subgroups as well as characters is countable.

The fact that $\SLi$ and $\GLi$ have a countable number of characters implies that they have few von Neumann factorial representations (see \cite[Chapter 11]{bekka_la_harpe2020} for details, though this point of view will not be pursued in the sequel).
This has led Gorin, Kerov, and Vershik to consider the group $\GLB$ of infinite invertible matrices over $\FFq$ which have finitely many non-zero entries below the diagonal and no restriction on or above the diagonal. 
In \cite{Gorin_Kerov_Vershik2014} they show that $\GLB$ has a rich structure of semi-finite traces, i.e.\ traces that may attain the value $\infty$.
This suggests that the space of invariant random subgroups of $\GLB$ may be similarly interesting. To study the connection between invariant random subgroups and traces an appropriate generalization of the Vershik map to this non-discrete group is needed. 

The characters of several other families of infinite-dimensional groups are known, opening the door to a classification of  their invariant random subgroups. 
These include the infinite-dimensional real Lie groups $\U_\infty(\RR)$ \cite{Voiculescu1976, Boyer1983, VershikKerov_U_infty_1982}, $\SO_\infty(\RR)$ and $\Sp_\infty(\RR)$ \cite{Boyer1992}. 
Note that $S_\infty$ can be viewed as the Weyl group of $\GL_\infty(\RR)$. 
The characters of the other infinite Weyl groups are classified in \cite{HiraiHirai2003} and it may be of interest to classify their invariant random subgroups.

Lastly, we mention the infinite-dimensional Lie groups over a finite field, i.e. $\U_\infty(\FFq)$, $\SO_\infty(\FFq)$ and $\Sp_\infty(\FFq)$. 
The characters of these groups (with $q$ odd) have been classified in \cite{Leinen2007}.
Thomas \cite{thomas2018URS} has shown that $\Sp_\infty(\FF_{2^n})$ does admit  a non-trivial uniformly recurrent subgroup which might indicate that the theory of invariant random subgroups there might be more delicate. 

\subsection{Permutation stability}

A discrete group $G$ is called  \textit{permutation stable} if every almost-action of $G$ on a finite set is close to an honest action. The precise formulation of this notion is given in \Cref{def:permutation stable} below. 

The connection between invariant random subgroups and permutation stability for amenable groups was established by Becker, Lubotzky and Thom \cite{blt2019}. 
The precise statement is the following.
An invariant random subgroup of $G$ is called \textit{co-sofic} if it is a limit in the weak-$*$ topology of invariant random subgroups supported on finite index subgroups.

\begin{theorem}[{\cite{blt2019}}]\label{thm:blt_stability}
    A finitely generated amenable group $G$ is permutation stable if and only if every invariant random subgroup of $G$ is co-sofic.
\end{theorem}

The work {\cite{blt2019}} is one of the motivations for the study of invariant random subgroups of amenable groups.
For example, in \cite{levit_lubotzky2022}, the authors use \Cref{thm:blt_stability} to show that the uncountable family of the B.~H.~Neumann groups are permutation stable. The analysis of the invariant random subgroups of these groups crucially depends on the above-mentioned classification of the invariant random subgroups of $S_\infty$ by Vershik.
Our classification of the invariant random subgroups of the group $\SLi$ in 
\Cref{thm:main_theorem} can be used to establish the permutation stability of an analogous family of groups.

\begin{theorem}\label{thm:elementary-enrichment-stability}
    Let $\cP = (n_i)$ be a monotone sequence of odd integers with $n_1 \geq 3$. The diagonal product $\bigotimes_i \PSLj{n_i}$ is permutation stable.
\end{theorem}

The abstract definition of the diagonal product operation is involved, and we do not reproduce or invoke it here. See e.g.\ \cite{dogon2026characters} where the same group is defined and shown to be Hilbert--Schmidt stable. We only consider here a concrete and hands-on construction of the group $\bigotimes_i \PSLj{n_i}$ as a finitely generated dense subgroup of the profinite group $\prod_i \PSLj{n_i}$, analogous to the original definition of the B.~H.~Neumann groups as finitely generated dense subgroups of $\prod_i A_{n_i}$. See \Cref{section:elementary-enrichment-stability} for details on all this.

\subsection{Proof outline}
\label{sub:proof outline}
Our strategy to prove \Cref{thm:main_theorem} is as follows. Fix a non-trivial ergodic invariant random subgroup $\mu$ of $\SLi$. We find suitable integers $k, \l \in\NN_0 $ and construct a \emph{fixed-point map} 
\[
    \FP : \Sub{\SLi} \to V^{k, \l}
\]
which associates to a $\mu$-generic subgroup $H$ a tuple $ \FP(H) \in V^{k, \l}$ such that $H = \Stab_\SLi(\FP(H))$. 
Further, we show that if the subgroup $H$ is distributed according to $\mu$ then the distribution of the tuple $\FP(H)$ is $\Haar^{k, \l}$. Hence
\[
    \mu = (\Stab_\SLi)_* \FP_* \mu = (\Stab_\SLi)_* \Haar^{k, \l} = \mu_{k, \l}
\]
as required.
To construct the map $\FP$ we use Vershik's classification of the invariant random subgroups of the finitary alternating group $\Ai$.

Let us provide a bit more detail.
A \textit{finitary basis} is  a basis for the vector space $\bigoplus_{n \in \NN} \FFq$ that agrees with the standard basis in all but finitely many members.
Given each finitary basis $f$ let $i_f : \Ai \hookrightarrow \SLi$ be the embedding of the finitary alternating group via  permutation matrices in the basis $f$. 
The image of this embedding will be denoted $\Ai(f)$.
The map $\Sub{\SLi} \to \Sub{\Ai(f)}$ given by intersection induces an invariant random subgroup $\nu_f$ of $\Ai(f)$ by pushing forward the measure $\mu$.
The proof has three  conceptual stages:
\begin{enumerate}
    \item Identify the invariant random subgroup $\nu_f$ for each finitary base $f$.
    \item Use the joint distribution of the measures $\nu_f$ as $f$ ranges over all possible finitary bases $f$ to identify the parameters $k, \l \in \NN_0$ and to construct a candidate fixed point map $\FP : \Sub{\SLi} \to V^{k, \l}$.
    \item Show that the map $\FP$ satisfies $\mu = (\Stab_\SLi)_* \FP_* \mu$ and $\FP_* \mu = \Haar^{k, \l}$.
\end{enumerate}
Before elaborating further, we  need to  sketch the invariant random subgroup structure of the group $\Ai$.
Choose a probability vector $\alpha = (\alpha_i)_{i \in \NN}$. Color the set $\NN$ of natural numbers using (possibly) infinitely many colors in a random  i.i.d.\ manner according to the vector  $\alpha$. The probability that each natural number in $\NN$ has the color $i$ is $\alpha_i$. 
Let $H$ be the random subgroup of $\Ai$ which preserves this random coloring. 
The distribution of $H \in \Sub{\Ai}$ is an ergodic invariant random subgroup of $\Ai$. All ergodic invariant random subgroups of $\Ai$ are of this form (up to small refinement; the exact details are given in \Cref{subsection:irs_a_infty_description}).

By the above classification, the invariant random subgroup $\nu_f$ obtained by intersecting $\mu$ with the subgroup of permutation matrices in each finitary basis $f$ is associated with a random coloring of the set $\NN$.
We conclude that  a permutation matrix $\sigma \in \Ai(f)$ belongs to our chosen $\mu$-generic subgroup $H$ if and only if $\sigma$ preserves the (random) color vector.

For the invariant random subgroup $\mu_{k, \l}$ corresponding to some known choice of $k, \l \in \NN_0$, there is a natural description of this random coloring. 
If $H \in \Sub{\SLi}$ is the stabilizer of the tuple $\ux = (\uv, \up) = (\uv^1, \ldots, \uv^k, \up^1, \ldots, \up^\l) \in V^{k, \l}$ and $f$ is some finitary basis, identify our set of colors with $\FFq^{r}$ where $r = k+\l$.  Each natural number $n \in \NN$ will be assigned the color represented by the tuple
\begin{equation}\label{eq:intro_coloring}
    \ux_n = \left(\uv^1_n, \ldots, \uv^k_n, \up^1_n, \ldots, \up^\l_n\right) \in \FFq^{r}.
\end{equation}
Here the subscripts $v_n$ and $\varphi_n$ respectively stand for the $n$-th coordinate of the vector $v$ and the functional $\p$ in the finitary basis $f$.
Indeed, a permutation matrix $\sigma$ fixes the tuple $\ux$ (equivalently, the matrix $\sigma$ belongs to $H$) exactly when it preserves this coloring. 
We see  that the resulting invariant random subgroup  $\nu_f$ of the group $\Ai$ corresponds  to a coloring of the set $\NN$ in $q^{r}$ colors, where each color has an equal probability of $q^{-r}$.

Let us return to the analysis of a general ergodic invariant random subgroup $\mu$. 
To show that $\mu$ is of the form $ \mu_{k, \l}$, we reverse the analysis of the above paragraph.
In stage (1) we show that for each finitary basis $f$, the invariant random subgroup $\nu_f$ of the group $\Ai$ corresponds to a  uniform coloring using $q^r$ (in particular, only finitely many)  colors for some parameter $r \in \NN_0$.
In stage (2) for each $\mu$-generic subgroup $H \in \Sub{\SLi}$ we find 
\enquote{candidate} values  $k, \l \in \NN_0$ and a tuple $\ux \in V^{k, \l}$ such that $k + \l = r$ and such that the associated coloring is of the form \cref{eq:intro_coloring} for \textit{all} finitary bases $f$. The correspondence $H \mapsto \ux$ gives rise to the map $\FP$ discussed above.
Lastly, we conclude with stage (3).

\subsection{Organization}

In \Cref{sec:the group SLi} we define the finitary special and general linear groups and consider their properties and subgroup structure. This section contains all the algebraic information to be used throughout the work. In \Cref{section:analysis_of_action} we prove abstractly that any action of $\GLi$  satisfying certain axioms is isomorphic to the action on $V^{k, \l}$.

The discussion of the main ideas of the work begins in earnest 
in \Cref{section:notations_and_preliminaries}, where we recall  the definitions of invariant random subgroups and traces, generally speaking and for the finitary alternating group $\Ai$ in particular.
In
\Cref{section:singer_cycles} we use Skudlarek's classification of the characters of the group $\SLi$ to analyze its  invariant random subgroups  (\Cref{sub:IRS contain SLn}) and implement stage (1) of the strategy outlined in  
\Cref{sub:proof outline}
(\Cref{subsection:irs_restriction}).
In 
\Cref{section:generalities_on_colors} we introduce the theory of the ``coloring equivalence relation'', a random coloring given by an invariant random subgroup $\mu$. This gives rise to  an action of the group $\GLi$ on vectors of colors (i.e. the equivalence classes of the coloring relation), which satisfies the axiomatic properties studied in \Cref{section:analysis_of_action}.

\Cref{section:measures_on_v} is standalone and involves ergodic theory. In that section  we prove the measure classification result \Cref{thm:xkl_measures}. As a corollary, we classify the $\SLi$-invariant measures on $V^{k, \l}$ and on the infinite Grassmannian spaces (\Cref{subsection:invariant_measures_on_related_spaces}). In 
\Cref{section:proof_end} we define the map $\FP$ and finish the proof of our main results Theorems \ref{thm:main_theorem} and
\Cref{thm:gl_inf_theorem}. 
This section contains a discussion of topological  structure of the  space of ergodic invariant random subgroups of $\SLi$ and $\GLi$. 
Finally, in 
\Cref{section:elementary-enrichment-stability} we recall the notion of permutation stability and use the classification of $\IRS(\SLi)$ to prove \Cref{thm:elementary-enrichment-stability}, that is, show that the diagonal products $\bigotimes \PSLj{n_i}$ are permutation stable.

\subsection{Standing assumptions and notations}
All our probability measures are Borel.

We let $\NN_0 = \NN \cup \{0\}$ denote the set of non-negative integers. We let $\NN_{>n}$ denote  the set of natural numbers greater than $n$. 
Let $\FFq$ denote the finite field of order $q$. We let $0_{n,k}$ denote the all-zeros matrix in $M_{n,k}(\FFq)$.

We use the following two standard asymptotic notations. The notation $o(1)$ stands for  arbitrary function which decays to zero as $n \to \infty$. The notation $O(n)$  stands for an arbitrary function $f$ satisfying $\abs{f(n)} \leq c n$ for some constant $c > 0$ which is independent of $n$ (though it may depend on other parameters, such that the field $\FFq$ in our case).
 
\subsection{Grant support}

AL is supported by the ISF grants  1788/22 and 3553/26 and the BSF grant 2022105. 

\subsection{AI statement}
Most of the research was done prior to the year 2026. No AI was used towards the research. AI was used to look up relevant references, including properties of Singer cycles.

\section{The finitary special linear group}
\label{sec:the group SLi}

Let $q$ be a prime power and $\FFq$ be a finite field. In this section we introduce the finitary special and general linear groups $\SLi$ and $\GLi$. The notations introduced here will be maintained and useful to us in the forthcoming sections. We also discuss various algebraic topics from finite group theory, such as Singer cycles and maximal subgroups.

\subsection{The groups $\SLi$ and $\GLi$}
Consider the $\FFq$-vector space $\bigoplus_{n \in \NN} \FFq$ equipped with the standard ordered basis $e = (e^n)_{n=1}^\infty$. A \textit{finitary basis} for the $\FFq$-vector space $\bigoplus_{n \in \NN} \FFq$ is  any ordered basis $f = (f^n)_{n=1}^\infty$ satisfying $f^n = e^n$
for all but finitely many indices $n \in \NN$. 

The finitary general linear group $\GLi$ consists of all invertible $\FFq$-linear transformations of $\bigoplus_{n \in \NN} \FFq$ fixing all but finitely many members of some (equivalently, any) finitary basis. 
The determinant map is defined naturally on $\GLi$. We denote $\det : \GLi \to \FFq^\x$ and let $\SLi = \ker \det$. Clearly $\SLi$ is a finite-index normal subgroup of $\GLi$ with quotient $\FFq^\x$.

An element of the group $\GLi$ can be represented in any finitary basis by an infinite matrix $A$ over the field $\FFq$ with coordinates indexed by $\NN$ such that $A_{ij} = \delta_{ij}$ for all but finitely many pairs $(i,j) \in \NN \times \NN$. An element of the group $\SLi$ can be likewise represented, with the additional requirement it has determinant one.

Let $f$ be a finitary basis and $I \subset \NN$ an index set. We define $\SL_I(f)$ to be the subgroup of all matrices $g \in \SLi$ keeping invariant the subspace ${\spn_{\FFq}}\setcomp{f^i}{i \in I}$  and fixing pointwise the complementary subspace ${\spn_{\FFq}}{\setcomp{f^i}{i \in \NN \setminus I}}$. An element of $\SL_I(f)$ is said to be \textit{supported on the index set $I$ in basis $f$}. 
The group $\SLi$ is isomorphic to the direct limit
\[ \SLi = \varinjlim_{I \subset \NN, |I| <\infty} \SL_I(f).\]
We will sometimes write $\SL_n(f)$ for $\SL_{\set{1, \ldots, n}}(f)$. If the index set $I \subset \NN$ is finite with $|I| =n $ then we obtain a natural embedding
\[\SLn \to \SL_I(f) \subseteq \SLi
\]
which is an isomorphism onto its image. The image of a matrix $A \in \SLn$ under this embedding with respect to the standard basis $e$ will be denoted by $[A]_I$. If the index set $I \subset \NN$ is infinite then  the subgroup $\SLj{I}$ is isomorphic to the group $\SLi$. 

All the notations and remarks in the previous paragraph  apply mutatis mutandis to the general linear group $\GLi$.

For each finitary basis $f$ there is an embedding $i_f$ of the finitary alternating group into the group $\SLi$ by permutation matrices acting in the basis $f$. We denote the image of this embedding by $\Ai(f) = i_f(\Ai)$. This embedding is an isomorphism onto its image.
Given  a subset $I \subset \NN$ we denote by $A_I(f)$ the subgroup of $\Ai(f)$ consisting of permutation matrices which are supported on $\setcomp{f^n}{n\in I}$.

For each element $g \in \GLi$ the transformation $g - I$ has finite rank on the vector space $\bigoplus_{n \in \NN} \FFq$. 
This motivates the following handy notation: we write \[\rkt(g) = \rk(g - I) \quad \forall g \in \GLi.\] 
Note that $\rkt(g)$ is zero if and only if $g$ is the identity transformation. If $g$ is not the identity then $\rkt(g)$ equals the size of the smallest index set $I$ which supports $g$ in some finitary basis.

Lastly, 
we set some common notation for specifc matrices.
Let $n \geq 2$ and  $1 \le i, j \leq n$ be a pair of distinct indices. Let $\sigma_{ij} \in \GLn$ be the permutation matrix swapping the $i$-th and $j$-th coordinates.
Let $E_{ij} \in \GLn$ be the elementary matrix with $1$ in its $(i, j)$-th entry.
For any scalar $\alpha \in \FFq^\x$ we let $D_i(\alpha) = \diag(1, \ldots, \alpha, \ldots, 1) \in \GLn$ be a diagonal matrix,
where $\alpha$ is in the $i$-th position. 

\subsection{The standard and dual actions} 
\label{subsec:standard and dual}

We let $V$ and $\widehat{V}$  denote the vector space $\prod_{n \in \NN} \FFq$ equipped respectively with the  standard and the dual actions of the finitary group $\GLi$.  The standard action is given by $A \cdot v = Av$ and the dual action is given by $A \cdot \varphi = (A^{-1})^t \varphi$ for all $A \in \GLi, v \in V$ and $\varphi \in \widehat{V}$. We often restrict these actions to the subgroup $\SLi$. Note that the underlying space for both actions is the same, and indeed, the two actions are equivariantly isomorphic. Further, note that $V$ and $\widehat{V}$ are compact profinite abelian groups equipped with a Haar probability measure.

For each fixed $k,\l \in \NN_0$ we let $V^{k, \l}$ denote the product space $V^k \x \hV^\l$ with the suitable product $\GLi$-action.
This space is $\GLi$-equivariantly isomorphic to the matrix product space $\M_{\infty \x k}(\FFq) \x \M_{\l \x \infty}(\FFq)$ with an action by matrix multiplication on the left in the first factor and on the right in the second factor. The two actions on $V^{k,\l}$ and $V^{\l,k}$ are $\GLi$-equivariantly isomorphic.

We set up some additional notation. For each $n \in \NN$ set $V_n = \FFq^n$ and $V_{>n} = \prod_{i=n+1}^\infty \FFq$. The space $V$ splits as a product $V = V_n \x V_{> n}$. 
For each vector $v \in V$ we denote by $v_{\leq n} \in V_{n}$ and $v_{> n} \in V_{> n}$ its corresponding projections. We also use analogous notations for the space $\hV$ and write $\hV = \hV_n \x \hV_{>n}$. We view the space $\hV_n$ as the dual of the space $V_n$.
We will use the notations $v_{\leq n}$ and $v_{> n}$ also for vectors $v$ in the finite-dimensional spaces $\FFq^m$ (provided $m > n$ of course).

It is natural to interpret elements in $\hV$ as functionals on the vectors in $V$. Keeping in line with this interpretation, we use the letters $v$ and $\varphi$ for vectors in $V$ and in $\widehat V$ respectively.
While a global pairing $V \x \hV \to \FFq$ is not well-defined, there do exist finite pairings
\[
\angles{\cdot, \cdot}_n : V \x \hV \to \FFq,
\quad \angles{v, \p}_n = \angles{v_{\leq n}, \p_{\leq n}}
\]
for each $n \in \NN$. Observe that the action of the finite group $\GLm$ preserves the pairings $\angles{\cdot, \cdot}_n$ for all $n \geq m$, that is
$\angles{g v, g \p}_n = \angles{v, \p}_n$ for any element $g \in \GLm$. 
For each tuple  $\ux = (\uv,\up) \in V^{k, \l}$ and $n \in \NN$, we define the pairing matrix $M_n(\uv, \up) \in \M_{\l \x k}(\FFq)$ by
\[
    (M_n(\uv, \up))_{i j} = \angles{\uv^j, \up^i}_n.
\]

\subsection{Extension by permutations} 
\label{sub:extension by permutations}
Let $H \leq \GLi$ be a fixed subgroup. 
An element $g \in \GLi$ is \textit{extendable by permutations to $H$} if there exists a subset of indices $I \subset \NN$ with $|I| = n$  and a matrix $A \in \GLn$ with $g = \left[A\right]_I$ such that for every sufficiently large $m \in \NN$ there is an additional subset of indices $J \subset \NN$ with $|J| =k $ and $\max I < m < \min J$ as well as a permutation matrix  $P \in \GLj{k}$ satisfying
\[ \left[\begin{pmatrix}
    A & 0 \\
    0 & P
\end{pmatrix}\right]_{I \cup J} \in H.
\]
The subset of all matrices $g \in \GLi$ extendable by permutations to $H$ is a subgroup of $\GLi$. We denote this subgroup by $H^*$.
Notice that $H \subseteq H^*$. The operation taking a subgroup $H$ to the subgroup $H^*$ is $\GLi$-equivariant, where $\GLi$ is acting by conjugation on subgroups.

\subsection{Singer cycles}\label{subsection:singer_cycles}
We discuss Singer cycles, a special type of group element. These are well-known in the theory of finite linear groups.
See \cite{Hestenes1970SingerGroups} for a general overview. 

\begin{definition}\label{def:singer_cycles}
    Let $q$ be a prime power. A Singer cycle in $\GL_n(\FFq)$ is an element of order $q^n-1$. 
    A Singer cycle in $\SL_n(\FFq)$ is an element of order $\frac{q^n-1}{q-1}$.
\end{definition}

We list basic properties of Singer cycles.
A reference for these can be found in \cite{Hestenes1970SingerGroups}.
Singer cycles have a maximal order in their corresponding group.
In the $\GL_n$ case, these are exactly the semisimple matrices whose eigenvalues are Galois conjugate primitive elements of $\FF_{q^n}$.
If $s \in \GLn$ is a Singer cycle then the unital ring it generates in $\M_n(\FFq)$ is isomorphic to the finite field $\FF_{q^n}$.
The centralizer of $s$ in $\M_n(\FFq)$ consists only of the matrices which lie in this ring. Hence the size of the conjugacy class of $s$ in $\GLn$ is
\begin{equation}
    \abs{\operatorname{Cl}_{\GLn} (s)} = \frac{\abs{\GLn}}{q^n-1} = q^{n^2 + O(n)}.
\end{equation}
 
For $\SL_n$, Singer cycles are exactly the semisimple matrices whose eigenvalues are Galois conjugates, each of which generates the kernel of the norm $\FF_{q^n}^\x \to \FFq^\x$.
The size of the conjugacy class of a Singer cycle in $\SL_n$ is  
\begin{equation}\label{eq:singer_cycle_cc_size}
    \abs{\operatorname{Cl}_{\SLn} (s)} = \frac{\abs{\SLn}}{(q^n-1)/(q-1)} = q^{n^2 + O(n)}.
\end{equation}

Singer cycles have been useful in the study of finite linear groups. In particular, the types of maximal subgroups of $\GLn$ or $\SLn$ containing a Singer cycle have been classified. 
See \cite{kantor1980linear, guralnick1999linear, Guralnick_Malle_2012} for this theory and for subgroups containing related types of elements.
We will need the following consequence of this classification (we encourage the reader interested in the fastest route towards our main result Theorem \ref{thm:main_theorem} to skip the proof and move onward to Sections \ref{section:analysis_of_action} and \ref{section:notations_and_preliminaries} upon first reading).

\begin{proposition}\label{prop:singer_cycle_subgroup_size}
    Let $M $ be a proper subgroup of $\SLn$ containing a Singer cycle. Then
    \[
    \abs{M} \leq q^{\frac 1 2 n^2 + O(n)}.
    \]
\end{proposition}

\begin{proof}
    By Lemma 2.3 in \cite{Guralnick_Malle_2012}, when $n > 20$ the proper maximal subgroups of $\SLn$
    containing a Singer cycle are conjugate to one of the following groups:
    \begin{enumerate}
        \item The normalizer of the special unitary group $\SU_n(\FF_{q^{1/2}})$ for $n$  odd and $q$  a square.
        \item The normalizer of the derived subgroup of the special orthogonal group $\Omega_n^-(\FFq)$ for $n$ even and $q$  odd.
        \item The normalizer of the symplectic group $\Sp_n(\FFq)$ for $n$ even.
        \item The normalizer of $\GL_{n/f} (\FF_{q^f}) \cap \SLn$ for $f$ a prime divisor of $n$, where we view the group $\GL_{n/f} (\FF_{q^f})$ as a subgroup of $\GLn$ via restriction of scalars.
    \end{enumerate}
    The definitions of these groups are given in \cite{Kleidman_Liebeck_1990}.
    We bound the sizes of these normalizers using known orders of linear groups.    
    Let $K$ be either $\SU_n(\FF_{q^{1/2}})$, $\Omega_n^-(\FFq)$ or $\Sp_n(\FFq)$. By Corollary 2.10.4 and Proposition 2.10.6 in \cite{Kleidman_Liebeck_1990}
    we have $N_{\GLn}(K) \leq \Delta$, where $\Delta \leq \GLn$ is the group of matrices that preserve the appropriate
    bilinear form up to a scalar.
    Table 2.1.C there gives $\Delta$'s size in each case, and we can see that it is $q^{\frac 1 2 n^2 + O(n)}$. Hence 
    \[
        \abs{N_\SLn (K)} \leq \abs{N_\GLn (K)} \leq \abs{\Delta} = q^{\frac 1 2 n^2 + O(n)}.
    \]
    Now let $K = \GL_{n/f} (\FF_{q^f}) \cap \SLn$.
    Then  $N_{\GLn}(K) \cong \GL_{n/f}(\FF_{q^f}) \rtimes \Gal(\FF_{q^f} / \FFq)$.
    In this case, we obtain again
    \[
        \abs{N_{\SLn} (K)} \leq \abs{\Gal(\FF_{q^f} / \FFq)} \cdot \abs{ \GL_{n/f}(\FF_{q^f})} = f (q^f)^{\left( \frac{n}{f} \right)^2 + O(n/f)} = q^{\frac{n^2}{f} + O(n)} \leq q^{\frac 1 2 n^2 + O(n)},
    \]
    as required.
\end{proof}

\subsection{On maximal subgroups of $\GLn$}

We state some results in the theory of finite general linear groups and their subgroups. These are to be used in \Cref{section:analysis_of_action} below. 

\begin{proposition} \label{prop:gl_n_1_overgroups}
    Suppose $n \geq 3$. Let $S$ be a proper subgroup of $\GLn$ containing the subgroup
    \[
    \GLno = \begin{pmatrix} 1 & 0_{1, n-1} \\ 0_{n-1, 1} & \GLno  \end{pmatrix}.
    \]
    Then there is a subgroup $K \leq \FFq^\times$ such that $S$ is equal to one of the following three groups:
    \begin{enumerate}
        \item $ Q(K)=\begin{pmatrix}K & 0_{1, n-1} \\ 0_{n-1, 1} & \GLno\end{pmatrix}$,
        \item $ P_V(K) = \begin{pmatrix}K & (\FFq^{n-1})^t \\ 0_{n-1, 1} & \GLno\end{pmatrix}$, or
        \item $ P_H(K) = \begin{pmatrix}K & 0_{1, n-1} \\ \FFq^{n-1} & \GLno\end{pmatrix} = P_V(K)^t $.
    \end{enumerate}
\end{proposition}

\begin{proof}
The case $n \ge 4$ suffices for our needs and follows directly from \cite[Theorem 2]{borevich1984arrangement}. 
The sharper case $n \ge 3$ can be deduced with some extra effort from the methods of \cite{McLaughlin1967}.
\end{proof}

\begin{corollary}
\label{cor:on subgroups of specific index}
 Suppose $n \geq 3$. Let $S$ be a proper subgroup of $\GLn$ containing the subgroup $\GLno$ as in Proposition \ref{prop:gl_n_1_overgroups}. If $q^n-1$ divides $[\GLn:S]$ then either
 \begin{enumerate}
    \item $S = P_V = P_V(\set{1})$, namely $S = \begin{pmatrix} 1 & (\FFq^{n-1})^t \\ 0_{n-1, 1} & \GLno \end{pmatrix}$,
    \item $S = P_H = P_H(\set{1}) = P_V^t$, or
    \item $S = P_V \cap P_H = Q(\set{1}) = \begin{pmatrix}1 & 0_{1, n-1} \\ 0_{n-1, 1} & \GLno\end{pmatrix}$.
\end{enumerate}
\end{corollary}
\begin{proof}
Let us look at the indices of the proper subsgroups of $\GLn$ which are overgroups of $\GLno$ listed in \Cref{prop:gl_n_1_overgroups}.  
These depend on the choice of a subgroup $K \le \FFq^\x$.
We see that
\[ [\GLn : Q(K)] = q^{n-1}(q^n-1) / \abs{K}\]
and
\[ [\GLn : P_V(K)] = [ \GLn : P_H(K) ] = (q^n-1)/\abs{K}.\]
Since $\abs{K}$ is coprime to $q$, the only possible subgroup $K$ that meets the divisibility constraint is $K = \set{1}$. This concludes the proof.
\end{proof}

\section{Permutational actions of \texorpdfstring{$\GLi$}{\GLipdf}}
\label{section:analysis_of_action}

Fix a prime power $q$. We introduce an abstract notion of permutational actions of the group $\GLi$. The main result of the section says that such an action must be equivariantly isomorphic to an action of the form $V^{k,\l}$. Given an arbitrary set $C$, a point $\gamma \in C^{\NN}$ in the product set and an index set $I =\{i_1<\cdots<i_n§\}\subset \NN$ we denote $\gamma_I = (\gamma_{i_k})_{k=1}^n \in C^n$.

\begin{definition}\label{def:permutational_action}
    Let $C$ be a finite set. Consider an arbitrary action of the group $\GLi$ on the product set $C^\NN$. The action is called \emph{local} if for every $n \in \NN$,   finite subset $I \subset \NN$ with $|I|=n$ and matrix $A \in \GLn$ we have
    \[\gamma_I =\delta_I \; \Rightarrow \; ([A]_I \gamma)_I = ([A]_I\delta)_I \quad \forall \gamma,\delta \in C^\NN \]
    and
    \[([A]_I \gamma)_{\NN \setminus I} = \gamma_{\NN\setminus I} \quad \forall \gamma \in C^\NN.\]
    The action is called \emph{permutational} if permutation matrices act in the standard manner. Namely, if $P_\sigma \in \GLi$ is a permutation matrix corresponding to the finitary permutation $\sigma \in S_\infty$ then 
    \[(P_\sigma \cdot \gamma)_i = \gamma_{\sigma(i)}\] 
    for every $\gamma \in C^\NN$.
\end{definition}

    We will refer to elements of $C$ as \emph{colors}, for reasons that will become clearer in the next sections. A local action as in Definition \ref{def:permutational_action} induces for every $n \in \NN$ an action of the group $\GLn$ on the product set $C^n$ of $n$-tuples of colors, by means of identifying this group with the subgroup $\GLn = \GLj{\{1,\ldots,n\}} \le \GLi$. 
    The main result of the current section is the following.



\begin{theorem}\label{thm:action_classification}
     Let $C$ be a set with $|C| = q^r$ where $r \in \NN$. Let $\red \in C$ be a distinguished color. Suppose $\GLi$ admits a local permutational action on the product set $C^\NN$. Assume that the action satisfies the properties
    \begin{enumerate}
        \item \label{item:red_fixed_assumption} the constant vector $\rho =  (\red)_{n \in \NN} \in C^\NN$ is a fixed point for the action, and
        \item \label{item:no_fixed_points_assumption} if $A \in \GLn $ is a matrix not admitting $1$ as an eigenvalue then the only fixed point of the induced action of $A$ on $C^n$ is the tuple $\rho_n = (\red)_{i=1}^n$.
    \end{enumerate}
    There exist $k,\l \in \NN_0$ with $k + \l = r$ and a bijection $\pi : C \to \FFq^k \x \FFq^\l$ such that the bijection $\pi^\NN : C^\NN \to V^{k, \l}$ is $\GLi$-equivariant. The bijection $\pi$ satisfies $\pi(\red) = (0,0) \in \FFq^k \x \FFq^\l$.
\end{theorem}

For the rest of this section, we will assume that the objects which appear in \Cref{thm:action_classification} are fixed.  In particular, we are given a local permutational action on color vectors satisfying the above two additional assumptions.

\begin{remark}
The material of the current \Cref{section:analysis_of_action} is relevant to the rest of the paper only through the application of \Cref{thm:action_classification}. A reader interested in the fastest route to understand Theorems \ref{thm:main_theorem} and \ref{thm:gl_inf_theorem} is encouraged to skip the rest of \Cref{section:analysis_of_action} upon first reading.
\end{remark}

We will be making the following identification implicitly: 
for any pair $n,m \in \NN \cup \{\infty\}$ where $n \leq m$ we identify the product set $C^n$ with the subset $\setcomp{\gamma \in C^m}{\gamma_i = \red \text{ for all } i > n}$ of the larger product set $C^m$.
Since the action of $\GLi$ is local, this subset is invariant under the  action of the subgroup $\GLn$, and this action coincides with the restriction of the action of $\GLm$ on $C^m$.

\begin{remark}
The largest matrix which plays a role in the proof of \Cref{thm:action_classification} has size $4 \x 4$. Hence the proof applies in the same way for actions of the finite group $\GLn$ on the set $C^n$ provided $n \geq 4$.
We have not made any attempt to optimize the bound $n \ge 4$, and it may work for actions of $\GLj{3}$ or $\GLj{2}$ too.
\end{remark}

\subsection{Color types}
The first step in the proof of \Cref{thm:action_classification} is to exhibit a product structure on the set of colors $C$, i.e. write $C = C_V \x C_H$ for some $k + \l = r$ and some sets $C_V,C_H$ of size $\abs{C_V} = q^k$ and $\abs{C_H} = q^\l$.

We will use the distinguished element $\red \in C $ to construct the two subsets $C_V \x \set{\red}$ and $\set{\red} \x C_H$.
The property that enables us to detect which colors belong to these two subset sets is the type of the stabilizer of a certain vector related to the color. 

Let $c \in C$ be any color other than the distinguished color $\red$. Let $n \ge 3$. We examine the $\GLn$-orbit of the vector 
\[v_c = (c, \red, \ldots, \red) \in C^n.\] 
Let $S_c = \Stab_{\GLn}(v_c)$ be its stabilizer.
The subgroup $S_c$ contains an isomorphic copy of $\GLno$ embedded as block-diagonal matrices in the last $n-1$ coordinates.
Let $s \in \GLn$ be any Singer cycle (see \Cref{subsection:singer_cycles}). Under Assumption (\ref{item:no_fixed_points_assumption}) of \Cref{thm:action_classification}, the cyclic group $\angles{s}$ acts freely on the $\GLn$-orbit of $v_c$. This means that the order $\abs{\angles{s}} = q^n-1$ divides the orbit's size, which is equal to $[\GLn : S_c]$. In particular $S_c$ is a proper subgroup of $\GLn$. 
We conclude by Corollary \ref{cor:on subgroups of specific index} that $S_c$ must be one of three subgroups  $P_V, P_H$  or $P_V \cap P_H$.

\begin{definition}
    A color $c \in C $ other than the distinguished color $\red$ will be called \textit{vertical} (resp. \textit{horizontal}) if 
    $S_c = \Stab_{\GLn}(v_c) = P_V$  (resp. $S_c = P_H$).
    We let $C_V$ (resp. $C_H$) be the set of vertical (resp. horizontal) colors, together with the $\red$ color.
    A color is said to be \textit{pure} if it is either vertical or horizontal.
\end{definition}

It is not hard to check that the above definition does not depend on the choice of $n$.
In what follows, we will show that $C \cong C_V \times C_H$, and that there are non-negative integers $k$ and $ \l$ with $k + \l = r$ such that $C_V^\NN \cong V^{k, 0}$ and $C_H^\NN \cong V^{0, \l}$. 
These isomorphisms justify the terminology ``vertical'' and ``horizontal''.


\begin{lemma}\label{lemma:elementary_action_basic}
    Let $c, d \in C$ be a pair of colors. Write
    \[
    E_{21} \cdot (c, d) = (c', d') \in C^2 \quad \text{and} \quad
    E_{12}^{-1} \cdot (c, d) = (c'', d'') \in C^2.
    \]
    We have
    \begin{enumerate}
        \item \label{item:elementary_action_basic_first} 
        If $d \in C_V$ then $d' \in C_V$.
        If in addition $c \in C_V$ then $c' \in C_V$.
        \item 
        If $d \in C_H$ then $d'' \in C_H$. If in addition $c \in C_H$ then $c'' \in C_H$.
    \end{enumerate}
\end{lemma}

\begin{proof}
    We only show item (1), as the proof of the second item (2) is analogous. Assume that $d \in C_V$. We need to show $d'  \in C_V$. That is, we need to check that every matrix $B \in P_V$ fixes the vector $v_{d'} = (d', \red, \red) \in C^3$. 
    Since the action is local, this can be done by showing that the larger matrix $\begin{psmallmatrix}1 & 0\\ 0 & B \end{psmallmatrix}$ fixes the vector $(c', d', \red, \red)$. 
    We may write $B = \begin{psmallmatrix} 1 & y \\ 0_{2, 1} & A \end{psmallmatrix}$ for some matrix  $A \in \GLj{2}$ and vector $y \in \Mj{1}{2}$. We calculate:
    \begin{align*}
        \begin{pmatrix}1 &  &  \\  & 1 & y \\ & 0_{2, 1} & A \end{pmatrix} \cdot (c', d', \red, \red)
        & = \begin{pmatrix}1 &  &  \\  & 1 & y \\ & 0_{2, 1} & A \end{pmatrix} \begin{pmatrix} E_{21} & 0_{2}\\ 0_2 & I_2\end{pmatrix} \cdot (c, d, \red, \red)
        \\ &= \begin{pmatrix} E_{21} & 0_{2}\\ 0_2 & I_2\end{pmatrix} 
        \begin{pmatrix}1 &  &  \\  & 1 & y \\ & 0_{2, 1} & A \end{pmatrix} \cdot (c, d, \red, \red)
        \\ &= \begin{pmatrix} E_{21} & 0_{2}\\ 0_2 & I_2\end{pmatrix} \cdot (c, d, \red, \red) = (c', d', \red, \red).
    \end{align*}
    We get $d' \in C_V$ as required.
    If in addition $c \in C_V$, we follow a similar strategy and calculation:
    \begin{align*}
        \begin{pmatrix}1 & & y \\ & 1 & \\ 0_{2, 1} & & A \end{pmatrix} \cdot (c', d', \red, \red)
        & = \begin{pmatrix}1 & & y \\ & 1 & \\ 0_{2, 1} & & A \end{pmatrix} \begin{pmatrix} E_{21} & 0_{2}\\ 0_2 & I_2\end{pmatrix} \cdot (c, d, \red, \red)
        \\ &= 
        \begin{pmatrix}E_{21} & 0_{2}\\ 0_2 & I_2\end{pmatrix} 
        \begin{pmatrix}1 & & y \\ & 1 & \\ 0_{2, 1} & & A \end{pmatrix} 
        \begin{pmatrix}1 & 0 & 0_{1,2} \\ 0 & 1 & -y \\ 0_{2,1} & 0_{2,1} & I_2 \end{pmatrix} \cdot (c, d, \red, \red)
        \\ &= (c', d', \red, \red).
    \end{align*}
    In the last equation, the second and third matrices act trivially because $c, d \in C_V$.
    Since the action is local, the above computation implies that the stabilizer of the vector $v_{c'} = (c', \red, \red)$ is $P_V$. Therefore $c' \in C_V$ as required.
\end{proof}

\subsection{Orbits of pure colors} \label{subsection:single_color_orbit}
The goal of this section is  classify the restriction of the $\GLn$-action to the orbit of a single vector of the form $v_c = (c, \red, \ldots, \red) \in C^n$ for a pure color  $c \in C$.

The action of the group $\GLi$ on the product space $C^\NN$ induces an action of its commutative subgroup $\GLj{1} \cong \FFq^\x$ on the set $C$. We write this action as $\alpha \cdot c$ for $\alpha \in \FFq^\x$ and $c \in C$.
By Assumption (\ref{item:red_fixed_assumption}) of   \Cref{thm:action_classification} and by fact the action is local, we have $\alpha \cdot \red = \red$ for all $\alpha  \in \FFq^\x$. We extend this by  setting $0 \cdot c = \red$ for all $c \in C$.

We now claim that the $\GLj{1}$-action on the set of colors $C$ preserves vertical and horizontal colors. To see this, consider the scalar matrix $\alpha I = D_1(\alpha) \cdots D_n(\alpha) \in \GLn$ with $\alpha \in \FFq^\x$. On the one hand, we have $\alpha I \cdot v_c = v_{\alpha \cdot c}$ for any and color $c \in C$, as each matrix $D_i(\alpha)$ acts locally on the $i$-th coordinate via  the above action of $\GLj{1}$. On the other hand, the matrix $\alpha I$ normalizes the three subgroups $P_V, P_H$ and $\GLno$. The claim follows. 
In addition, this action is free on $C \setminus \set{\red}$ by Assumption (\ref{item:no_fixed_points_assumption})  of   \Cref{thm:action_classification}. 
We later extend it to an $\FFq$-vector space structure on $C_V$ and on $C_H$.

\begin{proposition}\label{prop:single_orbit_action}
    Let $n \geq 3$. Let $c \in C$ be a pure color. Let
    \[
        O_c = \set{\rho} \union \GLn \cdot v_c \subseteq C^n
    \]
    be the orbit of $v_c$ in $C^n$ together  with the distinguished vector  $\rho = (\red,\ldots,\red)$.
    Let 
    \[
        T_c = \FFq \cdot c  \subseteq C
    \]
    be the set of multiples of $c$. Let $\theta_c : T_c \to \FFq$ be the bijection sending $\alpha \cdot c$ to $\alpha$.
    If $c$ is vertical (resp. horizontal) then 
   $T_c \subseteq C_V$ (resp. $T_c \subseteq C_H$) and
    $O_c$ equals the Cartesian power $T_c^n$ as a subset of $C^n$.
    In addition, 
    the bijection $\theta_c^n : O_c = T_c^n \to \FFq^n$ is $\GLn$-equivariant,  where the group acts on $\FFq^n$ via  the standard (resp. dual) action.
\end{proposition}

\begin{proof}
Let $c \in C_V$ be a vertical color.
The orbit $\GLn \cdot v_c$ is $\GLn$-equivariantly isomorphic to $\GLn / P_V \cong \FFq^n \setminus \set{0}$ with the standard $\GLn$-action. This isomorphism sends to vector $v_c$ to $e^1$.
Both points $\rho \in O_c$ and $0 \in \FFq^n$ are fixed by $\GLn$. We may extend the above and define a $\GLn$-equivariant isomorphism $\p_c : O_c \to \FFq^n$ with $\p_c(\rho) = 0$ and $\p_c(v_c) = e^1$.

The fact that $T_c \subseteq C_V$ (resp. $T_c \subseteq C_H$) has been established in the discussion preceding this proposition.
We will show that $O_c = T_c^n$ and that the map $\p_c$ coincides with $\theta_c^n$. We start with two claims.

\begin{claim*}
    Let $1 \leq i, j \leq n$. A point $\gamma \in O_c \subseteq C^n$ has $\gamma_i = \gamma_j$ if and only if $\p_c(\gamma)_i = \p_c(\gamma)_j$.
\end{claim*}
\begin{proof}
    The $\GLi$-action is permutational. This means that  $\gamma_i = \gamma_j$ if and only if $\sigma_{ij} \cdot \gamma = \gamma$.
    As the map $\p_c$ is equivariant, this happens if and only if $\sigma_{ij} \p_c(\gamma) = \p_c(\gamma)$. The latter is equivalent to $\p_c(\gamma)_i = \p_c(\gamma)_j$.
\end{proof}
 
\begin{claim*} 
    Any pure color $d \in C_V$ satisfies $
        E_{21} \cdot (d, \red) = (d, d)$.
\end{claim*}
\begin{proof}
    Let $d \in C_V$ be a pure color. If $d = \red$ then $E_{21} \cdot (\red, \red) = (\red, \red)$ by Assumption (\ref{item:red_fixed_assumption}) of  \Cref{thm:action_classification}, in which case the claim holds. We assume from now on $d \neq \red$.
    Recall that we identify $\GLj{2}$-equivariantly the set $C^2$ with the subset of vectors in $C^n$ whose $3,\ldots,n$-th entries are $\red$. In particular  $(d, \red)$ is identified with the vector $v_d$.
    Let $\p_d : O_d \to \FFq^n$ be the isomorphism constructed above with respect to the pure color $d$. So $\p_d\left((d, \red)\right) = \p_d(v_d) = e^1$.
    Since $E_{21} e^1 = e^1 + e^2$ we may write
    \[
        E_{21} \cdot (d, \red) = (d', d')
    \]
    for some color $d' \in C$. By \Cref{lemma:elementary_action_basic} we have $d' \in C_V$. By the same argument, there is some color $d'' \in C_V$ such that
    \[
        E_{21} \cdot (d', \red) = (d'', d'').
    \]
    Altogether, we get
    \begin{align}\label{eq:e21_e31_calc}
        E_{21}
        E_{31} \cdot
        (d, \red, \red)
        = E_{21} \cdot
        (d', \red, d')
        = (d'', d'', d').
    \end{align}
    On the other hand, applying $\p_d$ to this equation and using equivariance gives
    \[
        \p_d\left( E_{21} E_{31} \cdot (d, \red, \red) \right) = E_{21} E_{31} e^1 = (1, 1, 1)^t.
    \]
    Thus the colors appearing on the right-hand side of \cref{eq:e21_e31_calc} are equal to each other, namely $d' = d''$. This gives
    \[
        E_{21} \cdot (d, \red) = E_{21} \cdot (d', \red).
    \]
    We obtain $d = d'$ as required.
\end{proof}

    To ease notations, and having the color $c$ fixed, we will write from now onward $O=O_c$, $T = T_c$, $\p = \p_c$ and $\theta = \theta_c$. We proceed to show that $O = T^n$ and that $\p = \theta^n$. Notice that by size considerations, it is enough to show the containment $T^n \subseteq O$.

    Consider a point $\gamma \in T^n$. We wish to show that $\gamma \in O$ and $\p(\gamma) = (\theta(\gamma_i))_i$ for all $1 \le i \le n$. We will show this by induction on the number of non-$\red$ entries in $\gamma$. 
    The case where there are no such entries (i.e. $\gamma = \rho$) is trivial. If $\gamma$ contains a single non-$\red$ color $d \in T_c$, say at the index $1 \le i \le n$, then
    \[
    \gamma = \sigma_{1i} \cdot v_d = \sigma_{1i} D_1(\theta(d)) \cdot v_c \in O,
    \]
    and
    \[
    \p(\gamma) = \p {\left( \sigma_{1i} D_1 (\theta(d)) \cdot v_c \right)} = \sigma_{1i} D_1(\theta(d)) e^1 = (0, \ldots, \theta(d), \ldots, 0)^t
    \]
    as required.

    For the induction step, suppose two entries $\gamma_i$ and $\gamma_j$ are non-$\red$ for some $1 \le i\neq j\le n$. 
    Let $\alpha \in \FFq^\x$ be a scalar such that $\alpha \cdot \gamma_j = \gamma_i$. 
    Then the application $D_j(\alpha)\cdot \gamma$ has equal entries in its $i$-th and $j$-th indices. 
    By the second claim above, we see that the vector $E_{ij}^{-1} D_j(\alpha) \cdot \gamma$ has exactly one more $\red$ entry than the vector $\gamma$.
    By the induction hypothesis, we get $E_{ij}^{-1} D_j(\alpha) \cdot \gamma \in O$, and so also $\gamma \in O$.
    The induction hypothesis gives:
    \begin{align*}
        \p(\gamma) &= D_{j}(\alpha^{-1}) E_{ij} \p(E_{ij}^{-1} D_j(\alpha) \cdot \gamma) \\ &= 
        D_{j}(\alpha^{-1}) E_{ij} (\theta(\gamma_1), \ldots, 0, \ldots, \theta(\gamma_i), \ldots, \theta(\gamma_n))^t
         = (\theta(\gamma_k))_k,
    \end{align*}
    as required.

    Lastly, the case where $c \in C_H$ is entirely analogous. This time $\GLn / P_H \cong \FFq^n \setminus \set{0}$, and the group $\GLn$ is acting on the latter space via the dual action.
\end{proof}

\subsection{Vector space structure on pure colors}
We define an addition operator on the set $C_V$ of vertical colors and on the set $C_H$ of horizontal colors.
For each pair of vertical colors $c, d \in C_V$ we set $c+d = d'$ where
\[
E_{21} \cdot (c, d) = (c', d') \in C^2.
\]
For a pair of horizontal colors $c, d \in C_H$ we set $c + d = d'$ where
\[
E_{12}^{-1} \cdot (c, d) = (c', d') \in C^2.
\]
Notice that by \Cref{lemma:elementary_action_basic}, in both cases the result $c+d$ has the same color type as $c$ and $d$.

\begin{proposition} \label{prop:pure_color_action}
    The operations $ + : C_V \times C_V \to C_V$  and $\FFq \times C_V \to C_V$ endow $C_V$ with an $\FFq$-vector space structure. 
    The $\red$ color is the zero vector.
    The action of $\GLn$ on $C_V^n$ corresponds to the standard action. The same result holds for $C_H$, and the action of $\GLn$ on $C_H^n$ corresponds to the dual action.
\end{proposition}

\begin{proof}
    We verify the vector space axioms over $\FFq$ on the set $C_V$. Let $c, d \in C_V$ be a pair of pure colors. These colors satisfy
    \[
        E_{31} E_{32} \cdot (c, d, \red)
        = E_{31} \cdot (c, d, d) = (c', d, c+d)
    \]
    for some $c \in C_v$
    In the first equality we have used the second claim from the proof of \Cref{prop:single_orbit_action}
    for some $c'_1,c'_2 \in C_V$.
    Since the matrices $E_{32}$ and $E_{31}$ commute, applying $E_{32} E_{31}$ to the same vector gives $d + c$ in the third entry. It follows that addition is commutative.

    We check associativity of addition. Let $c_1, c_2, c_3 \in C_V$ be three pure colors. They satisfy
    \[
        E_{41} E_{42} E_{43} \cdot (c_1, c_2, c_3, \red) =
        (c_1', c_2', c_3, c_1 + (c_2 + c_3)).
    \]
    Since all these elementary matrices commute, the application of $E_{43} E_{41} E_{42}$ to $(c_1, c_2, c_3, \red)$ gives the same result. It follows that 
    \[
        c_1 + (c_2 + c_3) = c_3 + (c_1 + c_2) = (c_1 + c_2) + c_3
    \]
    as required.
    Further, 
    by \Cref{prop:single_orbit_action}, the distinguished color $\red \in C_V$ is the additive identity and $(-1) \cdot c$ is the additive inverse. This shows that $(C_V, +)$ is an abelian group. 

    We consider the axioms of the scalar multiplication. Let $\alpha, \beta \in \FFq$ be a pair of scalars. The facts that $\alpha\cdot(\beta \cdot c)=(\alpha\beta)\cdot c$ and that $1 \cdot c = c$ follows from the construction of the action. Next, note that \Cref{prop:single_orbit_action} gives
    \[\begin{pmatrix}\alpha\\ \alpha+\beta\end{pmatrix} = \begin{pmatrix} \alpha & 0 \\ \alpha+\beta & \beta \end{pmatrix} e_1 =\theta_c(E_{21}D_1(\alpha)D_2(\beta) E_{21} \cdot v_c) = \theta_c((c',\alpha \cdot c + \beta \cdot c)).\]
    This implies the axiom $(\alpha + \beta) \cdot c =\alpha \cdot c + \beta \cdot c$. Similarly, the fact that the diagonal matrix $\alpha I$ commutes with $E_{21}$ gives
    \[ (\alpha \cdot c'', \alpha \cdot (c+d))=(\alpha I_2)E_{21} \cdot (c,d)=E_{21} (\alpha I_2) \cdot(c,d) = E_{21}\cdot(\alpha \cdot c, \alpha \cdot d) = (c',\alpha \cdot c+ \alpha \cdot d).\]
    The distributivity axiom $\alpha \cdot (c+d)= \alpha \cdot c + \alpha \cdot d$ follows.     
    Hence $C_V$ is a vector space over $\FFq$.

    It remains to show that the subset $C_V^n$ is $\GLn$-invariant, and that the $\GLn$-action on $C_V^n$ corresponds to the  standard action on $\FFq^n$.
    Both facts are true for permutation matrices by the assumption that the action is permutational. They are also true for diagonal matrices of the form $D_i(\alpha)$ by definition.
    The subset $C_V^n$ is invariant under elementary matrices by \Cref{lemma:elementary_action_basic}.
    Since the group $\GLn$ is generated by the $E_{ij}$'s and the $D_i(\alpha)$'s for $1 \leq i, j \leq n$ and $\alpha \in \FFq^\x$, we conclude that $C_V^n$ is indeed $\GLn$-invariant.
    
    It remains to show that the action of the elementary matrices $E_{ij}$ on $C_V^n$ corresponds to the standard action. 
    Since all elementary matrices are conjugate by permutation matrices, it suffices to show this for a single elementary matrix $E_{21}$.
    Recall that $c, d \in C_V$ is a pair of pure colors. Let $c' \in C$ be a color such that
    \[
        E_{21} \cdot (c, d) = (c', c+d).
    \]
    We need to show that $c' = c$.
    By \Cref{lemma:elementary_action_basic} we know that $c' \in C_V$. Now calculate:
    \[
        E_{21}^{-1} E_{31} E_{21} \cdot (c, d, \red)
        = E_{21}^{-1} E_{31} \cdot (c', c+d, \red)
        = E_{21}^{-1} \cdot (c', c+d, c')
        = (c, d, c').
    \]
    The two matrices $E_{21}$ and $E_{31}$ commute. Hence the above is  equal to
    \[
        E_{31} \cdot (c, d, \red) = (c, d, c).
    \]
    We conclude that $c = c'$, as required.
    This finishes the case of vertical pure colors.

    The case of the horizontal pure  colors $C_H$ is entirely analogous, except that the action on $C_H^n$ corresponds to the dual action of $\GLn$.
\end{proof}

\subsection{Colors as products of pure colors}
In this section we prove \Cref{thm:action_classification}.
By \Cref{prop:pure_color_action} the sets $C_V$ and $C_H$ of pure colors are vector spaces over $\FFq$, and $\GLn$ acts on $C_V^n$ and on $C_H^n$ by the standard and the dual action respectively. It  remains is to identify the set of all colors $C$ with the direct product $C_V \x C_H$.
The following lemma is a first step in this direction. It adds more information to \Cref{lemma:elementary_action_basic}.

\begin{lemma}\label{lemma:elementary_standard_red}
    Let $c \in C$ be any color. Then
    \[
    E_{21} \cdot (c, \red) = (c, d) \quad \text{and} \quad E_{12}^{-1} \cdot (c, \red) = (c, d')
    \]
    for some pure colors $d \in C_V$ and $d' \in C_H$.
\end{lemma}

\begin{proof}
    We prove only the first assertion for $E_{21}$ and $C_V$, the second assertion for $E_{12}^{-1}$ and $C_H$ being analogous (and its proof is obtained by applying transpose inverse to all the matrices involved).
    
    Write $E_{21} \cdot (c, \red) = (c_1, d)$ for some $c_1 \in C$ color. 
    By \Cref{lemma:elementary_action_basic} we know $d \in C_V$. Our goal is to show $c = c_1$. 
    Define the colors $c_2, c_3, d_2, d_3 \in C$
    according to the calculation
    \begin{equation} \label{eq:e_21_31_conjugation}
        E_{21}^{-1} E_{31} E_{21} \cdot (c, \red, \red)
        = E_{21}^{-1} E_{31} \cdot (c_1, d, \red)
        = E_{21}^{-1} \cdot (c_2, d, d_2) = (c_3, d_3, d_2).
    \end{equation}
    The matrices $E_{21}$ and $E_{31}$ commute, so that $E_{21}^{-1} E_{31} E_{21} = E_{31}$. We get
    \[
    (c_3, d_3, d_2) = (c_1, \red, d).
    \]
    Hence, the right-most  equality in \cref{eq:e_21_31_conjugation} implies that
    \[
    E_{21} \cdot (c_1, \red) = (c_2, d).
    \]
    Next, we have the following, where the right-most  equality uses the second claim from the proof of \Cref{prop:single_orbit_action} applied for $d \in C_V$:
    \[
    E_{32}^{-1} E_{31} E_{21} \cdot (c, \red, \red) =
    E_{32}^{-1} E_{31} \cdot (c_1, d, \red)
    = E_{32}^{-1} \cdot (c_2, d, d)
    = (c_2, d, \red).
    \]
    The standard commutator identity $[E_{32}, E_{21}] = E_{31}$ can be rewritten as
    \[
    E_{32}^{-1} E_{31} E_{21} = E_{21} E_{32}^{-1}.
    \]
    Acting with the matrix on right-hand side of the above on the vector $(c, \red, \red)$ gives $(c_1, d, \red)$.
    Hence $c_1 = c_2$. Further
    \[
        E_{21} \cdot (c, \red) = (c_1, d) = (c_2, d) = E_{21} \cdot (c_1, \red)
    \]
    so that $c = c_1$, as required.
\end{proof}

Note that \Cref{lemma:elementary_standard_red} shows in particular that pure colors do exist.

\begin{proof}[Proof of \Cref{thm:action_classification}]

We define a map 
\[\pi : C \to C_V \times C_H\]
as follows. For each color $c \in C$ set
\[
E_{21} \cdot (c, \red) = (c, \pi_V(c)) \quad \text{and} \quad 
E_{12}^{-1} \cdot (c, \red) = (c, \pi_H(c)).
\]
Note that the first entries in the above pairs are equal to  $c$ by \Cref{lemma:elementary_standard_red}. The desired map $\pi$ is defined by
\[\pi(c) = (\pi_V(c), \pi_H(c))\]
for all $c \in C$.
We claim that $\pi$ is a bijection, and that its product $\pi^n : C^n \to C_V^n \x C_H^n$ is $\GLn$-equivariant. 
    Recall that we know by \Cref{prop:pure_color_action} that $\GLn$ acts via the standard and dual action on $C_V^n$ and on $C_H^n$, respectively. Hence we will be able to deduce that $\pi^\NN : C^\NN \to C_V^\NN \x C_H^\NN \cong V^{k, \l}$ is a $\GLi$-equivariant isomorphism, as required.

    Let us first show that the product map $\pi^n$ is $\GLn$-equivariant.
    Consider the matrix $\Pi_V^n = \begin{psmallmatrix} I_n & 0_n \\ I_n & I_n \end{psmallmatrix} \in \GLj{2n}$. Denote as usual $\rho = (\red,\ldots,\red) \in C^n$ so that $\rho$ corresponds to the zero vector. Let $\gamma \in C^n$ be arbitrary. Notice that we have
    \[
        \Pi_V^n \cdot (\gamma, \rho) = (\gamma, \pi_V^n(\gamma)).
    \]
    Let $A \in \GLn$ be any matrix. We have on the one hand
    \[
        (\gamma, A \pi_V^n(\gamma))
        = \begin{pmatrix} I_n & 0_n \\ 0_n & A \end{pmatrix} \Pi_V^n \cdot (\gamma, \rho),
    \]
    while on the other hand
    \[
        (\gamma, \pi_V^n(A \cdot \gamma))
        = \begin{pmatrix} A^{-1} & 0_n \\ 0_n & I_n \end{pmatrix} \Pi_V^n \begin{pmatrix} A & 0_n \\ 0_n & I_n \end{pmatrix} \cdot (\gamma, \rho).
    \]
    Denote by $M_1$ and $M_2$  the two matrices being applied to the vector $(\gamma, \rho)$ in the first and second equation above, respectively.
    $M_1$ and $M_2$ agree on $(\gamma, \rho)$ if and only if $M_2^{-1} M_1$ fixes it. And indeed the matrix $M_2^{-1} M_1 = \begin{psmallmatrix} I_n & 0_n \\ 0_n & A\end{psmallmatrix}$  fixes $(\gamma, \rho)$ by Assumption (\ref{item:red_fixed_assumption}). We conclude that $A \pi_V^n(\gamma) = \pi_V^n(A \cdot \gamma)$, as required.
    An analogous argument applies to $\pi_H^n$.

    Next, we wish to define a map $\phi : C_V \x C_H \to C$ which will be a candidate for the inverse of $\pi$. To do so, consider the matrices
    \[
        \Phi_1 = E_{23}^{-1} E_{31} 
        \quad 
        \text{and}
        \quad
        \Phi_2 = E_{31} E_{23}^{-1}
        .
    \]
    Given a vertical $c_V \in C_V$ and a horizontal $ c_H \in C_H$  pure color, we write
    \[
        \Phi_1 \cdot (c_V, c_H, \red) = (c_V, c_H', \phi_1(c_V, c_H)),
        \quad \Phi_2 \cdot (c_V, c_H, \red) = (c_V', c_H, \phi_2(c_V, c_H)),
    \]
    for some pair of maps $\phi_1, \phi_2 : C_V \x C_H \to C$ and colors $c_V', c_H' \in C$.
    In this notation, we have used the fact that $\Phi_1$ fixes the first entry, and $\Phi_2$ the second.
    Indeed, by \Cref{lemma:elementary_standard_red} the factor $E_{31}$ fixes the first entry and the factor $E_{23}^{-1}$ fixes the second, and the remaining factor is supported away from that fixed entry.

    We claim that $\phi_1 = \phi_2$. 
    This is because
    \[
    (c_V, c_H', \phi_1(c_V, c_H)) = \Phi_1 \cdot (c_V, c_H, \red) = \Phi_1 \Phi_2^{-1} \Phi_2 \cdot (c_V, c_H, \red) = E_{21}^{-1} \cdot (c_V', c_H, \phi_2(c_V, c_H))
    \]
    and the support of $E_{21}^{-1}$ is contained in the first two indices. This function will be denoted by $\phi$.

    Now we show that $\phi = \pi^{-1}$.
    Let $c_V \in C_V, c_H \in C_H$, and $c = \phi(c_V, c_H) \in C$. We want to show that $\pi_V(c) = c_V$ and $\pi_H(c) = c_H$.
    We apply $E_{43}$ composed with $\Phi_1$ to $(c_V, c_H, \red, \red)$ to get
    \[
        E_{43} \Phi_1 \cdot (c_V, c_H, \red, \red) =
        (c_V, c_H', c, \pi_V(c)).
    \]
    We argue that $c_V = \pi_V(c)$. To see this, we show that the permutation $\sigma_{14}$ stabilizes the above vector.
    We calculate:
    \[
        (E_{43} \Phi_1)^{-1} \sigma_{14} E_{43} \Phi_1
        = \begin{pmatrix}
            1 & 0 & 1 & 1 \\ 0 & 1 & 0 & 0 \\ 0 & 0 & 0 & -1 \\ 0 & 0 & -1 & 0
        \end{pmatrix}.
    \]
    This matrix fixes the second coordinate and acts on the remaining ones, which all lie in $C_V$, by the standard action. Hence it stabilizes $(c_V, c_H, \red, \red)$, as required.
    An analogous calculation shows $\pi_H(c) = c_H$. And so $\pi \circ \phi = \id$.

    To show $\phi \circ \pi = \id$, consider a color $c \in C$. The matrix
    \[
        \Pi = E_{13}^{-1} E_{21} = \begin{pmatrix}
            1 & 0 & -1 \\ 1 & 1 & 0 \\ 0 & 0 & 1
        \end{pmatrix}
    \]
    has the effect of applying both maps $\pi_V$ and $\pi_H$, in the sense that
    \[
    \Pi \cdot (c, \red, \red) =
    (c, \pi_V(c), \pi_H(c)).
    \]
    Take the vector $v_c = (c, \red, \red, \red)$. Apply $\Pi$ to its first three coordinates, and then $\Phi_1$ to its last three coordinates. The resulting vector's first entry is $c$ and its last entry is $\phi(\pi(c))$. To show that these two are equal, it suffices to show that the permutation matrix $\sigma_{14}$ stabilizes this vector. This amounts to showing that the matrix
    \[
    \left(\begin{pmatrix} 1 & 0_3 \\ 0_3 & \Phi_1 \end{pmatrix} \begin{pmatrix} \Pi & 0_3 \\ 0_3 & 1 \end{pmatrix}\right)^{-1} \sigma_{14} \left(\begin{pmatrix} 1 & 0_3 \\ 0_3 & \Phi_1 \end{pmatrix} \begin{pmatrix} \Pi & 0_3 \\ 0_3 & 1 \end{pmatrix}\right) = \begin{pmatrix}
        1 & 0 & 0 & 0 \\ 0 & 1 & 0 & 0 \\ 0 & -1 & 0 & -1 \\ 0 & -1 & -1 & 0
    \end{pmatrix}
    \]
    stabilizes $v_c$.
    The latter fact holds, since this matrix is of the form $\begin{psmallmatrix} 1 & 0_3 \\ 0_3 & B \end{psmallmatrix}$, and  $B$ acts only on $\red$ entries. We conclude that $\phi \circ \pi = \id$ and that  $\pi$ is a bijection.

    By \Cref{prop:pure_color_action} the sets $C_V$ and $C_H$ are vector spaces over $\FFq$. So $\abs{C_V} = q^k$ and $\abs{C_H} = q^\l$ for some $k, \l \geq 0$.
    As $\pi$ is a bijection, we must have $q^{k + \l} = \abs{C} = q^r$, so that $k + \l = r$ as required. 
    Finally, the map $\pi$ takes $\rho$ to $(\red, \red) \in C_V \x C_H$, which is the additive identity of the vector space.
\end{proof}

\section{Invariant random subgroups and traces}
\label{section:notations_and_preliminaries}
In this section we recall the notions of invariant random subgroups and Thoma traces. In addition, we discuss the classification of these objects for the finitary alternating group $\Ai$ and analyze its Vershik trace map. This section is a summary of known results.

\subsection{Invariant random subgroups}

Let $G$ be a countable discrete group. The \textit{Chabauty space} of the group $G$ is the set of its subgroups. It is denoted by 
\[
\Sub{G} = \setcomp{H}{H \leq G}.
\]
The space $\Sub{G}$ is naturally a subset of the power set $\set{0, 1}^G$. The latter is a metrizable compact space with the product topology.  Being a closed subset, $\Sub{G}$ inherits a metrizable compact topology from $\set{0, 1}^G$. It is called the \textit{Chabauty topology}.
The group $G$ acts on $\Sub{G}$ continuously by conjugation, namely
\[
g \cdot H = g H g^{-1} \quad \forall H \in \Sub{G}, \; g \in G.
\]

Let $\cM(\Sub{G})$ be the convex set of Borel probability measures on $\Sub{G}$ endowed with the weak-$*$ topology.
The conjugation action induces a dual action of $G$ on $\cM(\Sub{G})$ by affine transformations. 
An \textit{invariant random subgroup} is a $G$-invariant measure in $\cM(\Sub{G})$. The set of all such measures will be denoted by $\IRS(G)$. 

The set $\IRS(G)$ forms a convex and weak-$*$ compact subspace of $\cM(\Sub{G})$. An invariant random subgroup is called \textit{ergodic} if it is an extreme point of the space $\IRS(G)$. 
This happens if and only if every measurable $G$-invariant subset of $\Sub{G}$ is either null or conull. 
The invariant random subgroup $\delta_\triv \in \IRS(G)$ is called  \textit{trivial}.

Assume that the countable group $G$  acts on some Borel probability space $(X, \nu)$ by measure-preserving maps.
The stabilizer map $\Stab_G : X \to \Sub{G}$ is Borel and the pushforward of $\nu$ under this map gives an invariant random subgroup of $G$.
We note that every invariant random subgroup of $G$ is a stabilizer distribution of some such action \cite{AGV2014}.

\begin{definition}
Let $H$ be a subgroup of the countable discrete map $G$.    
The \textit{intersection map} is given by
\[
\Sub{G} \to \Sub{H}, \quad K \mapsto K \inter H \quad \forall K \in \Sub{G}.
\]
The intersection map is $H$-equivariant and continuous. Hence it induces  a pushforward map 
$
\IRS(G) \to \IRS(H)
$, which we also call an intersection map by a slight abuse of notations. 
\end{definition}

See e.g. \cite{AGV2014, gelander2015lectureinvariantrandomsubgroups, levit_lubotzky2022_metabelian} for  more details about invariant random subgroups.
 
\subsection{Thoma and Vershik traces}
\label{subsection:thoma and vershik traces}
We recall the notion of a trace in the sense of Thoma \cite{thoma1964unzerlegbaren,thoma1964_choquet_simplex}. 
These are generalizations of the normalized trace functions of finite-dimensional unitary representations of a group.
A very good reference for the theory of Thoma traces is \cite{bekka_la_harpe2020}.

\begin{definition}\label{def:thoma_trace}
    Let $G$ be a countable discrete group. A function $\chi : G \to \CC$ is called a \textit{trace of $G$} (in the sense of Thoma) if it satisfies:
    \begin{enumerate}
        \item $\chi$ is normalized, i.e.\ $\chi(e) = 1$, 
        \item $\chi$ is central, i.e.\ $\chi(g^{-1} h g) = \chi(h)$ for all $g, h \in G$, and
        \item $\chi$ is positive definite, i.e. $$\sum_{i, j = 1}^n \overline{a}_j a_i \chi(g_j^{-1} g_i) \geq 0$$ for all $n \in \NN$ and for all $a_1, \ldots, a_n \in \CC$ and $g_1, \ldots, g_n \in G$.
    \end{enumerate}
    The set of traces $\Tr(G)$ of the group $G$ forms a convex compact subset of $\ell^\infty(G)$. The set of extreme points of this set is called the \textit{Thoma dual of $G$} and is denoted $\Ch(G)$. The elements of $\Ch(G)$ will be called \textit{characters of $G$}.
\end{definition}
We mention that the terminology in Definition \ref{def:thoma_trace} is not completely standardized, and varies across the literature.
There is a natural map $\IRS(G) \to \Tr(G)$, as follows:

\begin{definition}\label{def:vershik_trace}
    Let $G$ be a countable group and  $\mu$ be an invariant random subgroup of $G$. The map 
    \[
        \chi_\mu : G \to [0, 1], \quad \chi_\mu(g) = \mu \left( \setcomp{H \in \Sub{G}}{g \in H} \right) 
    \]
    is a trace of the group $G$. It is called the \textit{Vershik trace} associated to $\mu$. The map $\mu \mapsto \chi_\mu$ is called the \textit{Vershik trace map}.
\end{definition}

To see that the Vershik trace $\chi_\mu$ associated to $\mu$ is positive definite, note that it can be expressed as  $\chi_\mu = \int_{\Sub{G}} \one_H \diff \mu(H)$. A straightforward computation shows that each indicator function $\one_H$ is positive definite.

We recall that a metrizable compact convex subset $X$ of a locally convex topological vector space is called a \textit{Choquet simplex} if for every point $x_0 \in X$ there is a unique Borel probability measure $\lambda$ supported on the extreme points of $X$ such that $x_0 = \int_X x \diff \lambda(x)$. 
Both $\Tr(G)$ and $\IRS(G)$ are Choquet simplices. See \cite{thoma1964_choquet_simplex} and \cite[Chapter 12]{phelps2002lectures} respectively.
Given a trace or an invariant random subgroup of $G$, we call the measure provided by Choquet's theorem its \textit{ergodic decomposition}.

\subsection{Thoma characters and invariant random subgroups of the finitary alternating group}
\label{subsection:irs_a_infty_description}
Let $\Ai$ be the infinite alternating group of finitary even permutations.
We review the theory of its Thoma characters, its invariant random subgroups and their relationship via the Vershik map.

As mentioned, the Thoma characters of $\Ai$ were classified by Thoma in \cite{thoma1964unzerlegbaren}.
The classification of the invariant random subgroups of the group $S_\infty$ of all finitary permutations was done by Vershik in \cite{Vershik2012} with a slight inaccuracy. In \cite{thomas2022_sinfty_irs_and_characters} Thomas presents a corrected argument for this classification.
To get the argument for $\Ai$ one can follow the proof of the classification of $\IRS(S_\infty)$ and check that it works the same, with appropriate modifications, for this subgroup. 
Alternatively, \cite{levit_lubotzky2022} presents a derivation of the invariant random subgroups of $\Ai$ from those of $S_\infty$.
The calculation of the Vershik traces and the description of these as a weighted average on the set of Thoma characters is done in \cite{thomas2022_sinfty_irs_and_characters} for the group $S_\infty$ of all finitary permutations. 
These calculations still hold with minor modifications for the case of $\Ai$, as we now explain.

\subsubsection{Thoma characters of $\Ai$}
Let $(\beta_i)_{i \in \NN}$ and $(\gamma_i)_{i \in \NN}$ be two non-increasing sequences of non-negative real numbers such that $$\sum_{i=1}^\infty \beta_i + \sum_{i=1}^\infty \gamma_i \leq 1.$$ 
Define the function $\chi_{(\beta; \gamma)} : \Ai \to \CC$ by
\[
    \chi_{(\beta; \gamma)}(\sigma) = \prod_{k=2}^\infty \left( \sum_{i=1}^\infty \beta_i^k + (-1)^{k+1} \sum_{i=1}^\infty \gamma_i^k \right)^{c_k(\sigma)} \quad \forall \sigma\in \Ai
\]
where for each $k \ge 2$ we let $c_k(\sigma)$ denote the number of cycles of length $k$ in the cycle decomposition of the permutation $\sigma \in \Ai$.
\begin{theorem}[{\cite[Satz 6]{thoma1964unzerlegbaren}}]
    The functions $\chi_{(\beta; \gamma)}$ described above are precisely the characters of $\Ai$.
    The identity $\chi_{(\beta; \gamma)}=\chi_{(\beta'; \gamma')}$ holds if and only if $(\beta,\gamma) = (\beta',\gamma')$ or $(\beta,\gamma) = (\gamma',\beta')$.
\end{theorem}

\subsubsection{Invariant random subgroups of $\Ai$} 
\label{subsub:IRS of Ai}
We describe the ergodic invariant random subgroups of the group $\Ai$ using stabilizers of ergodic actions.
Let $\cD \subseteq [0, 1]^{\NN_0}$  be the set of sequences of non-negative real numbers $\alpha=(\alpha_i)_{i\in \NN_0}$ which satisfy:
\begin{itemize}
    \item The sequence $(\alpha_{i})_{i=1}^\infty$ is non-increasing, and
    \item $\sum_{i=0}^\infty \alpha_i = 1$.
\end{itemize}

Fix a sequence $\alpha \in \cD$. Let  $I_\alpha = \setcomp{i \in \NN}{\alpha_i > 0}$ and consider the group $S_\alpha = \bigoplus_{i \in I_\alpha} \set{\pm 1}$.
Let $s_s$ be the multiplication homomorphism
\[
s_\alpha : S_\alpha \to \set{\pm 1}, \quad s_\alpha(t) = \prod_{i \in I_\alpha} t_i \quad \forall t \in S_\alpha.
\]
Denote $E_\alpha = \ker s_\alpha \le S_\alpha$.
Let $T$ be any fixed subgroup of $E_\alpha$. We will define the invariant random subgroup $\nu_\alpha^T$ of $\Ai$ corresponding to the data consisting of the pair $\alpha \in \cD$ and $T \in \Sub{E_\alpha}$.

The fixed sequence $\alpha \in \cD$ defines a probability measure $p_\alpha$ on $\NN_0$ by $p_\alpha(i) = \alpha_i$.
Let $\mu_\alpha$ be the corresponding product measure on $(\NN_0)^\NN$.
We think of elements $c \in (\NN_0)^\NN$ as colorings of the set of natural numbers $\NN$, where the probability of coloring a number $n \in \NN$ using the color $i \in \NN_0$ is $\alpha_i$, in an i.i.d.\ manner over all $n \in \NN$.
Given a coloring $c \in (\NN_0)^\NN$ and a color $i \in \NN_0$, consider the subset 
$$B_i(c) = \setcomp{n \in \NN}{c_n = i}.$$ 
Let $H_{\alpha,T}(c) \in \Sub{\Ai}$ be the subgroup consisting of all permutations $\sigma \in \Ai$ which satisfy:
\begin{itemize}
    \item $\sigma(B_i(c)) \subseteq B_i(c)$ for all $i \in \NN$,
    \item $\sigma |_{B_0(c)} = \id$, and
    \item $\left( \sgn (\sigma|_{B_i(c)})\right)_{i \in {I_\alpha}} \in T$.
\end{itemize}
In other words, the subgroup $H_{\alpha,T}(c)$ consists of all permutations that preserve  each color setwise, fix  the distinguished color $0 \in \NN_0$ pointwise, and whose sign pattern on the colors belongs to $T$.
The pushforward of the probability measure $\mu_{\alpha}$ under the $\Ai$-equivariant Borel map $H_{\alpha,T} : (\NN_0)^\NN \to \Sub{\Ai}$ yields an invariant random subgroup which we denote $\nu_{\alpha}^T$.
The following result says that these are, in fact, all the ergodic invariant random subgroups of $\Ai$.
\begin{theorem}[Vershik, {\cite[Theorem 9.1]{thomas2018_ainfty_irs}}]
For any ergodic invariant random subgroup $\nu$ of $\Ai$ there exist a sequence $\alpha \in \cD$ and a subgroup $T \in \Sub{E_\alpha}$  such that $\nu = \nu_\alpha^T$.
\end{theorem}

There are pairs $(\alpha, T)$ and $(\alpha, T')$ with $\nu_\alpha^T = \nu_\alpha^{T'}$ and $T \neq T'$. 
However, an ergodic invariant random subgroup determines the sequence $\alpha$ uniquely \cite[Remark 1.2]{thomas2022_sinfty_irs_and_characters}.

\subsubsection{The Vershik trace map for $\Ai$}
We describe the Vershik traces associated to the ergodic invariant random subgroups of $\Ai$.
As a consequence, it will follow  that the Vershik trace $\chi_{\nu_\alpha^T}$ is a character if and only if $T = E_\alpha$. 


Fix a sequence $\alpha \in \cD$. The group $S_\alpha$ is a discrete $2$-group,  so its characters are the homomorphisms $S_\alpha \to \set{\pm 1}$. These form the compact Pontryagin dual group $\widehat{S}_\alpha$. There is a natural identification
\[
    \widehat{S}_\alpha \cong \prod_{i \in I_\alpha} \set{0, 1}
\]
by means of the pairing 
\[
    \bigoplus_{i \in I_\alpha} \set{\pm 1} \x \prod_{i \in I_\alpha} \set{0, 1} \to \set{\pm 1}, \quad (t, u) \mapsto \prod_{i \in I_\alpha} t_i^{u_i} \quad \forall t \in S_\alpha, u\in \widehat{S}_\alpha.
\]
Under this identification, the multiplication homomorphism $s_\alpha$ corresponds to the vector $s_\alpha = (1)_{i \in I_\alpha} \in \widehat{S}_\alpha$.
The Pontryagin dual group of the subgroup $E_\alpha$ is naturally isomorphic to the quotient $\widehat{E}_\alpha \cong \widehat{S}_\alpha / \angles{s_\alpha}$.

Given the fixed sequence $\alpha \in \cD$ and for each  character $u \in \widehat{E}_\alpha$, we define the Thoma character $\chi_\alpha^u $ of the group $\Ai$ by
\[
    \chi_\alpha^u (\sigma) = \prod_{k=2}^\infty \left( \sum_{i \in I_\alpha} (-1)^{(k+1) u_i} \alpha_i^k \right) ^{c_k(\sigma)} \quad \forall \sigma \in \Ai.
\]
The character $\chi_\alpha^u $ coincides with the character $\chi_{(\beta; \gamma)}$ defined above, where 
\[
\beta_i = \begin{cases}
    \alpha_i, & \text{if } u_i = 0 \\
    0, & \text{if } u_i = 1
\end{cases}
\quad\text{and}\quad
\gamma_i = \begin{cases}
    0, & \text{if } u_i = 0 \\
    \alpha_i, & \text{if } u_i = 1.
\end{cases}
\]
We note that this formula is independent of the choice of representative for $u$ in $\widehat{S}_\alpha$, since if $\chi_{\alpha}^u = \chi_{(\beta; \gamma)}$ then $\chi_{\alpha}^{u + s} = \chi_{(\gamma; \beta)}$.
Thomas has calculated the Vershik trace of the invariant random subgroup $\nu_{\alpha}^T$:

\begin{theorem}[{\cite[Theorem 1.7]{thomas2022_sinfty_irs_and_characters}}] \label{thm:ainfty_irs_vershik_character}
Let $\alpha \in \cD$ be a  sequence and $T \le E_\alpha$ be a subgroup. 
Identify the Pontryagin dual group $\widehat{E_\alpha / T}$ with the subgroup of characters of $E_\alpha$ which are trivial on $T$. The Vershik trace of the invariant random subgroup $\nu_\alpha^T$ is given by
    \[
        \chi_{\nu_\alpha^T}  = \int_{\widehat{E_\alpha / T}} \chi_\alpha^u \; \diff u
    \]
    where the measure $\diff u$ is the  Haar probability measure on the dual group $\widehat{E_\alpha/T}$.
\end{theorem}

The proof in \cite{thomas2022_sinfty_irs_and_characters} is for  invariant random subgroups of $S_\infty$, and it works the same in this case by replacing the occurrences of $S_\alpha$ with $E_\alpha$.

\subsection{The orbit equivalence map}\label{subsubsection:the_orbit_map}
We take this opportunity to introduce several objects that involve the finitary alternating group and that will be useful in \Cref{subsection:color_definition} below:

Let $\cE(\NN) \subseteq 2^{\NN \x \NN}$ be the set of equivalence relations on $\NN$, with the standard Borel structure. The action of the finitary alternating group $\Ai$ on the set $\NN$ induces an action $\Ai \curvearrowright 2^{\NN \x \NN}$, which restricts to an action $\Ai \curvearrowright \cE(\NN)$.
\begin{definition}\label{def:orbit_map}
    Let $\cO : \Sub{\Ai} \to \cE(\NN)$  be the \emph{orbit equivalence map} given by
    \[
    \cO(H) = \setcomp{(i, \sigma \cdot i)}{i \in \NN, \sigma \in H} \quad \forall H \in \Sub{\Ai}.
    \] 
    Let $\cN : \cE(\NN) \to \NN \union \set{\infty}$  be the map which takes an equivalence relation and returns the number of its equivalence classes.
\end{definition}

We define a map $\varrho : (\NN_0)^\NN \to \cE(\NN)$, as follows. The equivalence relation $\varrho(c)$ associated to a given coloring $c \in (\NN_0)^\NN$ corresponds to the partition $P_c$ of the set $\NN$ given by
\[
P_c = \setcomp{B_k(c)}{k \in \NN, \, B_k(c) \neq \varnothing} \union \setcomp{\set{i}}{ i \in B_0(c) }.
\]

Let $\alpha \in \cD$ be a fixed sequence.
The pushforward of the probability measure $\mu_\alpha$ under the map $\varrho$ will be denoted by $\lambda_\alpha$ (the collection $\cD$ and the measure $\mu_\alpha$ were defined in \Cref{subsub:IRS of Ai}).
The following proposition (whose proof is immediate from the definitions) summarizes the results that we will need on the distributions of the maps $\cO$ and $\cN$.

\begin{proposition}\label{prop:orbit_and_num_of_orbits_properties}
    The maps $\cO$ and $\cN$ are Borel. The map $\cO$ is $\Ai$-equivariant, and the map $\cN$ is $\Ai$-invariant.
    Let $\alpha \in \cD$ be a sequence and $T \in \Sub{E_\alpha}$ be a subgroup. The associated random subgroup $\nu_\alpha^T$ satisfies $\cO_* \nu_{\alpha}^T = \lambda_\alpha$.
    The map $\cN$ is $\lambda_\alpha$-a.s. equal to the constant value $\abs{I_\alpha}$ if $\alpha_0 = 0$, or to the constant value $\infty$ if $\alpha_0 > 0$. 
\end{proposition}

We remark that Kingman \cite{Kingman1978RandomPartitions} has shown that the $\lambda_\alpha$'s are exactly all the ergodic $S_\infty$-invariant measures on $\cE(\NN)$. We will use a slight variation of this theorem in the next discussion on pointed partitions.

\subsection{Pointed partitions}
\label{subsubsection:pointed_partitions}
We  digress further to overview the notion of pointed partitions of $\NN$ and of the classification of ergodic $\Ai$-measures on the space of such partitions. This is a special case of the notion of signed partitions, as developed by Vershik in \cite{Vershik2012}.  

\begin{definition}
    A \textit{pointed partition} of the set $\NN$ is a finite or countable partition $P$ together with a distinguished block $B_P \in P$.
    We denote the set of all pointed partitions of $\NN$ by $\cE^*(\NN)$.
\end{definition}

The set $\cE^*(\NN)$ is a Borel space admitting a natural action of the group  $S_\infty$ of finitary permutations through its action on the set $\NN$. We will describe all the ergodic  $S_\infty$-invariant probability measures on $\cE^*(\NN)$. These should be compared with the ergodic invariant measures on the set $\cE(\NN)$ discussed in the previous  \Cref{subsubsection:the_orbit_map}.

Here is one way to construct pointed partitions. Given a point $c \in (\ZZ_{\geq -1})^\NN$, we denote $B_k(c) = \setcomp{n \in \NN}{c_n = k}$ for all $k \in \NN_{\ge -1}$.
Let $P_c \in \cE^*(\NN)$ be the pointed partition
\[
P_c = \setcomp{B_k(c)}{k \in \NN \union \set{-1}} \union \setcomp{\set{i}}{ i \in B_0(c) },
\]
whose distinguished block is $B_{-1}(c)$.

Let $\cD^* \subseteq [0, 1]^{\ZZ_{\geq -1}}$ be the set of sequences of non-negative real numbers $(\alpha_i)_{i\in \ZZ_{\geq -1}}$ which satisfy:
\begin{itemize}
    \item The subsequence $(\alpha_{i})_{i=1}^\infty$ is non-increasing, and
    \item $\sum_{i=-1}^\infty \alpha_i = 1$.
\end{itemize}
Given each sequence $\alpha \in \cD^*$ let $p_\alpha$ be the  probability measure on $\ZZ_{\geq -1}$ defined by $p_\alpha(i) = \alpha_i$ for all $i \in \ZZ_{\ge -1}$. Let $\mu_\alpha$ be the corresponding product measure on $(\ZZ_{\geq -1})^\NN$.
The pushforward of the probability measure $\mu_\alpha$ under the map $c \mapsto P_c$ will be denoted $\lambda_\alpha^*$.

\begin{proposition}\label{prop:e_star_ergodic_measures}
    Any ergodic $\Ai$-invariant probability measure on the space $\cE^*(\NN)$ of pointed partitions is of the form $\lambda_\alpha^*$ for some sequence $\alpha \in \cD^*$.
\end{proposition}

    This  is a slight refinement of Kingman's classification of the ergodic $S_\infty$-measures on $\cE(\NN)$ \cite{Kingman1978RandomPartitions}.
    We choose to reduce the proposition to Vershik's classification of invariant measures on \textit{signed} partitions \cite{Vershik2012}.
    A \textit{signed partition} of the set $\NN$ is a partition together with a designation of ``positive'' or ``negative'' to each non-singleton block of the partition.
    Denote the set of signed partitions of $\NN$ by $\cE^\pm(\NN)$.
    Given a fixed probability vector $\alpha \in [0, 1]^\ZZ$, one can define a probability measure $\lambda_\alpha^{\pm}$ on the space $\cE^\pm(\NN)$ in a manner very similar to the measures $\lambda_\alpha^*$ and $\lambda_\alpha$: 
    choose a random ``color vector'' $c \in \ZZ^\NN$ in an i.i.d.\ manner according to $\alpha$.
    The corresponding random partition $P_c$ is
    \[
        P_c = \setcomp{B_k(c)}{k \in \ZZ \setminus \set{0}} \union \setcomp{\set{i}}{ i \in B_0(c) }
    \]
    The positive blocks are the $B_k(c)$'s with $k > 0$, and the negative blocks are the $B_k(c)$'s with  $k < 0$.
    Vershik has shown that every ergodic $S_\infty$-invariant  measure on $\cE^\pm(\NN)$ is of the form $\lambda_\alpha^\pm$ for some such fixed probability vector $\alpha$. His proof works in the same way for $\Ai$-invariance.
    
\begin{proof}[Proof of \Cref{prop:e_star_ergodic_measures}]
    Let $\lambda$ be an ergodic $\Ai$-invariant probability measure on the space $\cE^*(\NN)$. We claim that the distinguished block is $\lambda$-a.s not a singleton.  Otherwise,  the pushforward of $\lambda$ under the $\lambda$-measurable map $\cE^*(\NN) \to \NN$ which takes a pointed partition to the unique element in the singleton distinguished block would give rise to  an $\Ai$-invariant probability measure on $\NN$, which is a contradiction.

    It follows that the  map $s: \cE^*(\NN) \to \cE^\pm(\NN)$ which designates the distinguished block of the partition as negative and all other non-singleton blocks as positive is $\lambda$-a.s.\ well-defined. 
    It is also $\Ai$-equivariant. Hence $s_* \lambda = \lambda^\pm_\alpha$ for some probability vector $\alpha \in \left[0,1\right]^\ZZ$ according to Vershik's classification.
    Observe that $\alpha$ is supported on $\ZZ_{\geq -1}$. We conclude that $\lambda = \lambda^*_\alpha$ as required.
\end{proof}

Given a pointed partition $P \in \cE^*(\NN)$ and an index $k \in \NN_{\ge -1}$ let $d_k(P) \in [0, 1]$ be the asymptotic density of the $k$-th block of $P$ (whenever it exists). The non-distinguished blocks have $k \geq 1$ and are ordered by non-increasing density. The distinguished block has $k = -1$. Write $D(P) = (d_k(P))_{k \geq -1}$ for the asymptotic density vector (whenever these exist). 

 Consider any fixed sequence  $\beta \in \cD^*$. For $\lambda_\beta^*$-a.e.\ pointed partition $P \in \cE^*(\NN)$, the asymptotic density vector $D(P)$ exists and equals $\beta$. This is immediate from the strong law of large numbers applied to the i.i.d.\ sequence of labels $c \in (\ZZ_{\geq -1})^\NN$ used to construct the random pointed partition $P$ with respect to the probability measure $\mu_\beta$.

\begin{lemma}\label{lemma:density_determines_ergodic}
    Let $\lambda$ be an $\Ai$-invariant probability measure on $\cE^*(\NN)$ and $\alpha \in \cD^*$ be a fixed sequence. If the asymptotic density vector $D(P)$ is $\lambda$-a.s.\ equal to $\alpha$ then $\lambda = \lambda_\alpha^*$ and $\lambda$ is ergodic.
\end{lemma}
\begin{proof}
    By \Cref{prop:e_star_ergodic_measures}, the  probability measure $\lambda$ admits an ergodic decomposition of the form $\lambda = \int_{\cD^*} \lambda_\beta^* \diff m(\beta)$ for some probability measure $m$ on $\cD^*$. Combined with the  observation preceding the lemma, this decomposition gives $D_* \lambda = \int \delta_\beta \diff m(\beta) = m$.
    If $D(P)$ is $\lambda$-a.s.\ equal to the sequence $\alpha$ then necessarily $m = \delta_\alpha$. This implies that $\lambda = \lambda_\alpha^*$ and in particular that $\lambda$ is ergodic.
\end{proof}

\section{The Vershik trace of an invariant random subgroup of \texorpdfstring{$\SLi$}{\SLipdf}}
\label{section:singer_cycles}
Let $\mathbb{F}_q$ be a finite field of order $q$ and $\SLi$ the finitary special linear group. Let $\mu$ be an invariant random subgroup of  $\SLi$. 
In this section we use Skudlarek's classification of Thoma traces of $\SLi$ (see \Cref{thm:sl_infty_characters}) to get information on $\mu$. 
This is done by taking advantage of the ergodic decomposition of the associated Vershik trace $\chi_\mu$.

\begin{example}
Observe that any character of $\SLi$ other than $\delta_I$ is strictly positive on the entire group $\SLi$. As long as $\mu \neq \delta_\triv $ we have $\chi_\mu \neq \delta_I$ and hence the Vershik trace $\chi_\mu$ is also strictly positive.
This means that \emph{every} element of $\SLi$ belongs to a $\mu$-random subgroup with a positive probability. 
Consequently
each conjugacy class of $\SLi$ intersects $\mu$-almost every random subgroup.
\end{example}

Here is a more elaborate use of the Vershik trace. Regard  $\SLn$ as a subgroup of $\SLi$ in the standard way (i.e. in the \enquote{upper-left corner}).  
We will see in \Cref{prop:sln_probability} that the probability that a $\mu$-random subgroup $H$ contains $\SLn$ is approximately the same as the probability that a matrix $g$ with $\rkt(g) = \rk(g-I) = n$ belongs to $H$. 
This will imply in \Cref{prop:sl_n_exists_in_some_basis} that $\mu$-almost every subgroup contains a copy of $\SLn$ in some finitary basis. See \Cref{sub:IRS contain SLn} for more details on this direction.

In \Cref{subsection:irs_restriction}, we use the theory of invariant random subgroups and Thoma characters of the group $\Ai$ discussed in \Cref{subsection:irs_a_infty_description} to identify the image of the invariant random subgroup $\mu$ under the intersection map $\IRS(\SLi) \to \IRS(\Ai)$. We think of $\Ai$ as a subgroup of $\SLi$ consisting of permutation matrices in some finitary basis.

\subsection{The degree of an invariant random subgroup}
Let  $\chi_\mu$ be the Vershik trace associated to the invariant random subgroup $\mu$ of the group $\SLi$.
By Skudlarek's classification (\Cref{thm:sl_infty_characters}), the ergodic decomposition of the trace $\chi_\mu$ has the form
\begin{equation}\label{eq:chi_formula}  
    \chi_\mu = \sum_{k = 0}^\infty a_k(\mu) \chi_k + a_\infty(\mu) \chi_\infty,
\end{equation}
where $a_k(\mu) \geq 0$ are non-negative real numbers satisfying $\sum_{k=0}^\infty a_k(\mu) + a_\infty(\mu) = 1$. 
 Recall that \[\chi_k(g) = q^{-k \cdot \rkt g} \quad \forall g \in \SLi\] for all $k \in \NN$, and $\chi_\infty = \delta_I$.
 
\begin{definition}
\label{def:degre}
    The \emph{degree} of the invariant random subgroup $\mu$ is the minimal number $r(\mu) \in \NN \cup \{\infty\}$  such that $a_{r(\mu)}(\mu) > 0$. 
\end{definition}

Note that the degree $r(\mu)$ is infinite if and only if $\chi_\mu = \delta_I$. Since the group $\SLi$ is countable, this happens if and only if $\mu = \delta_\triv$.

We observe several additional useful properties of the Vershik trace $\chi_\mu$.
First, the value of any character of the group $\SLi$ at any element $g \in \SLi$ depends solely on the quantity $\rkt g$. Being a linear combination of such characters, the Vershik trace $\chi_\mu$ is consequently also a function of this quantity.

\begin{definition}
Given any  $n \in \NN_0$, we denote by $\chi_\mu(n)$ the common value of the Vershik trace $\chi_\mu$ on any  element $g \in \SLi$ having $\rkt(g) = n$.
\end{definition}

The second property is $\chi_\mu$'s growth rate, regarded as a function of $ n = \rkt g$ in the above sense. 
Asymptotically, as $n = \rkt g \to \infty$ we get
\begin{equation} \label{eq:chi_asymptotic}
    \chi_\mu(n) = a_{r(\mu)}(\mu) \cdot q^{-r(\mu) n} (1 + o(1)).
\end{equation}
In particular, we see that $\chi_\mu(n)$ decays exponentially in $n$, with exponential rate $q^{-r(\mu)}$.

\subsection{Invariant random subgroups contain $\SLn$}
\label{sub:IRS contain SLn}
Singer cycles were discussed in \Cref{subsection:singer_cycles}.
We observe that the size of the conjugacy class of a Singer cycle in $\SLn$ is larger than the order of any maximal proper subgroup containing it. Given an invariant random subgroup $\mu$, each of these conjugates has an equal probability of belonging to a $\mu$-random subgroup. By the pigeonhole principle, sufficiently many conjugates must appear in any $\mu$-random subgroup to generate an isomorphic copy of $\SLn$.
\begin{proposition}\label{prop:sln_probability}
    Let $\mu$ be an invariant random subgroup of $\SLi$ and  $f$ be a finitary basis.
    Then
    \begin{equation}
        \mu\left( \setcomp{H \in \Sub{\SLi}}{\SL_n(f) \leq H} \right) = \chi_\mu(n)(1 - o(1)).
    \end{equation}
\end{proposition}

\begin{proof}
    For every $n \in \NN$ let $s_n \in \SL_n(f)$ be a arbitrary Singer cycle. Since $1$ is not an eigenvalue of $s_n$, it satisfies $\rkt(s_n) = n$. Hence
    \[
    \mu\left( \setcomp{H}{\SL_n(f) \leq H} \right) \leq \mu\left( \setcomp{H}{s_n \in H} \right) = \chi_\mu(n).
    \]
    We need to show that the other inequality holds, up to a factor of $(1-o(1))$.
    For $\mu = \delta_\triv$ both sides vanish for all $n \in \NN$. So we may assume that $r(\mu) < \infty$.

    Let $N_n$ be the size of the conjugacy class  of $s_n$ in $\SL_n(f)$. We know that $N_n$  is $q^{n^2 + O(n)}$ by \cref{eq:singer_cycle_cc_size}.
    By \Cref{prop:singer_cycle_subgroup_size} any $k_n=q^{\frac 3 4 n^2}$ distinct conjugates of $s_n$ generate $\SL_n(f)$, for all large enough $n$.
    Apply the quantitative pigeonhole principle (see \Cref{lemma:pigeonhole_principle} below) to get
    \begin{align*}
        \mu\left( \setcomp{H}{\SL_n(f) \leq H} \right) &\geq \mu\left( \set{k_n + 1 \text{ distinct conjugates of } s_n \text{ belong to } H} \right)
        \\&\geq \frac{N_n \chi_\mu(s_n) - k_n} {N_n - k_n} \geq \chi_\mu(n) - \frac{k_n}{N_n}
        \\&= \chi_{\mu}(n) - q^{-\frac 1 4 n^2 + O(n)} = \chi_\mu(n)(1 - o(1)).
    \end{align*}
    In the last equality 
     we use the fact that $\chi_\mu(n)$ decays exponentially in $n$ (\cref{eq:chi_asymptotic}).  This concludes the proof.
\end{proof}

In the above proof, we have used the following quantitative version of the pigeonhole principle for probability spaces.

\begin{lemma} \label{lemma:pigeonhole_principle}
Let $(X, \cB, \PP)$ be a probability space. Let $N \in \NN$ and $E_1, ..., E_N$ be events in $\cB$. 
Then for any $0 \leq k < N$ we have
\begin{equation} \label{eq:pigeonhole_eq}
    \PP\paren{\sum_{i=1}^N \one_{E_i} > k} \geq \frac{\sum_{i=1}^{N} \PP\paren{E_i} - k}{N - k}.
\end{equation}
\end{lemma}

\begin{proof}
    Let $F = \setcomp{ x \in X}{\sum_{i=1}^N \one_{E_i}(x) > k}$. Then we have
    \begin{align*}
            \sum_{i=1}^N \PP\paren{E_i} &= \EE{\sum_{i=1}^N \one_{E_i}} =
            \EE{\sum_{i=1}^N \one_{E_i} \one_{F} + \sum_{i=1}^N \one_{E_i} \one_{\comp{F}}}\\
            &\leq \EE{N \one_F + k \one_{\comp{F}}}
             = N \PP\paren{F} + k(1 - \PP\paren{F}).
    \end{align*}
    The desired \cref{eq:pigeonhole_eq} follows by rearranging.
\end{proof}

We proceed to consider some of the consequences of \Cref{prop:sln_probability}.

\begin{proposition}\label{prop:sl_n_exists_in_some_basis}
    Let $\mu \in \IRS(\SLi)$ be a non-trivial ergodic invariant random subgroup. Then for every $n \in \NN$ and for $\mu$-almost every subgroup $H \in \Sub{\SLi}$ 
    there exists a finitary basis $f$ such that $\SL_n(f) \leq H$.
\end{proposition}

\begin{proof}
For each $n \in \NN$ the subset 
\[
E_n = \setcomp{H \in \Sub{\SLi}}{\text{there exists a finitary basis $f$ such that $\SL_n(f) \leq H$}}
\]
is $\SLi$-invariant. Hence by ergodicity $\mu(E_n) \in \set{0, 1}$. \Cref{prop:sln_probability} implies that  $\mu \left( \setcomp{H}{\SLn \leq H} \right) > 0$ for all $n$ which are greater or equal to some $n_0 \in \NN$.
    For $n < n_0$ we note that
    \[\mu \left( \setcomp{H}{\SLn \leq H} \right) \geq \mu \left( \setcomp{H}{\SL_{n_0}(\FFq) \leq H} \right) > 0.\]
    Hence  $\mu(E_n) \geq \mu \left( \setcomp{H}{\SLn \leq H} \right) > 0$. So $\mu(E_n) = 1$ for all $n \in \NN$.
\end{proof}

\begin{remark}
    We obtain the tentative result saying that a $\mu$-random subgroup $H$ with respect to a non-trivial invariant random subgroup $\mu$ of $\SLi$ is not solvable.
Of course, this conclusion would a posteriori follow from our main result \Cref{thm:main_theorem}, and can be seen as partial progress towards it.
\end{remark}

We complete the classification of ergodic invariant random subgroups of degree $0$.

\begin{proposition}
    Let $\mu \in \IRS(\SLi)$ be an  ergodic invariant random subgroup. If $r(\mu) = 0$ then $\mu = \delta_{\SLi}$.
\end{proposition}
\begin{proof}
    Assume that the degree $r(\mu)$ is $0$. By \Cref{prop:sln_probability} we observe that
    $$
    \mu\left( \setcomp{H}{\SL_\infty(\FFq) \leq H} \right) = \lim_{n \to \infty} \mu\left( \setcomp{H}{\SL_n(\FFq) \leq H} \right)
    = \lim_{n \to \infty} a_0(\mu) (1 - o(1)) = a_0 > 0.
    $$
    Because the event $\setcomp{H}{ \SL_\infty(\FF_q) \leq H }$ is $\SLi$-invariant, it must occur with probability one by ergodicity.
\end{proof}

\subsection{The intersection map from $\SLi$ to $\Ai$} \label{subsection:irs_restriction}

We use the Vershik trace of the given invariant random subgroup $\mu$ to identify its image $\nu_f$  under the intersection map $\IRS(\SLi) \to \IRS(\Ai(f))$ given a finitary basis $f$. 

We use notations from \Cref{subsection:irs_a_infty_description}.
For each $k \in \NN$ define $\alpha^k \in \cD$ to be the sequence given by
    \[
        \alpha^k_i = 
        \begin{cases}
            1/q^k, & 1 \leq i \leq q^k, \\
            0, & \text{$i=0$ or $i > q^k$}.
        \end{cases}
    \]
     Define $\alpha^\infty \in \cD$ to be the sequence given by $\alpha^\infty_0 = 1$ and $\alpha^\infty_i = 0$ for all $i \in \NN$. For each $k \in \NN_0 \union \set{\infty}$ consider the invariant random subgroup \[\nu_k = \nu_{\alpha^k}^{E_{\alpha^k}} \in \IRS(\Ai).\]
     In other words $\nu_k$ is the stabilizer distribution of a random coloring of the set $\NN$ using $q^k$ distinct colors independently and with equal probabilities.
\begin{proposition}\label{prop:irs_intersection_result}
    Let $\mu$ be an invariant random subgroup of $\SLi$ of  degree $r(\mu) < \infty$. 
    For every finitary basis $f$ the image $\nu_f$ of $\mu$ under the  intersection map $\IRS(\SLi) \to \IRS(\Ai(f))$ admits  the ergodic decomposition\begin{equation}\label{eq:nu_decomposition}
        \nu_f = \sum_{k = r(\mu)}^\infty a_k(\mu) \nu_k + a_\infty(\mu) \nu_\infty.
    \end{equation}
\end{proposition}
    The coefficients  $a_k(\mu)$ appearing in \cref{eq:nu_decomposition} correspond to the ergodic decomposition of the Vershik trace $\chi_\mu$ as in \cref{eq:chi_formula}.

\begin{proof}[Proof of \Cref{prop:irs_intersection_result}]
Fix a finitary basis $f$. To ease our notations, set $r = r(\mu)$, $a_k =a_k(\mu)$ for all $k$, and $\nu = \nu_f$.
    We first calculate the Vershik trace $\chi_\nu$ associated to $\nu$. Notice that $\chi_\nu$ is by definition the restriction to the subgroup $\Ai(f)$ of  the Vershik trace $\chi_\mu$ associated to  $\mu$. Let $\sigma \in \Ai$ be any non-trivial element. Let $\sigma_1, ..., \sigma_t$ be the distinct cycles of $\sigma$ and let $\ell_i$ be the length of $\sigma_i$. In the permutation representation we have  $\rkt \sigma _i = \ell_i - 1$.
    Therefore 
    \begin{align*}
        \chi_\nu(\sigma) = \chi_\mu(\sigma) 
        &= \sum_{k=r}^\infty a_k q^{-k \cdot\rkt \sigma} = \sum_{k=r}^\infty  a_k \prod_{\ell = 2}^\infty \left( q^{-k(\ell - 1)} \right)^{c_\ell(\sigma)}
        \\ &= \sum_{k=r}^\infty  a_k \prod_{\ell = 2}^\infty \left( \sum_{i=1}^{q^k} \left( \frac{1}{q^k} \right)^\ell \right)^{c_\ell(\sigma)}
        = \sum_{k=r}^\infty a_k \chi_{(\alpha^k; 0)} (\sigma)
        \\ &= \sum_{k=r}^\infty a_k \chi_{(\alpha^k; 0)} (\sigma) + a_\infty \chi_{(\alpha^\infty; 0)}(\sigma).
    \end{align*}
    This relation clearly holds also for the trivial permutation $\sigma = \id$. Overall,   we get
    \begin{equation} \label{eq:ainfty_character_decomposition}
        \chi_\nu = \sum_{k=r}^\infty a_k \chi_{(\alpha^k; 0)} + a_\infty \chi_{(\alpha^\infty; 0)}.
    \end{equation}
    By Choquet theory, there exists a probability measure $m$ on the ergodic invariant random subgroups of $\Ai$ such that 
    \[
        \nu = \int_{\Extr(\IRS(\Ai))} \lambda \diff m(\lambda).
    \]
 Let $\tr : \IRS(\Ai) \to \Tr(\Ai)$ denote the Vershik trace map, i.e.\ $\tr(\lambda) = \chi_\lambda$. Applying this map to the above, we get    \begin{equation} \label{eq:ainfty_irs_character_decomposition}
        \tr(\nu) = \chi_\nu = \int_{\Extr(\IRS(\Ai))} \chi_{\lambda} \diff m(\lambda) = \int_{\Tr(\Ai)} \chi \diff (\tr_* m) (\chi).
    \end{equation}   
    Recall from \Cref{subsection:thoma and vershik traces} that the set of traces of a group forms a Choquet simplex.
    That is, the ergodic decomposition of a trace is unique.
    Comparing \cref{eq:ainfty_irs_character_decomposition} to \cref{eq:ainfty_character_decomposition}, we see that $\tr_* m$ must be supported on characters 
    whose ergodic decomposition is contained in the set $\setcomp{\chi_{(\alpha^k; 0)}}{k \in \NN_{\geq r} \union \set{\infty}}$.
    It is also supported on the Vershik traces of ergodic invariant random subgroups of $\Ai$.
    By \Cref{thm:ainfty_irs_vershik_character} and the formula preceding it, we must have $\nu = \sum_{k = r}^\infty a_k \nu_{\alpha^k}^{E_{\alpha^k}} + a_\infty \nu_{\alpha^\infty}^{E_{\alpha^\infty}}$, as required.
\end{proof}

Let us work out some consequences of the previous proposition.

\begin{proposition}\label{prop:color_basic_definitions_strengthened}
    Let $\mu \in \IRS(\SLi)$ be an invariant random subgroup.
    Then $\mu$-almost every non-trivial subgroup $H \in \Sub{\SLi}$ satisfies  $H \inter \Ai(f) \neq \triv$ for every finitary basis $f$. 
\end{proposition}

In particular, if $\mu$ is ergodic and non-trivial then $\mu$-almost every subgroup $H$ satisfies  $H \inter \Ai(f) \neq \triv$ for all finitary bases $f$.

\begin{proof}[Proof of \Cref{prop:color_basic_definitions_strengthened}]
Let $f$ be any fixed finitary basis.
    Making use of the ergodic decomposition, we may assume without loss of generality that  $\mu$ is ergodic and that  $\mu \neq \delta_\triv$. 
We know by \Cref{prop:sl_n_exists_in_some_basis} that $\mu$-almost every subgroup $H$ admits some finitary basis $f'$ satisfying $\SL_3(f') \leq H$. The permutation matrices in $\SL_3(f') \cap \Ai(f')$ show that the intersection $H \cap \Ai(f')$ is non-trivial. \Cref{prop:irs_intersection_result} implies that $\mu$-almost surely there is an infinite subset  $c \subseteq \NN$ such that $A_c(f') \leq H$ (such a subset will be called a \emph{color} in \Cref{section:generalities_on_colors}).
    The two bases $f$ and $f'$ agree starting from some basis element, namely $f^n = f'^n$ for all $n \ge n_0$ for some $n_0 \in \NN$.
    Thus
    \[
        A_{c \inter \NN_{\geq n_0}}(f) = A_{c \inter \NN_{\geq n_0}}(f') \leq H.
    \]
    This is a non-trivial subgroup containing permutation matrices in the basis $f$. We conclude that  $\Ai(f) \inter H \neq \triv$, as required.
\end{proof}

\begin{corollary}
Let $\mu \in \IRS(\SLi)$ be an ergodic invariant random subgroup.
If $\mu$ is non-trivial then  $a_\infty(\mu) = 0$.
\label{cor:non-trivial ergodic has no k=infty}
\end{corollary}

\begin{proof}
    The ergodic decomposition of the image of the intersection map (see \Cref{prop:irs_intersection_result}) shows that $a_\infty(\mu)$ is equal to the probability that 
    $H \inter \Ai(e) = \triv$, where $e$ is the standard basis.
    By \Cref{prop:color_basic_definitions_strengthened} this condition is $\mu$-almost surely equivalent to saying that $H = \triv$. The ergodicity of $\mu$ implies that  $\mu(\set{\set{I}}) \in \set{0,1}$. Since $\mu$ is non-trivial, this probability must be $0$.
\end{proof}

\section{The coloring equivalence relation}
\label{section:generalities_on_colors}
In this section we define the  coloring equivalence relation corresponding to a subgroup of $\SLi$ and a finitary basis $f$. We define the corresponding ``colors'' as the equivalence classes of this relation. 
The basic properties of the coloring, combined with the results from \Cref{section:singer_cycles}, will allow us to show in \Cref{prop:num_of_colors_constant} that the Vershik trace of a non-trivial ergodic invariant random subgroup $\mu$ is precisely $\chi_{r(\mu)}$  where $r(\mu)$ is the degree. Moreover, we define an action of the group $\GLi$ on color vectors and examine its properties. This action will be shown to satisfy the assumptions of our abstract \Cref{thm:action_classification}, thereby allowing us to equivariantly identify the space of color vectors $\Col(H)^\NN$ with the vector space $V^{k,\l}$.

\subsection{Definition of colors} \label{subsection:color_definition}
The primary focus of this section is the following object. We use notations from \Cref{subsubsection:the_orbit_map}.

\begin{definition}\label{def:color_definition}
    Let $f$ be a finitary basis.
    \begin{itemize}
    \item The \textit{coloring relation map} is given by
    \[
        R_f : \Sub{\SLi} \to \cE(\NN), \quad R_f(H) = \cO(H \inter \Ai(f))
    \]
    where $\cO$ is the orbit map (\Cref{def:orbit_map}).
 In other words, the map $R_f$ associates to each subgroup of $\SLi$ the orbit equivalence relation of the even permutation matrices in basis $f$ contained in that subgroup.
    
    \item The \textit{coloring partition of $H$ in the basis $f$} is the partition of $\NN$ corresponding to the equivalence relation $R_f(H)$. It will be denoted by $\Col_f(H)$. The elements of $\Col_f(H)$ are the \textit{colors of $H$ in the basis $f$}.

       \item An index $i \in \NN$ is \textit{colored in the color $c$ for the subgroup $H$ in the basis $f$} if  $c \in \Col_f(H)$ and $i \in c$. 

\end{itemize}
\end{definition}

Given a non-trivial ergodic invariant random subgroup $\mu \in \IRS(\SLi)$ and a finitary basis $f$,  we can control  the distribution of the equivalence relation $R_f$. For each $k \in \NN$ let $\alpha^k \in \cD$ be the sequence defined in \Cref{subsection:irs_restriction}. Let $\lambda_k = \lambda_{\alpha^k}$ be the $\Ai$-invariant ergodic probability measure on $\cE(\NN)$ associated to the sequence $\alpha^k$ (see \Cref{subsubsection:the_orbit_map}).
Combining  
\Cref{prop:orbit_and_num_of_orbits_properties}, the information from \Cref{subsection:irs_restriction} and taking into account all of the involved definitions gives the decomposition
\begin{equation} \label{eq:color_distribution}
    \left(R_f\right)_* \mu = \sum_{k=r(\mu)}^{\infty} a_k(\mu) \lambda_k.
\end{equation}
We immediately obtain the following.
\begin{proposition} \label{prop:color_basic_definitions}
    Let $\mu$ be an invariant random subgroup of $\SLi$ and  $f$ be a finitary basis. Then for $\mu$-almost every non-trivial subgroup $H \in \Sub{\SLi}$ 
    \begin{enumerate}
        \item  any color in $\Col_f(H)$ is an infinite subset of $\NN$, and
        \item  a permutation matrix $\sigma \in \Ai(f)$ belongs to $H$ if and only if every cycle of $\sigma$ is colored the same color in the basis $f$ for the subgroup $H$. In particular, if $c \in \Col_f(H)$ then $A_c(f) \leq H$. 
    \end{enumerate}

\end{proposition} 

Any finitary basis $f$ coincides with  the standard basis $e$ in all but a finite number of vectors. Hence for a $\mu$-generic non-trivial subgroup $H$,  every color in $\Col_f(H)$ differs by a finite number of elements from a color of $\Col_e(H)$.
Since all colors are infinite, we obtain a natural identification of the colors in $\Col_e(H)$ with the colors in $\Col_f(H)$. This observation allows us to talk of colors without having to specify a particular finitary basis. 

\begin{definition} 
\label{def:color_definition part 2}
Let $\mu \in \IRS(\SLi)$ be an invariant random subgroup. Let $H$ be a  $\mu$-generic non-trivial subgroup.
\begin{itemize}
    \item The colors $\Col(H)$ of the subgroup $H$ are its colors in any finitary basis, up to the above natural identification.

    \item  
    The \textit{color vector of the subgroup $H$ in the finitary basis $f$} is the point in $\prod_{i \in \NN} \Col(H)$ whose $i$th entry is the color of the index $i \in \NN$ in the basis $f$. It will be denoted by $\col_f(H)$. 
\end{itemize}
\end{definition}

Note that the color vector of a subgroup does depend on the choice of a finitary basis (even if its set of colors does not). 
The number of colors $\abs{\Col(H)}$ of a $\mu$-generic non-trivial subgroup is a power of $q$.

\begin{remark}
For the purpose of $\mu$-measurability, we may identify the set of colors $\Col(H)$ of a $\mu$-generic non-trivial subgroup with its set of colors $\Col_e(H)$ in the standard basis $e$. The latter partition may be identified with the equivalence relation it induces, the coloring relation $R_e(H)$. 
\end{remark}




\subsection{The number of colors} 
We will now define a $\sigma$-algebra on $\Sub{\SLi}$ consisting of events which depend only on colors of far-away elements of finitary bases.
We show that the number of colors is measurable with respect to this $\sigma$-algebra. 
This allows us to prove that an ergodic invariant random subgroup has Vershik trace equal precisely to $\chi_{r(\mu)}$.

\begin{definition}\label{def:tail_color_sigma_algebra}
Fix an arbitrary finitary basis $f$. 
Let $R^{>n}_f : \Sub{\SLi} \to \cE(\NN_{>n})$ be the function that maps each subgroup $H \in \Sub{\SLi}$ to the restriction of the equivalence relation $R_f(H)$ to $\NN_{>n}$.
    Define the $\sigma$-algebra $\cT_n$ to be the $\sigma$-algebra on $\Sub{\SLi}$ generated by the map $R^{> n}_f$. The \textit{tail color $\sigma$-algebra} is the intersection $\cT = \bigcap_{n \in \NN} \cT_n$.
\end{definition}

Note that the tail color $\sigma$-algebra is independent of the choice of $f$.

\begin{proposition}\label{prop:tail_is_trivial}
    Let $\mu$ be an invariant random subgroup of $\SLi$. 
    For each $n \in \NN$ the $\sigma$-algebra $\cT_n$ is  $\SLn$-invariant modulo $\mu$-null sets, namely, for all $E \in \cT_n$ and $g \in \SLn$ we have 
    $\mu(E \triangle g E g^{-1}) = 0$. 
\end{proposition}

It follows that 
the $\sigma$-algebra $\cT$ is  $\SLi$-invariant modulo $\mu$-null sets.
In particular, if $\mu$ is ergodic then $\cT$ is $\mu$-trivial, in the sense that $\mu(E) \in \set{0, 1}$ for all $E \in \cT$. 

\begin{proof}
[Proof of \Cref{prop:tail_is_trivial}]
Fix an arbitrary finitary basis $f$.
First, we claim that $\mu$-almost every  subgroup $H \in \Sub{\SLi}$ satisfies the following property:
    if $n < i < j \in \NN$ are three indices and $g \in \SLn$ is an element 
    then $(i, j) \in R_f(H)$ if and only if $(i, j) \in R_f(g H g^{-1})$.
    Indeed, let $H$ be a $\mu$-generic subgroup, and suppose $i$ and $j$ are colored in the same color in the basis $f$ for the subgroup $H$. 
    The intersection $H \inter \Ai(f)$ cannot be trivial, since in that case all colors are singletons. Hence the common color $c$ of $i$ and $j$ is infinite, and $A_c(f) \subseteq H$. 
    Let $\sigma \in A_{c \inter \NN_{> n}}(f)$ be any permutation matrix which swaps $i$ and $j$.
    Then $\sigma$ commutes with $g$, i.e.\ $\sigma = g \sigma g^{-1} \in g H g^{-1}$.
    This implies $(i, j) \in R_f(g H g^{-1})$, which shows one direction of the claim. 
    The reverse direction is established analogously, by using the fact that the conjugate $g^{-1} H g$ is also $\mu$-generic.

    We now exhibit a $\pi$-system which generates the $\sigma$-algebra $\cT_n$ and is  $\SLn$-invariant modulo $\mu$-null sets. 
    We look at events of the following form.
    Let $(i_1, j_1), \ldots, (i_a, j_a)$ and $(s_1, t_1), \ldots, (s_b, t_b)$ for some $a,b \in \NN$ be pairs of indices  all greater than $n$. 
    We define the event
    \[
        E = \setcomp{H \in \Sub{\SLi} }{\forall 1\le k \le a, (i_k, j_k) \in R_f(H)\, \wedge\, \forall 1 \le \l \le b, (s_\l, t_\l) \notin R_f(H)}.
    \]
    The set of all such events $E$ is a generating $\pi$-system for $\cT_n$. 
    The previous claim shows that $\mu(g E g^{-1} \triangle E) = 0$ for all $g \in \SLn$, as required.
\end{proof}

\begin{proposition}\label{prop:number_of_colors_is_as_const}
    Let $\mu \in \IRS(\SLi)$ be a non-trivial ergodic invariant random subgroup.
    The number of colors $\abs{\Col(H)}$ of the subgroup $H$  is $\mu$-almost surely constant.
\end{proposition}

\begin{proof}
We claim that the \emph{number-of-colors} map 
\[\mathcal{NC} : \Sub{\SLi} \to \NN \cup \{\infty\}, \quad \mathcal{NC} : H \mapsto |\Col(H)|\]
is $\cT$-measurable. 
Recall the function $\cN : \cE(\NN) \to \NN \union \set{\infty}$ which counts the number of equivalence classes of an equivalence relation (\Cref{def:orbit_map}).
Fix a finitary basis $f$. We may compose the map $\cN$ with the coloring relation map $R_f : \Sub{\SLi} \to \cE(\NN)$ to get the \textit{number-of-colors in the basis $f$} 
\[
    \cN_f : \Sub{\SLi} \to \NN_0 \union \set{\infty}, \quad \cN_f = \cN \circ R_f.
\]
Note that the map $\mathcal{NC}$ coincides $\mu$-almost surely with the map $\cN_f$.

For a $\mu$-generic  subgroup $H$, the equivalence classes of $R_f(H)$ are all infinite. Hence their number is the same as the number of equivalence classes of $R_f^{> n}(H)$ for any $n \in \NN$. 
    Thus $\cN_f  = \cN_f^{> n} $, where $\cN_f^{> n} : \Sub{\SLi} \to \NN_0 \union \set{\infty}$ is the map which counts the number of equivalence classes of $R_f^{> n}(H)$. 
    This map is $\cT_n$-measurable, so that $\cN_f $ is $\cT_n$-measurable for any $n \in \NN$.
    Therefore the map $\cN_f $ and hence also $\mathcal{NC}$ is $\cT $-measurable.
    
    Lastly, since $\cT$ is $\mu$-trivial, we conclude that the map $\mathcal{NC}$ is $\mu$-almost surely constant.
\end{proof}

\begin{proposition}\label{prop:num_of_colors_constant}
    Let $\mu$ be a non-trivial ergodic invariant random subgroup of $\SLi$. Let $r = r(\mu)$ denote the degree of $\mu$.
    Then $a_r(\mu) = 1$, that is, the Vershik trace associated to $\mu$ is $\chi_r$.
    The number of colors of $\mu$-almost every subgroup  is $q^r$.
\end{proposition}
\begin{proof}
    Propositions \ref{prop:orbit_and_num_of_orbits_properties}
    and \ref{prop:irs_intersection_result} together imply  that  the value $a_k(\mu)$ is equal to the probability that a $\mu$-random subgroup has $q^k$ colors. 
    Since the number of colors  is $\mu$-almost surely constant, it must be the case that  $a_{r}(\mu) = 1$. Hence $\chi_\mu = \chi_{r}$. 
    Thus the number of colors of a random subgroup is $\mu$-almost surely  equal to $q^{r}$. 
\end{proof}

\subsection{The color \texorpdfstring{$\red$}{\redpdf}}\label{subsection:the_color_red}
Fix throughout the current \Cref{subsection:the_color_red} a non-trivial ergodic invariant random subgroup $\mu \in \IRS(\SLi)$.  


\begin{proposition}\label{prop:red_exists_and_unique}
For $\mu$-almost every subgroup $H \in \Sub{\SLi}$ and for any finitary basis $f$ there exists a \emph{unique} color $c \in \Col_f(H)$ such that $\SL_c(f) \leq H$. The representative of the color $c$ in the set $\Col(H)$ is independent of the choice of the basis $f$.
\end{proposition}

Once established, this result will allow us to introduce the following key notion.

\begin{definition}
The unique color $\red(H) \in \Col(H)$ given in \Cref{prop:red_exists_and_unique} is called the \emph{red color} of the $\mu$-generic subgroup $H$.
\end{definition}
 
The proof of \Cref{prop:red_exists_and_unique} will require some preparation.

\begin{claim}\label{prop:if_sl_n_then_colored}
    Let $f$ be a finitary basis. Then  $\mu$-almost every subgroup $H$  has the following property: if $I \subseteq \NN$ is a subset with $\abs{I} \geq 3$ such that $\SL_I(f) \leq H$ then the indices in $I$ are colored in the same infinite color $c \in \Col_f(H)$  and 
    $\SL_c(f) \leq H$.
\end{claim}
\begin{proof}
Let $H \in \Sub{\SLi}$ be a generic non-trivial subgroup and $I \subset \NN$ be a subset with $\abs{I} \geq 3$ such that $\SL_I(f) \leq H$. 
    The permutation matrices in $\SL_I(f)$ show that 
    all indices $i \in I$  have the same color $c \in \Col_f(H)$, and that this color is infinite, see \Cref{prop:color_basic_definitions}.
    
    We wish to show  that  $\SL_c (f) \leq H$ in the above situation. 
    It suffices to establish the following statement: if $J$ is any finite subset of $c$ with $\abs{J} \geq 3$ such that $\SL_J(f) \leq H$ then for all indices $i \in c \setminus J$ we have $\SL_{J \union \set i}(f) \leq H$. 
    
    We prove the above statement.
    Let $\sigma_0$ be a permutation matrix in the $f$ basis corresponding to a full cycle on the indices $J \union \set{i}$. If $\sigma_0$ is even, we let $\sigma = \sigma_0 \in A_c(f)$.
    Otherwise, multiply $\sigma_0$ by a transposition swapping two arbitrary indices in $c \setminus (J \union \set{i})$, and call this permutation matrix $\sigma \in A_c(f)$.
    Since $A_c(f) \subseteq \Ai(f) \inter H$, it follows that $\sigma \in H$. Conjugate the subgroup $\SL_J(f)$ by the matrix $\sigma$ to get the subgroup $\SL_{J'}(f)$ of $H$, where $J' = \left(J \setminus \set{\sigma(i)} \right) \union \set{i}$. The two subgroups $\SL_J(f)$ and $\SL_{J'}(f)$ of $H$ generate the subgroup $\SL_{J \union J'}(f) = \SL_{J \union \set{i}}(f) \leq H$, as required.
\end{proof}

From what we have established already, the existence of the desired color is easy to show.

\begin{claim}\label{claim:red_exists}
For $\mu$-almost every subgroup $H \in \Sub{\SLi}$ and for any finitary basis $f$ there is a color $c \in \Col_f(H)$ such that $\SL_c(f) \leq H$.     
\end{claim}

\begin{proof}
Fix a finitary basis $f$.    By \Cref{prop:sl_n_exists_in_some_basis} and Claim \ref{prop:if_sl_n_then_colored}, for $\mu$-almost every subgroup $H \in \Sub{\SLi}$ there is some finitary basis $f'$ and some infinite color $c \in \Col(H)$ for which $\SL_c(f') \leq H$. The two  finitary bases $f$ and $f'$ differ only by a finite number of vectors. Say $f^i = f'^i$ for all $i > n$. In this case 
    \[
        \SL_{c \inter \NN_{> n}}(f) = \SL_{c \inter \NN_{> n}}(f') \leq H.
    \]
Another application of Claim \ref{prop:if_sl_n_then_colored} implies that $\SL_{c}(f) \leq H$, as desired.
\end{proof}

The argument to show the uniqueness of the desired color will be quantitative. 
The rough idea is the following: assume we had two distinct colors $c_1, c_2 \in \Col(H)$ with $\SL_{c_i}(f) \leq H$ for $i = 1, 2$. Intuitively speaking,  the probability of $\SL_n(f) \leq H$ would then be $2 q^{-r n}$ (the first $n$ coordinates are either all colored $c_1$ or all colored $c_2$). 
However, according to  \Cref{prop:sln_probability} this probability should be about $q^{-rn}$, a contradiction.

To make this quantitative argument rigorous, we make use of the notion of pointed partitions from \Cref{subsubsection:pointed_partitions}. For a $\mu$-generic subgroup $H \in \Sub{\SLi}$ and a finitary basis $f$ we consider the set of colors 
\[
    \Col_f^\red(H) = \setcomp{c \in \Col_f(H)}{\SL_c(f) \leq H}.
\]
We form the pointed partition $P^\red_f(H) \in \cE^*(\NN)$ which groups together all ``non-red'' colors. That is 
\[
P^\red_f(H) = \setcomp{c}{c \in \Col_f^\red(H)} \union \set{\NN \setminus \bigsqcup_{c \in \Col_f^\red(H)} c} \in \cE^*(\NN)
\]
with $\NN \setminus \bigsqcup_{c \in \Col_f^\red(H)} c$ being the distinguished block.

\begin{claim}
\label{claim:equivariance claim}
Let $f$ be a finitary basis.
The $\mu$-measurable map $P_f^\red : \Sub{\SLi} \to \cE^*(\NN)$ is 
    $\Ai$-equivariant.
    Here $\Ai$ acts on $\Sub{\SLi}$ by conjugation, through its embedding  via permutation matrices in the basis $f$, namely $\Ai \to \Ai(f) \subseteq \SLi$.
\end{claim}
\begin{proof}
    Let $H$ be a $\mu$-generic subgroup, and $\sigma \in \Ai$ be a permutation. Denote  by $\sigma_f$ the permutation matrix in $\Ai(f)$ corresponding to $\sigma$. 
    The claim is equivalent to showing that the condition $c \in \Col_f^\red(H)$ implies $\sigma \cdot c \in \Col_f^\red(\sigma_f H \sigma_f^{-1})$. Notice that $\sigma \cdot c$ is a color of $\sigma_f H \sigma_f^{-1}$. 
    Since $\SL_c (f) \leq H$, we have
    \[
        \SL_{\sigma \cdot c}(f) = \sigma_f \SL_c(f) \sigma_f^{-1} \leq \sigma_f H \sigma_f^{-1}.
    \]
    Hence $\sigma \cdot c \in \Col_f^\red(\sigma_f H \sigma_f^{-1})$, as required.
\end{proof}

We are now ready to conclude the  proof concerning the uniqueness of the red color.

\begin{proof}[Proof of \Cref{prop:red_exists_and_unique}]
    Let $\mu \in \IRS(\SLi)$ be a non-trivial  ergodic invariant random subgroup. Let $f$ be a finitary basis.
    Let 
    \[
        \cN_f^\red : \Sub{\SLi} \to \NN_0 \union \set{\infty}
    \]
    be the map given by $\cN_f^\red(H) = \abs{\Col_f^\red(H)}$ for each subgroup $H \in \Sub{\SLi}$, i.e. the number of candidate red colors.
    Similarly to the proof of \Cref{prop:number_of_colors_is_as_const}, one can show that $\cN_f^\red $ is $\cT $-measurable, where $\cT$ is the tail color $\sigma$-algebra. Hence $\cN_f^\red $ is $\mu$-almost surely constant. Denote this constant value by $N$. Claim \ref{claim:red_exists} implies that $N > 0$ while \Cref{prop:num_of_colors_constant} implies that $N < \infty$. We claim that in fact $N = 1$. 
    
    By Claim \ref{claim:equivariance claim}, the pushforward probability measure $\lambda = (P_f^\red)_* \mu$ on the space of pointed partitions $\cE^*(\NN)$ is $\Ai$-invariant. There are  $N$ non-distinguished blocks and one distinguished block. By \Cref{eq:color_distribution} and the the discussion preceding it we get $(R_f)_* \mu = \lambda_{\alpha^r}$ where the sequence $\alpha^r \in \cD$ is uniform on the first $q^r$ points. By the strong law of large numbers, for  $\mu$-almost every subgroup the asymptotic density of each  color in the basis $f$ is equal to $q^{-r(\mu)}$. Group the $N$ colors of $\Col_f^\red(H)$ as the non-distinguished blocks, and the rest of the colors into the distinguished block. For $\lambda$-almost every pointed partition $P \in \cE^*(\NN)$ the density vector $D(P)$ considered in  \Cref{lemma:density_determines_ergodic} is equal to the deterministic sequence $\alpha \in \cD^*$ given by
    \[
        \alpha_n = \begin{cases}
            q^{-r(\mu)}, & 1 \leq n \leq N, \\
            1 - N q^{-r(\mu)}, & n = -1, \\
            0, & \text{otherwise}.
        \end{cases}
    \]
    By \Cref{lemma:density_determines_ergodic} we have $\lambda = \lambda^*_\alpha$. In particular $\lambda$ is ergodic.
    By \Cref{prop:sln_probability} we get for any $n \in \NN$ that
    \begin{align*}
        q^{-rn} (1 - o(1)) &= \mu {\left( \set{\SL_n(f) \leq H} \right)}
        = \mu {\left( \set{\text{$1, \ldots, n$ have the same color from $\Col_f^\red(H)$}} \right)} 
        \\&= \lambda_\alpha^*{\left( \set{1, \ldots, n \textrm{ belong to the same non-distinguished block}} \right)} = N q^{-r n}.
    \end{align*}
    This implies $N = 1$, thereby completing the proof of the first part of the statement. 

    Lastly, let $f$ and $f'$ be a pair of finitary bases. Given a $\mu$-generic subgroup $H$, let $c \in \Col_f(H)$ and $c' \in \Col_{f'}(H)$ be the  unique colors so that $\SL_c(f) \le H$ and $\SL_{c'}(f') \le H$. There is a sufficiently large $n_0 \in \NN$ such that $f^n = (f')^n$ for all $ n > n_0$. This implies that 
    \[\SL_{c \cap \NN_{>n_0}}(f') = \SL_{c \cap \NN_{>n_0}}(f) \le \SL_{c }(f) \le H.\]
    It follows from Claim \ref{prop:if_sl_n_then_colored} and the uniqueness of the color $c'$ that the infinite intersection $c \cap \NN_{>n_0}$ is contained in $c'$. Therefore $c$ and $c'$ represent the same color in the set $\Col(H)$.
\end{proof}

\subsection{ $\GLj{n}$-action  on color vectors}
\label{subsection:action_on_coloring}

Let $\mu \in \IRS(\SLi)$ be a fixed non-trivial ergodic invariant random subgroup.
Recall that the distribution of the coloring relation is $(R_e)_* \mu = \lambda_{r(\mu)}$. Hence for $\mu$-almost every subgroup $H \in \Sub{\SLi}$ and for each $ n \in \NN$, every point in $\Col(H)^n$  appears as the tuple $\col_e(H)_I$ for an infinite number of subsets $I \subseteq \NN $ with $|I| = n$ and with $\min I$ being arbitrary large (this index subscript notation was introduced in \Cref{section:analysis_of_action}).





Recall that the color vector $\col_f(H) \in \Col(H)^\NN$  is \textit{not} the same between different finitary bases. We define an action of $\GLi$ on the product space $\Col(H)^\NN$. This action will have the property that if $f = A^{-1} e$ for some matrix $A \in \GLi$ then $\col_f(H) = A \cdot \col_e(H)$.

\begin{definition}
[Action on color vectors]
\label{def:action_on_coloring}
    Let $H$ be a $\mu$-generic subgroup. We define an action $\Theta_H$ of the group $\GLn$ on the product space $\Col(H)^n$ as follows.
    Let $A \in \GLn$ be a matrix and $\gamma \in \Col(H)^n$ be a color tuple.
    There exists a subset $I \subset \NN$ with $|I| = n$  such that $\gamma = \col_e(H)_{I}$.
    Let $f$ be the finitary basis given by $f = [A^{-1}]_I \, e$. 
    We set
    \[
    \Theta_H(A) \gamma = \col_f(H)_I \in \Col(H)^n.
    \]
When the context is clear, we shall suppress the notation $\Theta_H$ and write simply $A \cdot \gamma$ for $\Theta_H(A)\gamma$.
\end{definition}

The basic properties of the action $\Theta_H$ are established in the following result.

\begin{proposition} \label{prop:action_defined}
    Let $H$ be a $\mu$-generic subgroup. The actions $\Theta_H$ on the groups $\GLn$ on the product spaces $\Col(H)^n$ for each $n \in \NN$ have the following properties:
    \begin{enumerate}
        \item $\Theta_H$ is well-defined, i.e.\ its definition does not depend on the choice of subset $I$.
            \item $\Theta_H$ gives rise to a group action $\GLn \curvearrowright \Col(H)^n$.

         \item Fix $n_1, n_2 \in \NN$. Suppose $\gamma_1 \in \Col(H)^{n_1}$ and $\gamma_2 \in \Col(H)^{n_2}$. Let $A_1 \in \GLj{n_1}$ and $A_2 \in \GLj{n_2}$. Then
        \[
           \Theta_H \begin{pmatrix} A_1 & 0 \\ 0 & A_2 \end{pmatrix} (\gamma_1, \gamma_2) = (A_1 \gamma_1, A_2 \gamma_2).
        \]
        
    \end{enumerate}
\end{proposition}

The proof of \Cref{prop:action_defined} will require the following  \Cref{lemma:equal_colors_iff}, which is the basic tool to show that two color vectors are equal. It uses the notion of extension by permutations and of the subgroup $H^* \le \GLi$ associated to any given subgroup $H \le \GLi$ as in \Cref{sub:extension by permutations}. In most situations where \Cref{lemma:equal_colors_iff} is invoked, we will only need conditions (1) and (2). The more complicated condition (3) will be seldom used.
    
\begin{lemma}\label{lemma:equal_colors_iff}
    Let $H$ be a $\mu$-generic subgroup.
    Let $f$ be the finitary basis given by $f = A \cdot e$ for some matrix $A \in \GLi$. Let $I, J \subseteq \NN$ be two index sets with $|I|=|J|=n$ and with $\max I < \min J$. 
    Consider the matrix
    \[
    \Omega = \begin{pmatrix}
                    0 & I_n  \\ I_n & 0        \end{pmatrix} \in \SL_{2n}(\FFq).
    \]
    The following three conditions are equivalent:
    \begin{enumerate}
        \item $\col_f(H)_I = \col_f(H)_J$.
        \item $A \left[\Omega\right]_{I \union J} A^{-1} \in H^*$.
        \item Let $\omega \in \GLj{2}$ be the  matrix $\begin{psmallmatrix}
        1 & 0 \\ 0 & 1
    \end{psmallmatrix}$ in case  $n$ is even, and $\begin{psmallmatrix}
        0 & 1 \\ 1 & 0
    \end{psmallmatrix}$ in  case  $n$ is odd. For every index $m > \max J$ with $A \in \GLm$ and \textit{every} pair of indices $m < i < j$ with $\col_f(H)_i = \col_f(H)_j$ we have
        \[
        A \left[\begin{pmatrix}
            \Omega & 0 \\ 0 & \omega
        \end{pmatrix}\right]_{I \union J \union \set{i, j}} A^{-1} \in H.
        \]
    \end{enumerate}
\end{lemma}
\begin{proof}
Condition (2) implies condition (1) by the definition of the subgroup $H^*$ and by \Cref{prop:color_basic_definitions}.

Assume condition (1). Let $m > \max J$ and $m < i < j$ be indices as in condition (3).
The matrix appearing in condition (3) is a permutation matrix in the basis $f$ with determinant equal to $1$, that is a matrix in $\Ai(f)$.
It is a product of transpositions matching elements of $I$ and of $J$, together with the transposition $(i,j)$ in case $n$ is odd. These pair have the same color by assumption (1) and by the choice of the pair $i,j$.
By \Cref{prop:color_basic_definitions}, this matrix belongs to $H$, which shows condition (3).

Assume condition (3). Let $m > \max J$ be a sufficiently large index with $A \in \GLm$. There is some pair of indices $m < i < j$ with $\col_f(H)_i = \col_f(H)_j$.
As $A \in \GLm$, the matrix appearing in (3) is an extension by permutation of $A \left[\Omega\right]_{I \union J} A^{-1}$, regarded as a matrix in $\GLm$, by the permutation matrix $\omega$ supported on the pair $\set{i, j}$.
As this extension by permutation belongs to $H$, condition (2) follows by the definition of $H^*$.
\end{proof}

\begin{proof}[Proof of \Cref{prop:action_defined}]
    We start by showing item (1). Let $I$ and $I'$ be two index sets of cardinality $n$ such that $\gamma = \col_e(H)_I = \col_e(H)_{I'}$. Define the two finitary bases $f = [A^{-1}]_I e$ and $f' = [A^{-1}]_{I'} e$. We wish to show that $\col_{f}(H)_I = \col_{f'}(H)_{I'}$.
    Since by assumption $\gamma = \col_e(H)_{I''}$ for infinitely many sets of indices $I''$ with the index $\min I''$ being  arbitrary large, we may assume without loss of generality that $I \inter I' = \varnothing$ and that $\max I < \min I'$.

    Let $\delta = \col_f(H)_{I}$. Find an index set $J \subseteq \NN$ with $\min J > \max I'$ and of cardinality $n$ such that $\delta = \col_e(H)_J$. Notice that $\delta = \col_f(H)_J = \col_{f'}(H)_J$ by the natural identification of colors in different bases. Apply \Cref{lemma:equal_colors_iff} with respect to the two equalities of pairs of color vectors
    \[
        \delta = \col_f(H)_I = \col_f(H)_J
    \] and
    \[
        \gamma = \col_e(H)_I = \col_e(H)_{I'}.
    \]
    The first equality gives
    \[
    B_1 = [A^{-1}]_I \left[\Omega\right]_{I \union J} [A]_I \in H^*
    \] 
    and the second 
    \[
    C = \left[\Omega\right]_{I \union I' } \in H^*.
    \]
    What we are trying to show is $\delta = \col_{f'}(H)_{I'}$. This amounts to demonstrating that 
    \[
    B_2 = [A^{-1}]_{I'} \left[\Omega\right]_{I' \union J} [A]_{I'} \in H^*.
    \]
    One may check directly the matrix equation $B_2 = B_1 C B_1^{-1}$. It follows that $B_2 \in H^*$. We conclude by the converse direction of \Cref{lemma:equal_colors_iff}. 
    
    We next show item (3). Let $I_1, I_2 \subseteq \NN$ be a pair of finite sequences of indices with $\gamma_1 = \col_e(H)_{I_1}$ and $\gamma_2 = \col_e(H)_{I_2}$. We may assume without loss of generality that $\max I_1 < \min I_2$. Thus
    \[
        A \cdot \gamma = \col_f(H)_{I_1 \union I_2},
    \]
    where $A = \begin{psmallmatrix} A_1 & 0 \\ 0 & A_2 \end{psmallmatrix}$, $\gamma = (\gamma_1, \gamma_2)$ and $f = [A^{-1}]_{I_1 \union I_2} \cdot e$.
    We will show that 
    \[
        A_1 \cdot \gamma_1 = \col_f(H)_{I_1} \quad\text{and} \quad A_2 \cdot \gamma_2 = \col_f(H)_{I_2}.
    \]

    Let $\delta_1 = A_1 \cdot \gamma_1$. By definition $\delta_1 = \col_{f_1}(H)_{I_1}$, where $f_1 = [A_1^{-1}]_{I_1} e$.
    Choose a set of indices $J \subseteq \NN$ with $\min J > \max I_2$ such that $\delta_1 = \col_e(H)_J$. Then also $\delta_1 = \col_{f_1}(H)_J = \col_f(H)_{J}$. We may now invoke \Cref{lemma:equal_colors_iff} with respect to the equality
    \[
        \col_{f_1}(H)_{I_1} = \col_{f_1}(H)_J.
    \]
    This gives
    \[
        [A_1^{-1}]_{I_1} [\Omega]_{I_1 \union J } [A_1]_{I_1} \in H^*. 
    \]
    The above matrix also equals
    \[
        [A^{-1}]_{I_1 \union I_2} [\Omega]_{I_1 \union J } [A]_{I_1 \union I_2}. 
    \]
    Hence
    \[
        \col_f(H)_{I_1} = \col_f(H)_J = \delta_1,
    \]
    as required. The proof for the counterpart action $A_2 \cdot \gamma_2$ is exactly the same. 
    
    For item (2), the group action property, let $\gamma \in \Col(H)^n$ be a color tuple and $A, B \in \GLn$ be a pair of matrices. 
    Define $\delta = B \cdot \gamma$. We want to show that $A \cdot \delta = AB \cdot \gamma$. Let $I, J \subseteq \NN$ with $\max I < \min J$ be finite sequences of indices such that $\gamma = \col_e(H)_I$ and $\delta = \col_e(H)_J$.
    Let $T = \begin{psmallmatrix} A B & 0\\ 0 & A \end{psmallmatrix}$ and  $f = [T^{-1}]_{I \union J} \cdot e$. By item (3), which was established above, we know that
    \[
        AB \cdot \gamma = \col_f(H)_I, \quad A \cdot \delta = \col_f(H)_J.
    \]
    It suffices to show that 
    \[
    [T^{-1}]_{I \union J} [\Omega]_{I \union J } [T]_{I \union J} \in H^*,
    \]
    where the matrix $\Omega$ is as in \Cref{lemma:equal_colors_iff}. However, it can be seen that the above matrix is equal to
    \[
    C=[B^{-1}]_I [\Omega]_{I \union J } [B]_{I}.
    \] 
    We have $C \in H^*$ by another application of the lemma and since $\delta = B \cdot \gamma$. 
\end{proof}

\subsection{$\GLi$-action on color vectors}

Given a  non-trivial ergodic invariant random subgroup $\mu \in \IRS(\SLi)$ and a $\mu$-generic subgroup $H$, we have defined  a collection of actions $\Theta_H :\GLn \curvearrowright \Col(H)^n$  for each $n \in \NN$ (\Cref{def:action_on_coloring}). This collection induces a group action $\Theta_H : \GLi \curvearrowright \Col(H)^\NN$. This action is local, in the sense of \Cref{def:permutational_action}.
Here are two cases where we can calculate the action.

\begin{proposition} \label{prop:action_properties} Let $H$ be a $\mu$-generic subgroup. 
The action $\Theta_H : \GLi \curvearrowright \Col(H)^\NN$ satisfies the following two properties:
    \begin{enumerate}
        \item The constant vector $\rho=(\red(H))_{n \in \NN} \in \Col(H)^\NN$ is a fixed point of the action.
        \item The action is permutational, in the sense of \Cref{def:permutational_action}.
    \end{enumerate}
\end{proposition}

\begin{proof}
    We first show statement (1). Fix $n \in \NN$. Let $A \in \GLn$ be any matrix. Let $\rho_n \in \Col(H)^n$ be the constant $\red$ vector, namely $(\rho_n)_i = \red(H)$ for all $i \in \set{1,\ldots,n}$.
    Let $I \subset \NN$ be a subset of indices with $|I|=n$ such that $\col_e(H)_I = \rho_n$.
    To show that $A \cdot \rho_n = \col_{[A^{-1}]_I e}(H)_I$ is equal to $\rho_n$, it suffices by \Cref{prop:red_exists_and_unique} to show that $\SL_I([A^{-1}]_I e) \leq H$. But this subgroup is equal to $[A^{-1}]_I \SL_I(e) [A]_I = \SL_I(e)$, which is contained in $H$ by the choice of $I$.
        
    We now show statement (2). 
    Let $\sigma \in S_\infty$ be a permutation and  $P_\sigma \in \GLi$ be the corresponding  permutation matrix in the standard basis. We need to show that  \[
            (P_\sigma \cdot \gamma)_i = \gamma_{\sigma(i)} \quad \forall \gamma \in \Col(H)^\NN, \; i \in \NN.
        \]
    Fix some $n \in \NN$  such that $\sigma \in S_n$. Let $\delta = P_\sigma \cdot \gamma$. 
    Let $I, J \subseteq \NN$ be two finite sequences of indices with $|I|=|J|=n$  such that $\max I < \min J$ and $\gamma = \col_e(H)_I$, $\delta = \col_e(H)_J$. 
    It follows that $\delta = \col_f(H)_J$ with $f = [P_\sigma^{-1}]_I \cdot e$. 
    Apply \Cref{lemma:equal_colors_iff} with respect to the equation $\col_f(H)_I = \col_f(H)_J$. Let $\Omega$ be the matrix as in that Lemma. We get
    \[
    \left[ \begin{pmatrix} 0 & P_\sigma^{-1} \\ P_\sigma & 0  \end{pmatrix} \right]_{I \union J }
    = [P_\sigma^{-1}]_I [\Omega]_{I \union J} [P_\sigma]_I \in H^*.
    \]
    This is a permutation matrix. Therefore \Cref{prop:color_basic_definitions} implies that $\delta_k = \gamma_{\sigma(k)}$, as required.
\end{proof}

We now wish to relate the action $\Theta_H$ associated with a $\mu$-generic  subgroup $H$ to the action associated with its $\SLi$-conjugates.
Let $g \in \SLi$ be some element.
Since the element $g$ is finitary, the partitions $\Col(H)$ and $\Col(g H g^{-1})$ differ only on a finite number of elements.
This defines a natural identification of the colors in $\Col(H)$ with the colors in $\Col(g H g^{-1})$.

\begin{proposition} \label{prop:conjugate_action}
    Let $g \in \SLi$ be any element. For every $\mu$-generic subgroup $H$  we have that
    \[
        \Theta_H(g) \col(H) = \col(g H g^{-1}).
    \] 
\end{proposition}

\begin{proof}
    Let $n \in \NN$ be such that $g \in \SLn$. 
    Denote $I = \set{1, \ldots, n}$, $\gamma = \col(H)_I$ and $\delta = \col(g H g^{-1})_I$.
    Choose a finite subset of indices $J \subseteq \NN$ with $\min J > \max I$ and $\delta = \col(H)_J$. By the above identification of the colors in $\Col(H)$ with the colors in $\Col(g H g^{-1})$, we also get $\delta = \col(g H g^{-1})_J$. This means, using the  terminology of \Cref{lemma:equal_colors_iff}, that
    \[
        [\Omega]_{I \union J } \in (g H g^{-1})^* = g H^* g^{-1}.
    \]
    This is equivalent to 
    \[
        g^{-1} [\Omega]_{I \union J} g \in H^*,
    \]
    which shows that $\Theta_H(g) \gamma = \delta$.  
\end{proof}

\subsection{Fixed points for the $\GLn$-actions}
There is one  last property of the action $\Theta_H$ to be established in this section.

\begin{proposition} \label{prop:no_fixed_points}
    The action $\Theta_H$ associated to $\mu$-almost every generic subgroup $H $ has the following property:
    Fix $n \in \NN$ and suppose $A \in \GLn$ is such that $1$ is not an eigenvalue of $A$ (equivalently, $\rkt A = n$).
    Then the only fixed point of the action of $A$ on $\Col(H)^n$ is the constant $\red$ vector, that is, the point $\rho_n=(\red(H))_{i=1}^n$.
\end{proposition}

Before proceeding with the proof of the above, we note that it is not apriori clear that the event described  is measurable. We start by showing something stronger.

\begin{claim} 
\label{claim:measurable fixed point count}
Fix $n\in\NN$. Given a matrix $A \in \GLn$ let $\Phi_n^A : \Sub{\SLi} \to \NN$ be the $\mu$-measurable function that counts for a $\mu$-generic subgroup $H$ the number of fixed points of $A$ in its action  on $\Col(H)^n$ via $\Theta_H$. Then $\Phi_n^A$ is $\cT$-measurable, where $\cT$ is the tail color $\sigma$-algebra of $\mu$.
\end{claim}

The tail color $\sigma$-algebra was introduced in \Cref{def:tail_color_sigma_algebra}. Since the  vector $\rho_n$ is always a fixed point, \Cref{prop:no_fixed_points} refers to the measure of the event $\setcomp{H}{\Phi^A_n(H) = 1}$, which is seen to be measurable by Claim \ref{claim:measurable fixed point count}.

\begin{proof}[Proof of Claim \ref{claim:measurable fixed point count}]
    For each  $m \in \NN$, we shall show that $\Phi_n^A$ is $\cT_m$-measurable. Let $I, J \subseteq \NN_{> m}$ be two index sets with $|I|=|J|=n$ and $\max I < \min J$. 
    Given a matrix $B \in \GLn$ consider the event
    \[
        E^B_{I, J} = \setcomp{H  }{ \Theta_H(B)\col_e(H)_I = \col_e(H)_J}.
    \]
    Making use of \Cref{lemma:equal_colors_iff} and the notations therein, this event can be written as
    \[
        E^B_{I, J} = \setcomp{H  } 
        { [B^{-1}]_I [\Omega]_{I \union J} [B]_I \in H^*}.
    \] 
    This equation and statement (3) of \Cref{lemma:equal_colors_iff} show that for any such pair $I$ and $J$, the event $E^B_{I, J}$ is $\cT_m$-measurable. 
    Substituting the identity matrix $I_n$ for $B$, we define the event
    \[
        C_{I, J} = \setcomp{H  }{\col_e(H)_I = \col_e(H)_J} = E^{I_n}_{I, J} \in \cT_m,
    \]
    as well as the event
    \[
        F^A_I = \setcomp{H  }{A \cdot \col_e(H)_I = \col_e(H)_I}.
    \]
    The latter event $F^A_I$ is also $\cT_m$-measurable, since
    \[
        F^A_I = \setcomp{H  }{\exists J \subseteq \NN_{> m}, \abs{J} = n, \min J > \max I, H \in E_{I, J}^A \inter C_{I, J}} \in \cT_m.
    \] 
    Finally, the event $\setcomp{H}{\Phi_n^A \geq k}$ consists of the subgroups $H$ for which there are pairwise disjoint subsets of indices $I_1, \ldots, I_k \subseteq \NN_{> m}$ of cardinality $n$, with $\max I_s < \min I_{s + 1}$ for all $1 \leq s < k$, such that $H \in F^A_{I_s}$ for all $1 \leq s \leq k$ but $H \notin C_{I_s, I_t}$ for all $1 \leq s < t \leq k$.
    This shows that $\Phi_n^A$ is measurable  with respect to $\cT_m$. As $m \in \NN$ was arbitrary, $\Phi_n^A$ is also $\cT$-measurable.
\end{proof}

\begin{proof}[Proof of \Cref{prop:no_fixed_points}]
    Fix $n \in \NN$ and consider a matrix $A \in \GLn$. Let $H \in \Sub{\SLi}$ be a fixed subgroup. 
    By Claim \ref{claim:measurable fixed point count}, the cardinality of the set $\Fix_n(A)$ of fixed points of $A$ in its action on $\Col(H)^n$ via $\Theta_H$ is $\mu$-a.s.\ constant. 
    We will show the contrapositive of the desired statement, that is,  if this constant cardinality is at least $2$ then the matrix $A$ has $1$ as an eigenvalue.

    Suppose $\gamma \in \Col(H)^n$ is a fixed point of $A$. If $\gamma = \col_e(H)_{\{1, \ldots, n\}} = \col_e(H)_{\{n+1, \ldots, 2n\}}$ then \Cref{lemma:equal_colors_iff} implies in particular that
    \[
    [\Omega]_{\set{1, \ldots, 2n} } \in H^* 
    \quad \text{and} \quad
    [A^{-1}]_\set{1, \ldots, n} [\Omega]_{\set{1, \ldots, 2n} } [A]_\set{1, \ldots, n} \in H^* 
    .
    \]
    Denote by $B$ the product of these two matrices, so that
    \[B = \left[\begin{pmatrix}
        A^{-1} & 0 \\ 0 & A
    \end{pmatrix}\right]_{\set{1, \ldots, 2n}}.
    \]
    We have $B \in H$ (rather than just $B \in H^*$). Indeed, apply statement (3) of \Cref{lemma:equal_colors_iff} twice with the same pair of indices  of the same color $2n < i < j$.
    The two resulting matrices in $H$ have the same block $\omega$ in the $\set{i,j}$ position. This block cancels when taking their product, since $\omega^2 = I_2$.

    The probability that $B \in H$ is $q^{-r \cdot \rkt B}=q^{-2r \cdot \rkt A}$, where $r = r(\mu)$ is the degree of $\mu$. On the other hand, this probability is bounded from below by the probability that $\delta = \col_e(H)_{\{1, \ldots, n\}} = \col_e(H)_{\{n+1, \ldots, 2n\}}$ for some vector $\delta \in \Fix_n(A)$. Intuitively speaking, we expect the latter probability to be $\abs{\Fix_n(A)}q^{-2r n}$. Knowing this would imply the desired contrapositive statement.
    To compute the probability in question in a precise manner, we follow a similar route to the proof of the uniqueness of the $\red$ color in \Cref{subsection:the_color_red}.

    We convert the coloring partition $\Col(H)$ into another partition. Denote 
    \[\NN^{[n]} = \setcomp{(n(m-1)+1,\ldots,nm)}{m \in \NN}.\]
    We define an equivalence relation $L(H)$ on $\NN^{[n]}$ by setting $(b_1,b_2) \in L(H)$ if and only if 
    \[
        \col_e(H)_{b_1} = \col_e(H)_{b_2}.
    \]
    Each tuple of colors $\gamma \in \Col(H)^n$ determines a corresponding equivalence class of $L(H)$, which we denote by $c_\gamma$.
    From the relation $L(H)$ we form a pointed partition of $\NN^{[n]}$, to be denoted $P^A(H)$ (recall this notion from \Cref{subsubsection:pointed_partitions}).
    The blocks of $P^A(H)$ are 
    \[
        P^A(H) = \setcomp{c_\gamma}{\gamma \in \Fix_n(A)} \union \set{ \NN^{[n]} \setminus \bigunion_{\gamma \in \Fix_n(A)} c_\gamma},
    \]
    with the distinguished block being $\NN^{[n]} \setminus \bigunion_{\gamma \in \Fix_n(A)} c_\gamma$.
    The group $A_\infty(\NN^{[n]})$ of all finitary even permutations of the set $\NN^{[n]}$ naturally embeds in the group $\SLi$ as permutation matrices which permute blocks of size $n$ in the standard basis.
    The map $P^A : \Sub{\SLi} \to \cE^*(\NN^{[n]})$ is $\mu$-measurable and  $A_\infty(\NN^{[n]})$-equivariant, where $\cE^*(\NN^{[n]})$ is the set of all pointed partitions on $\NN^{[n]}$.
    
    Set $\lambda = P^A_* \mu$. The probability measure $\lambda$ is $A_\infty(\NN^{[n]})$-invariant. Recall that $\alpha^r \in \cD$ denotes the uniform distribution on $q^r$ points. Since $(R_e)_* \mu = \lambda_{\alpha^r}$ and as was observed in the proof of \Cref{prop:red_exists_and_unique} by the strong law of large numbers, $\mu$-a.s. every tuple $\gamma \in \Col(H)^n$ occurs as the color of a block with asymptotic density exactly $q^{-rn}$. Note that $N = \abs{\Fix_n(A)}$ is $\mu$-a.s.\ constant by Claim \ref{claim:measurable fixed point count}. We group the $N$-many color tuples in $\Fix_n(A)$ into non-distinguished blocks and the rest of the color tuples in $\Fix_n(A)$ into the distinguished block. This results in asymptotic density $q^{-rn}$ for each non-distinguished block, and asymptotic density $1 - N q^{-rn}$ for the distinguished one. Thus the asymptotic density vector $D(P^A(H))$ is $\mu$-a.s.\ equal to the deterministic vector $\alpha = (\alpha_m)_{m \in \ZZ_{\geq -1}}$ given by
    \[
        \alpha_m = \begin{cases}
            q^{-r n}, & 1 \leq m \leq \abs{\Fix_n(A)}, \\
            1 - \abs{\Fix_n(A)} q^{-r n}, & m = -1, \\
            0, & \text{otherwise}.
        \end{cases}
    \]
    By \Cref{lemma:density_determines_ergodic} we deduce that $\lambda = \lambda_\alpha^*$ and  that $\lambda$ is ergodic.
    To conclude, we make use of the fact that $P^A_* \mu=\lambda^*_\alpha $ and compute 
    \begin{align*}
        q^{-r \cdot \rkt B} &= \mu \left( \setcomp{H}{B \in H}  \right)
         \geq
        \mu \left(\setcomp{H}{\substack{\col_e(H)_{\{1, \ldots, n\}} = \col_e(H)_{\{n+1, \ldots, 2n\}}=\delta \\ \text{for some  $\delta \in \Fix_n(A)$}}}\right)
        \\ & = \lambda^*_\alpha \left( \set{\substack{\text{$(1,\ldots,n),(n+1,\ldots,2n) \in \NN^{[n]}$  are in} \\ \text{the same non-distinguished block of $P^A(H)$}}} \right)
        = \abs{\Fix_n(A)} \cdot q^{-2 r n}.        
    \end{align*}
    Since $\rkt B = 2 \rkt A$ we get from this inequality that $\abs{\Fix_n(A)} \leq q^{2 r (n - \rkt A)}$. Hence, the assumption $\abs{\Fix_n(A)} > 1$ implies that $\rkt A < n$, so that $1$ is an eigenvalue of $A$.
\end{proof}

\subsection{Summary}

Given a non-trivial ergodic invariant random subgroup $\mu$ of $\SLi$, we have constructed for a $\mu$-generic subgroup $H$ an action of the finitary group $\GLi$ on the product set $\Col(H)^\NN$ of color vectors, where $\Col(H)$ is the set of colors of the subgroup $H$ in any finitary basis. 
The set $\Col(H)$ is of size $q^r$ where $r = r(\mu)$ is the degree of the invariant random subgroup $\mu$. This action of $\GLi$ is local, permutational and enjoys the properties established in Propositions \ref{prop:action_properties} and \ref{prop:no_fixed_points}. 
These are precisely the necessary assumptions for our abstract \Cref{thm:action_classification} giving a classification of $\GLi$-actions. This will be the point of departure in our proof of \Cref{thm:main_theorem} in \Cref{section:proof_end}.

\section{Invariant measures on compact abelian groups}
\label{section:measures_on_v}

Fix a pair $k,\l \in \NN_0$. Recall from \Cref{sec:the group SLi} that $V^{k,\l}$ denotes the product space $V^k \times \widehat V^\l$ with the product of the standard and the dual actions of $\GLi$.
Let $X^{k, \l} \subseteq V^{k,\l}$ be the subset of all tuples $(v^1, \ldots, v^k), (\p^1, \ldots, \p^\l) \in V^k \x \hV^\l$ such that the collection $\set{v^1, \ldots, v^k}$ is linearly independent in $V$ and the collection $\set{\p^1, \ldots, \p^\l}$ is  linearly independent in $\widehat V$. Note that $X^{k,\l}$ is $\GLi$-invariant and co-null in $V^{k,\l}$. By a slight abuse of notation, we let $\Haar^{k,\l}$ be the restriction of the $\GLi$-invariant Haar measure of $V^{k, \l}$ to $X^{k, \l}$. 

Recall that a Borel action of a group $G$ on a Borel space  is called \textit{uniquely ergodic} if it admits exactly one $G$-invariant probability Borel measure.
In this section, we prove \Cref{thm:xkl_measures} which essentially states that $\SLi$ acts uniquely ergodically on $X^{k, \l}$, with the unique invariant measure being $\Haar^{k,\l}$.
Additionally, we classify all the $\SLi$-invariant measures on $V^{k, \l}$ and on the infinite Grassmannian over $\FFq$ (see \Cref{subsection:invariant_measures_on_related_spaces}). 

\subsection{Ergodic theory and dynamics}

We employ the pointwise ergodic theorem for amenable groups. 
The most general case was established by Lindenstrauss in \cite{lindenstrauss2001ergodic_theorem}. However, our specific setting involves the locally finite group $G=\SLi$, which can be expressed as the increasing union $G=\bigunion_{n=1}^\infty G_n$ of the finite groups $G_n = \SLn$. In this setting, the required result was established considerably earlier by Vershik \cite{vershik1974ergodic_thm}.
The version of the theorem that we will use states the following.
\begin{theorem}\label{thm:pointwise_ergodic}
    Let $G$ be a locally finite countable group, written as an increasing union of finite subgroups $G = \bigunion_{n=1}^\infty G_n$.
    Suppose $G$ acts ergodically on a Borel probability space $(X, \mu)$. 
    Let $E \subseteq X$ be a measurable subset. Then  $\mu$-a.e.\ point $x \in X$ 
    satisfies
    \[
        \mu(E) = \lim_{n \to \infty} \frac{1}{\abs{G_n}} \abs{\setcomp{g \in G_n}{g \cdot x \in E}}.
    \]
\end{theorem}
    
We use this theorem to reduce the question of classifying invariant measures 
to that of counting matrices with specific properties.
In \Cref{subsection:measures_xk0}, we use the ergodic theorem to show that $\SLi$ acts uniquely ergodically on both $X^{k, 0}$ and $X^{0, \l}$. In \Cref{subsection:joins_of_measured_systems} we reduce the problem of showing that $\SLi$ acts uniquely ergodically on $X^{k, \l} = X^{k, 0} \x X^{0, \l}$ to the problem of showing that for $\Haar^{0, \l}$-a.e.\ element $\up \in X^{0, \l}$ the stabilizer $\Stab_\SLi(\up)$ acts uniquely ergodically on $X^{k, 0}$.
The proof of the latter fact can be found in \Cref{subsection:measures_xkl}, where we invoke  the ergodic theorem for a second time. 
Lastly, in \Cref{subsection:invariant_measures_on_related_spaces}, we classify the ergodic measures on $V^{k, \l}$ and on the Grassmannian space $\Gr_k(V) \x \Gr_\l(\hV)$, to be defined precisely there.

\subsubsection{Joins of measured systems}
\label{subsection:joins_of_measured_systems}
We briefly recall some notions in the ergodic theory of dynamical systems. 
A reference for these definitions is \cite[Chapter 6]{glasner2003ergodic}.

Let $G$ be a countable group which acts on two Borel probability spaces $(X, \mu)$ and $(Y, \nu)$ by measure-preserving transformations.
A \textit{joining} of $(X, \mu)$ and $(Y, \nu)$ is a $G$-invariant measure $\lambda$ on $X \x Y$
that has $\mu$ and $\nu$ as its marginals; that is, $\proj_X(\lambda) = \mu$ and $\proj_Y(\lambda) = \nu$.
The two systems $(X, \mu)$ and $(Y, \nu)$ are called \textit{disjoint} if the only joining of the two systems is the product measure $\mu \x \nu$.

\begin{lemma}\label{lemma:disjoint_if_stabilizer_uniquely_ergodic}
    Let $G$ be a countable group which acts on two Borel probability spaces $(X, \mu)$ and $(Y, \nu)$ by measure-preserving transformations.
    Suppose that for $\mu$-a.e.\ point $x \in X$, the stabilizer $\Stab_G(x)$ is uniquely ergodic
    in its action on $Y$. Then $(X, \mu)$ and $(Y, \nu)$ are disjoint.
\end{lemma}

\begin{proof}
    Let $\lambda$ be a joining of $(X, \mu)$ and $(Y, \nu)$. Write $G_x = \Stab_G(x)$ for all $x \in X$.
    Disintegrate the measure $\lambda$ with respect to the projection map $\proj_X$. This gives a measurable family of probability measures $\set{\lambda_x}_{x \in X}$ on $Y$ such that
    \[
        \lambda = \int_X \left( \delta_x \x \lambda_x \right) \diff \mu(x).
    \]
    We claim that $\mu$-a.e.\  measure $\lambda_x$ is $G_x$-invariant.
    Indeed, for any $g \in G$ we have
    \[
        \lambda = g_* \lambda = \int_X \left( \delta_{g  x} \x g_*\lambda_x \right) \diff \mu(x).
    \]
    By the uniqueness of the disintegration, we have $g_* \lambda_x = \lambda_{g  x}$ for $\mu$-a.e.\ $x \in X$.
    As $G$ is countable, $\mu$-a.e.\ point $x \in X$ satisfies for all $g \in G$ that $g_* \lambda_x = \lambda_{g x}$.
    For all such points $x \in X$ we conclude that the measure $\lambda_x$ is $G_x$-invariant, as claimed.
    Since $G_x \curvearrowright Y$ is $\mu$-a.s.\ uniquely ergodic and $\nu$ is $G_x$-invariant (as it is $G$-invariant), the unique $G_x$-invariant measure on $Y$ must be $\nu$.
    Thus $\lambda_x = \nu$ holds $\mu$-a.s. So $\lambda = \mu \x \nu$, as required.
\end{proof}

\subsection{Invariant measures on \texorpdfstring{$X^{k, 0}$}{X\textasciicircum k,0}}
\label{subsection:measures_xk0}
Fix $k \in \NN$ (the index $\l$ will not used in this  subsection). We show that the actions $\SLi \curvearrowright X^{k, 0}$ and $\SLi \curvearrowright X^{0, k}$ are uniquely ergodic, with the unique invariant measure being the restriction of the corresponding Haar measure from $V^k$ and $\widehat V^k$. 
Note that these two actions are related by the automorphism $g \mapsto (g^{-1})^t$ of the group $\SLi$, so it suffices to deal  only with $X^{k, 0}$.

Given $m \in \NN$ and a sequence of vectors $a^1,\ldots, a^k \in V_{m}$, we define the set 
\[
    \left[ a^1, \ldots, a^k \right] = \setcomp{ \uv \in X^{k, 0}}{  v^i_{\leq m} = a^i \quad \forall 1 \le i \le k}.
\]
We refer to a set of this form as a \textit{box} and to the integer $m$ as its \textit{length}.
The collection of all boxes is a $\pi$-system which generates the Borel $\sigma$-algebra on $X^{k, 0}$.

A straightforward application of the pointwise ergodic theorem (\Cref{thm:pointwise_ergodic}) allows us to classify the ergodic invariant measures on $X^{k, 0}$.

\begin{proposition} \label{prop:xk0_measures}
    The only $\SLi$-invariant measure on $X^{k, 0}$ and on $X^{0, \l}$ is the Haar measure.
\end{proposition}

\begin{remark}
In the next sections, we will show the stronger statement that an action of some specific subgroup of $\SLi$ on $X^{k, 0}$ is uniquely ergodic (see \Cref{prop:stabilizer_dual_uniquely_ergodic}). This gives \Cref{prop:xk0_measures} as an immediate corollary. 
However, we choose to prove \Cref{prop:xk0_measures} directly to demonstrate a simple use of the ergodic theorem.
\end{remark}

\begin{proof}[Proof of \Cref{prop:xk0_measures}]
    Let $\mu$ be an $\SLi$-invariant measure on $X^{k, 0}$. Making use of the ergodic decomposition, we may assume without loss of generality that $\mu$ is ergodic.
    Let $E = [a^1, \ldots, a^k]$ be some box of length $m$ where
    $a^i \in V_m$.
    Let $\uv = \left( v^1, \ldots, v^k \right) \in X^{k, 0}$ be any generic point such that
    \[
        \mu(E) = \lim_{n \to \infty} \frac{\abs{\setcomp{g \in \SLn}{\left( g v^1, \ldots, g v^k \right) \in E}}}{\abs{\SLn}},
    \] 
    as per \Cref{thm:pointwise_ergodic}.
    Let $n > m + k$ be a large enough index such that the collection of vectors $\set{v^1_{\leq n}, \ldots, v^k_{\leq n}}$ is linearly independent.
    We want to estimate the number of elements $g \in \SLn$ satisfying $(gv^i) \in E$, that is $(g v^i)_{\leq m} = a^i$ for all $1 \le i \le k$.
    Complete $v^1_{\leq n}, \ldots, v^k_{\leq n}$ to a basis $f = \set{v_{\leq n}^1, \ldots, v_{\leq n}^k, u^{k+1}, \ldots, u^n}$ of $V_n$.
    Representing the transformation $g$ as a matrix going from basis $f$ to the standard basis, 
    this is equivalent to counting the number of matrices $g \in \SLn$ with $g_{ji} = a^i_j$ for all $1\le i \leq m$ and $1 \le j \leq k$.
    
    This is a standard counting argument. We proceed as follows. 
    Count the number of possibilities for each column of the matrix $g$ sequentially.
    The first column must start with $a^1$ and be non-zero, 
    so it has either $q^{n-m} - 1$ or $q^{n-m}$ possibilities, depending respectively on whether $a^1$ is zero or not.
    The second column must start with $a_2$ and be linearly independent from the first column,
    so it has $q^{n-m}$ minus at most $q$ possibilities. The exact number depends on the dimension of the span of $a^1$ and $a^2$. 
    Continuing this way, we see that for each $i \le k$ the $i$-th column of $g$ has $q^{n-m} - \beta_i$ possibilities, 
    where $\beta_i \leq q^{i-1}$. Since $n - m > k$ we also notice that this    number of possibilities is never zero.
    After choosing the first $k$ columns, the span of those columns has dimension exactly $k$.
    Hence for indices $i$ satisfying $k < i < n$, the $i$-th column has $q^{n} - q^{i-1}$ possibilities.
    The last column must account for $\det g = 1$, and as such has only $\frac 1 {q-1} (q^n - q^{n-1})$ possibilities.
    In summary, we have by this counting argument and by \Cref{thm:pointwise_ergodic}:
    \begin{align*}
        \mu(E) &= \lim_{n \to \infty} \frac{
            \frac{1}{q-1} ( q^{n-m} - \beta_1 ) \cdots ( q^{n-m} - \beta_k ) (q^n - q^k) \cdots (q^n - q^{n-1})
        }{
            \frac 1 {q-1} (q^n - 1) \cdots (q^n-q^{k-1}) (q^n-q^k) \cdots (q^n - q^{n-1})
        } = \\
        &= \lim_{n \to \infty} q^{-m k} \left( \frac{1 - \beta_1 q^{m-n}}{1 - q^{-n}} \right) \cdots \left( \frac{1 - \beta_k q^{m-n}}{1 - q^{k-1-n}} \right)
        = q^{-mk} = \Haar^{k, 0}(E).
    \end{align*}
    Since the boxes comprise a $\pi$-system, we get $\mu = \Haar^{k, 0}$, as required.

    As was mentioned above, the space $X^{0, k}$ is equivariantly isomorphic  to $X^{k, 0}$ under the intertwining automorphism $g \mapsto g^{-t}$ of $\SLi$, hence we obtain the result for it too.
\end{proof}

\subsection{Regular linear independence}
\label{sub:regular independence}

For each vector $v \in V$ and $m,t \in \NN$, let $v_{(m, m + t]}$ denote the vector $(v_i)_{i = m+1}^{m+t} \in \FFq^t$. Since $\hV$ and $V$ have the same underlying set, we use the same notation for elements of $\hV$.

\begin{definition} \label{def:regularly_linearly_independent}
    Fix $r \in \NN$. Let $\ubar v \in V^r$ and $m \in \NN$. The smallest  $t \in \NN$ such that the sequence of vectors \[
    \uv^1_{(m, m + t]}, \ldots, \uv^r_{(m, m+t]} \in \FFq^t
    \]
    is linearly independent is called the $m$-th \textit{index of independence} of $\ubar v$ and will be denoted $I_m(\ubar v)$. If no such $t$ exists, we let $I_m(\ubar v) = \infty$.
    If $T \in \NN$ is such that
    \[
        \liminf_{m \to \infty} I_m(\ubar v) = T,
    \]
    then $\ubar v$ is called \textit{$T$-regularly linearly independent}.
    If there is some $T$ such that $\ubar v$ is $T$-regularly linearly independent, then $\ubar v$ will be called \textit{regularly linearly independent}.
\end{definition}

Notice that a generic element of $X^{0, \l}$ is  $\Haar^{0, \l}$-almost surely $\l$-regularly linearly independent.
We aim to show that, provided $k \leq \l$, the action of the stabilizer of a regularly linearly independent functional $\up \in X^{0,\l}$ on $X^{k, 0}$ is uniquely ergodic, the only invariant measure being $\Haar^{k, 0}$.


\begin{lemma}\label{lemma:strongly_linearly_independent_dual}
    Let $\p^1, \ldots, \p^r \in \hV_n$ be linearly independent functionals. Let $1 \leq m \leq n$. The functionals $\p^1_{> m}, \ldots, \p^r_{>m}$ are linearly independent if and only if the restriction of the projection map
    \[
        \Pi^n_m:V_n \to V_m, \quad v \mapsto v_{\leq m}
    \]
    to the joint kernel $\set{\p^1, \ldots, \p^r}^\perp \subseteq V_n$ is surjective.
\end{lemma}

\begin{proof} Follows from basic linear algebra.
\end{proof}

\begin{lemma}\label{lemma:transitivity_of_stabilizer_dual}
    Fix $k, \l \in \NN$. Let $\up \in X^{0, \l}$ be a functional and $m \in \NN$.
    Assume that the functionals $\up^1_{>m}, \ldots, \up^\l_{>m}$ are linearly independent. Define
    \[
        K = \max\set{I_m(\up), \l + 2 k}.
    \]
    Then for any $n > m + K$ the following  holds:
    If $\uv \in V_{n}^k$ is linearly independent and $\uw \in V_m^k$
    then there exists an element $g \in \Stab_\SLn(\up)$ such that $(g \uv)_{\leq m} = \uw$.
\end{lemma}
\begin{proof}
    Fix some $n > m + K$. Let $\uv \in V_n^k$ be linearly independent and $\uw \in V_m^k$ be arbitrary.
    We will find an element $g \in \Stab_\SLn(\up)$ such that $(g \uv)_{\leq m} = \uw$.

    For each $i = 1, \ldots, k$ let $a^i = \uv_{\leq m}^i - \uw^i \in V_m$.
    By \Cref{lemma:strongly_linearly_independent_dual}, we can lift the $a^i$ to elements $\tilde{a}^i$ of the joint kernel subspace $\set{\up^1, \ldots, \up^\l}^\perp \subseteq V_n$ through the projection $\Pi^n_m$.
    The fibers of $\Pi^n_m$ are affine subspaces of dimension $n - \l - m$.
    Since this dimension is larger than $2 k$, we can choose the $\tilde{a}^i$ such that the set
    $\set{\uv^i}_{i=1}^k \cup \set{\tilde{a}^i}_{i=1}^k$ is linearly independent.
    Extend this linearly independent set to a basis
    \begin{equation} \label{eq:primal_basis}
        \uv^1, \ldots, \uv^k, \tilde{a}^1, \ldots, \tilde{a}^k, f^{2k+1}, \ldots, f^n
        \tag{$*$}
    \end{equation}
    of $V_n$.
    Let $g \in \SLn$ be the linear transformation which sends $\uv^i$ to $\uv^i - \tilde{a}^i$ for $1 \leq i \leq k$ and keeps the other basis vectors fixed.
    By definition $(g \uv^i)_{\leq m} = \uv^i_{\leq m} - a^i = \uw^i$ and $g$ has determinant $1$.

    Showing that the element $g$ fixes the functional $\up$ is equivalent to showing that $g^{-1}$ fixes $\up$. The latter fact holds if and only if
    \[
        \angles{gv, \up^j}_n = \angles{v, \up^j}_n \quad \text{for all } v \in V_n \text{ and } j = 1, \ldots, \l.
    \]
    These equations are satisfied for each of the basis vectors given in (\ref{eq:primal_basis}) by definition.
\end{proof}

The following proposition is the main technical tool that will help us apply the ergodic theorem to $\Stab_\SLi(\up)$-ergodic measures. We use the linear algebraic notations introduced in \Cref{subsec:standard and dual}.

\begin{proposition}\label{prop:infinitely_often_full_rank}
    Let $k, \l \in \NN$. Let $\up \in X^{0, \l}$ be regularly linearly independent and $\mu$ be an ergodic $\Stab_\SLi(\up)$-invariant measure on $X^{k, 0}$. Then
    \[
        \mu{\left( \setcomp{\uv \in X^{k, 0}} {{M_n}(\uv, \up) \text{ has full rank for infinitely many } n} \right)} = 1.
    \]
\end{proposition}

\begin{proof}
    Define $P = \Stab_\SLi(\up)$ and $P_n = P \cap \SLn$ for each $n \in \NN$. Let $E$ denote the event in question.
    Since the pairing $\angles{\cdot, \cdot}_n$ is $P_m$-invariant whenever $n \geq m$, the event $E$ is $P$-invariant.
    Hence by ergodicity of $\mu$ it will suffice to show that $\mu(E) > 0$.

    Let $T \in \NN$ be such that $\up$ is $T$-regularly linearly independent, and $J \subseteq \NN$ be the set of indices $m \in \NN$ such that the vectors $\up^j_{(m, m + T]}$ for  $j \in \set{1, \ldots, \l}$ are linearly independent. By assumption $J$ is infinite.
    Write $E$ as a nested intersection $E = \biginter_{m \in J} S_m$, where
    \[
        S_m = \setcomp{\uv \in X^{k, 0}} {\text{$\exists n \geq m$ such that $ {M_n}(\uv, \up)$  has full rank}}.
    \]
    We will show that $\mu(S_m)$ is bounded from below by some $\delta > 0$ independent of $m$.

    A box $[a^1, \ldots, a^k] \subseteq X^{k, 0}$ of length $m$ will be called \emph{independent} if the vectors $a^1,\ldots,a^k$ are linearly independent in $V_m$.
    Notice that $X^{k, 0} = \bigcup \set{B \text{ an independent box}}$. Hence, for all large enough $m$, there are disjoint independent boxes of length $m$ whose union has measure at least $1/2$.
    Fix such an index $m \in J$ together with a corresponding collection of independent boxes $B_1, \ldots, B_K$. We may also assume that $m$ is large enough such that the sequence $\up_{\leq m}$ is linearly independent.

    Define $Y \subseteq V_m^k$ to be the set of all vectors $\uw \in V_m^k$ such that $M_m(\uw, \up)$ has full rank.
    The ratio $\abs{Y} / q^{m k}$ is the probability that a uniformly random matrix in $\M_{\l \x k}(\FFq)$ has full rank.
    This probability is some positive value $\gamma > 0$ which is independent of $m$.

    Now, let $n = m + C$ where $C = \max\set{T, \l + 2k} + 1$.
    Each box of length $m$ is the disjoint union of $q^{k C}$ subboxes of length $n$. For each $i = 1, \ldots, K$ choose some subbox $B_i' \subseteq B_i$ of length $n$ satisfying $\mu(B_i') \geq \mu(B_i) q^{-k C}$.
    Up to reordering, we may assume without loss of generality that $\mu(B_i')$ is a decreasing sequence.
    Choose an arbitrarily ordering $\uw_1, \ldots, \uw_{\abs{Y}}$ of the set $Y$. Consider the pairs $(B_1', \uw_1), \ldots, (B_N', \uw_N)$ with $N = \min(K, \abs{Y})$.
    Since $K \leq q^{m k}$ it follows that $N/K \geq \gamma$.

    For each pair $(B_i', \uw_i)$ where $i \leq N$ select an element $g_i \in P_n$ such that $(g_i \uv)_{\leq m} = \uw_i$ for all $\uv \in B_i'$.
    This is possible to do by \Cref{lemma:transitivity_of_stabilizer_dual}.
    The subsets $g_i B_i'$ for $i=1, \ldots, N$ are pairwise disjoint.
    By definition we have $g_i B_i' \subseteq S_m$. 
    This gives:
    \[
        \mu(S_m) \geq \sum_{i=1}^N \mu(g_i B'_i)
        \geq \frac{N}{K} \sum_{i=1}^K \mu(B_i')
        \geq \gamma q^{-k C} \sum_{i=1}^K \mu(B_i) \geq \frac 1 2 \gamma q^{-k C}.
    \]
    Since the right-hand side is independent of $m$, this lower bound holds uniformly for arbitrarily large indices $m \in J$, as required.
\end{proof}

\subsection{Invariant measures on $X^{k,\l}$}
\label{subsection:measures_xkl}
We turn to our second use of the ergodic theorem.

\begin{proposition}\label{prop:stabilizer_dual_uniquely_ergodic}
    Fix $k,\l \in \NN$.
    Let $\up \in X^{0, \l}$ be a regularly linearly independent functional. Denote $P = \Stab_\SLi(\up)$.
    The only $P$-invariant probability measure on $X^{k, 0}$ is $\Haar^{k, 0}$.
\end{proposition}

\begin{proof}
    First, we may and will assume without loss of generality that $k \le l$. Indeed, replacing $\up$ by a functional $\up' \in X^{0,\l'}$ for some $\l' \ge \l$ whose first $l$ entries coincide with $\up$ replaces the
    stabilizer $P$ by a subgroup $P'$.
    Any $P$-invariant measure is, of course, $P'$-invariant. Hence, proving the proposition for the smaller group $P'$ would imply it for $P$, provided the functional $\up'$ is taken to be regularly linearly indepedendent, which is always possible.

    Let $\mu$ be a $P$-invariant measure on $X^{k, 0}$. By the ergodic decomposition we may assume without loss of generality that $\mu$ is ergodic.
    Let $E = [\ua] \subseteq X^{k, 0}$ be a box of some length $m$, where $\ua \in V_m^k$.
    Let $\uv = \left( \uv^1, \ldots, \uv^k \right) \in X^{k, 0}$ be a generic point with respect to $E$, in the sense that
    \begin{equation}\label{eq:ergodic_thm_for_stabilizer}
        \mu(E) = \lim_{n \to \infty} \frac{1}{\abs{P_n}} \abs{\setcomp{g \in P_n}{g \uv \in E}},
    \end{equation}
    as per \Cref{thm:pointwise_ergodic}. Here $P_n = P \cap \SLn$ for each $n \in \NN$.

    Choose an index $n > m + \l$ large enough such that each of the two collections $(\up^j_{\leq n})_{j=1}^\l$ and $(\uv^i_{\leq n})_{i=1}^k$ is linearly independent.
    Let $f = (f^j)_{j=1}^n$ be a basis of $\hV_n$ extending $\up^1_{\leq n}, \ldots, \up^\l_{\leq n}$, so that $f^j = \up^j_{\leq n}$ for $j = 1, \ldots, \l$. Viewing the elements of $\hV_n$ as row vectors, let $F$ be the $n \x n$ matrix whose $j$-th row is $f^j$. With respect to the standard basis, any element $g \in P_n$ takes the form
    \begin{equation}\label{eq:g_expression}
        g = F^{-1} \begin{pmatrix} I_{\l} & 0_{\l, n-\l} \\ B & A \end{pmatrix} F,
    \end{equation}
    where $A \in \SLj{n-\l}$ and $B \in \M_{(n-\l) \x \l}(\FFq)$ are some matrices.
    We obtain a group isomorphism $\SLj{n-\l} \ltimes \Mj{(n-\l)}{\l} \cong P_n$. In what follows, we identify each pair $(A, B)$ with its image in $P_n$ under this isomorphism.

    We look for conditions on the element $g \in P_n$ such that $g \uv \in E$, meaning $(g \uv^i)_{\leq m} = a^i$ for $i = 1, \ldots, k$.
    This calculation is naturally expressed in matrix notation. 
    Defining $v \in \Mj{\infty}{k}$ and $a \in \Mj{m}{k}$ as the matrices whose columns are $\uv$ and $\ua$ respectively, and letting $(\cdot)_{\leq m}$ denote the first $m$ rows of a matrix, we seek elements $g \in P_n$ such that $(g v)_{\leq m} = a$.

    Denote $u = F v_{\leq n} \in \Mj{n}{k}$. Split this matrix into the first $\l$ rows and the remaining rows,  as $u = \begin{pmatrix}
        u_{\leq \l} \\ u_{> \l}
    \end{pmatrix}$. Since $f^j = \up^j_{\leq n}$ for all $j \leq \l$, the matrix $u_{\leq \l}$ is equal to the pairing $M_n(\uv, \up)$  by definition.
    Writing $g v_{\leq n}$ in coordinates as $\begin{pmatrix} a \\ c \end{pmatrix}$, with $a \in \Mj{m}{k}$ and $c \in \Mj{(n-m)}{k}$,  expression (\ref{eq:g_expression}) gives
    \begin{equation} \label{eq:stabilizer_equation}
        F \begin{pmatrix} a \\ c \end{pmatrix} =
        \begin{pmatrix} I_{\l} & 0_{\l, n-\l} \\ B & A \end{pmatrix} \begin{pmatrix} u_{\leq \l} \\ u_{>\l} \end{pmatrix}
        = \begin{pmatrix} u_{\leq \l} \\ B u_{\leq \l} + A u_{>\l} \end{pmatrix}.
    \end{equation}

    The first $\l$ rows of this equation give constraints on $c$ which are independent of $g$:
    \begin{equation} \label{eq:c_equation}
        M_{n-m}(c, \up_{> m}) = u_{\leq \l} - M_m(\ua, \up).
    \end{equation}
    As $\up$ is assumed to be regularly linearly independent, for sufficiently large $n$  the vectors $\up^j_{(m, n]}$, $j = 1, \ldots, \l$ are linearly independent.
    For such $n$ the space of all matrices $c$ which are solutions to \cref{eq:c_equation} is a certain affine subspace $C(\uv, \ua) \subseteq \Mj{(n-m)}{k}$ of codimension $\l k$.

    Given $c \in C(\uv, \ua)$, the last $n - \l$ rows of \cref{eq:stabilizer_equation} give an equation for $g$ of the form
    \begin{equation} \label{eq:last_rows_constraint}
        B u_{\leq \l} + A u_{> \l} = d(\uv, \ua, c),
    \end{equation}
    where $d(\uv, \ua, c) \in \Mj{(n-\l)}{k}$.
    Rearranging \cref{eq:last_rows_constraint}, for each given $A \in \SLj{n-\l}$ we get linear equations on the matrix $B$:
    \begin{equation} \label{eq:b_equation}
        B u_{\leq \l} = d(\uv, \ua, c) - A u_{> \l}.
    \end{equation}
    In this form, each row of $B$ is constrained by $k$ affine equations. 
    For each row, the coefficient vectors of the homogeneous part of these equations are the columns of $u_{\leq \l} = M_n(\uv, \up)$.

    By \Cref{prop:infinitely_often_full_rank}, we may assume that $\uv$ is such that for infinitely many indices $n$, the corrsponding matrix $u_{\leq \l}$ has full rank. Since the limit in \cref{eq:ergodic_thm_for_stabilizer} necessarily exists, and our aim is to calculate it, we may restrict ourselves only to such $n$.
    
    In this case, we get that for each $c \in C(\uv ,\ua)$ and $A \in \SLj{n-\l}$, the set of matrices $B$ satisfying \cref{eq:b_equation} is an affine subspace of $\Mj{(n-\l)}{\l}$ of codimension $(n-\l)k$. Here we use the assumption $k \leq \l$ so that $\rk u_{\leq \l} = k$. Notice that the matrix $c$ is determined by the pair $(A, B)$ via \cref{eq:stabilizer_equation}.
    Thus we obtain the count
    \begin{align}\label{eq:g_count_good_n} 
    \begin{split}
        \abs{\setcomp{g \in P_n}{g \cdot \uv \in E}} &=
        \abs{C(\uv, \ua)} \cdot  \abs{\SLj{n-\l}} \cdot \abs{\Mj{(n-\l)}{\l}} q^{-k(n-\l)} = \\
        &= q^{k(n-m) - k \l} \cdot \abs{P_n} \cdot q^{-k(n-\l)} = \abs{P_n} q^{-km}.
    \end{split}
    \end{align}
    Substituting \cref{eq:g_count_good_n} in \cref{eq:ergodic_thm_for_stabilizer} we obtain
    \[
        \mu(E) = q^{- k m} = \Haar^{k, 0}(E).
    \]
    The set of boxes is a $\pi$-system which generates the product $\sigma$-algebra, so $\mu = \Haar^{k, 0}$ as required.
\end{proof}

We establish the main theorem of this section by combining the preceding results.

\begin{proof}[Proof of \Cref{thm:xkl_measures}]
    We may assume without loss of generality that $k \leq \l$, for the two actions on $X^{k,\l}$ and $X^{\l,k}$ are equivariantly isomorphic. Let $\mu$ be an $\SLi$-invariant probability measure on $X^{k, \l} = X^{k, 0} \x X^{0, \l}$.
    By \Cref{prop:xk0_measures} the measure $\mu$ is necessarily a joining of the two systems $(X^{k, 0}, \Haar^{k, 0})$ and $(X^{0, \l}, \Haar^{0, \l})$.

    Let $\up \in X^{0, \l}$ be a regularly linearly independent functional.  Let $P = \Stab_\SLi(\up)$ be its stabilizer. By \Cref{prop:stabilizer_dual_uniquely_ergodic} any such subgroup $P$ acts uniquely ergodically on $X^{k, 0}$. 
    
    Since a functional $\up \in X^{0, \l}$ is regularly linearly independent $\Haar^{0,\l}$-almost surely,  \Cref{lemma:disjoint_if_stabilizer_uniquely_ergodic} applied with respect to the two systems $X = X^{0,\l}$ and $Y = X^{k,0}$ says that these systems are disjoint. Hence $\mu = \Haar^{k, 0} \x \Haar^{0, \l} = \Haar^{k, \l}$, as required.
\end{proof}

\subsection{Invariant measures on related spaces}
\label{subsection:invariant_measures_on_related_spaces}
In this section, we use the unique ergodicity of the system $X^{k, \l}$ to describe invariant measures on two additional families of spaces. These are $V^{k, \l}$, of which $X^{k, \l}$ is a subspace, and the product of the infinite Grassmannians $\Gr_k(V) \x \Gr_\l(\hV)$, which is a quotient of $X^{k, \l}$.
The results in this section will not be used in the sequel, and we discuss them for their intrinsic interest.

\subsubsection{Invariant measures on $V^{k, \l}$}
We associate to each vector $\uv \in V^k$ the linear map 
\[T_\uv : \FFq^k \to V, \quad T_\uv(\alpha) = \sum_{i=1}^k \alpha_i \uv^i.\] 
Set $\ker(\uv) = \ker(T_\uv) \leq \FFq^k$.
Define $\ker(\up) \leq \FFq^\l$ for each functional $\up \in \hV^\l$ in an analogous manner.
For a given pair of subspaces $U \leq \FFq^k$, $W \leq \FFq^\l$ we define
\[
X^{k, \l}_{U, W} = \setcomp{(\uv, \up) \in V^{k, \l}}{\text{$\ker(\uv) = U$ and $\ker(\up) = W$}}.
\]
Every such $X^{k, \l}_{U, W}$ is a Borel subset of the space $V^{k, \l}$. We have the disjoint partition
\[
    V^{k, \l} = \bigsqcup_{U \leq \FFq^k, W \leq \FFq^\l} X^{k, \l}_{U, W}.
\]
For example, note that $X^{k, \l}_{0, 0} = X^{k, \l}$.
\begin{claim}
    Let $k, \l \in \NN_0$ and $U \leq \FFq^k, W \leq \FFq^\l$. The subset  $X^{k, \l}_{U, W}$ is $\SLi$-invariant. It is isomorphic as a $\SLi$-space to $X^{d,e}$ where $d = k - \dim U$ and $e = \l - \dim W$.
\end{claim}

\begin{proof}
    The invariance of $X^{k, \l}_{U, W}$ is clear.
    Choose arbitrarily surjective linear maps $\pi_1 : \FFq^k \to \FFq^d$ and $\pi_2 : \FFq^\l \to \FFq^e$ with kernels $U, W$ respectively. 
    The induced map $(\pi_1^\NN, \pi_2^\NN) : V^{k, \l} \to V^{d, e}$ is $\SLi$-equivariant and is an isomorphism from  $X_{U, W}^{k, \l}$ to $X^{d, e}$.
\end{proof}

The isomorphism constructed in the above claim endows $X^{k, \l}_{U, W}$ with an $\SLi$-invariant probability measure, induced by the unique $\SLi$-invariant probability measure on $X^{d,e}$. Its extension to $V^{k, \l}$ will be denoted $\tau^{k, \l}_{U, W}$.
This immediately gives the following classification of $\SLi$-invariant measures on $V^{k, \l}$.

\begin{theorem}\label{theorem:invariant_measures_on_v}
    Fix $k, \l \in \NN_0$. 
    Let $\mu$ be an ergodic $\SLi$-invariant measure on $V^{k, \l}$.
   There are subspaces $U \leq \FFq^k$ and $W \leq \FFq^\l$ such that $\mu$ is supported on $X^{k, \l}_{U, W}$ and $\mu = \tau^{k, \l}_{U, W}$. 
\end{theorem}

\begin{proof}
    Since each $X^{k, \l}_{U, W}$ is a Borel $\SLi$-invariant subset, and they form a partition of $V^{k, \l}$, the ergodic measure $\mu$ must be supported on exactly one such subset. 
    This subset is isomorphic to $X^{d, e}$ for some $d \leq k$ and $e \leq \l$ by the above claim. By \Cref{thm:xkl_measures}, the subset $X^{d, e}$ is uniquely ergodic. 
\end{proof}

\subsubsection{Invariant measures on Grassmannians}
\label{subsub:Grassmannians}
Given a vector space $W$, the $k$-dimensional Grassmannian $\Gr_k(W)$ for some $k\in\NN$ is the set of subspaces of $W$ of dimension $k$. 
In the case where $W$ has more structure, such as a topology or the structure of an algebraic variety, the Grassmannian $\Gr_k(W)$ usually inherits the corresponding structure.

We are interested in Grassmannians of the vector spaces $V$ and $\hV$. In the context of this section, we will give a definition using the spaces $X^{k, \l}$.
Notice that the space $V^{k, \l}$, in addition to the left action by the group $\SLi$ considered up to now, admits a right action by the group $\GLj{k} \x \GLj{\l}$. If we view $V^{k, \l}$ as the space of matrices $ \M_{\infty \x k}(\FFq) \x \M_{\l \x \infty}(\FFq)$ then this right action is the appropriate matrix multiplication. We define the \emph{$k$-dimensional Grassmannian} by
\[
    \Gr_k(V) = X^{k, 0} / (\GLj{k} \x 1).
\]
The projection map $X^{k, 0} \to \Gr_k(V)$ is the natural map sending a tuple $\uv \in X^{k, 0}$ to the $k$-dimensional subspace which it spans.
Since $X^{k, 0}$ is a standard Borel space, the Grassmannian $\Gr_k(V)$ is also standard Borel, being a quotient of $X^{k, 0}$ by the action of a finite group.
The Grassmannian $X^{k, 0}$ admit an $\SLi$-action induced from that of $X^{k, 0}$.

We can also define the dual Grassmannian $\Gr_\l(\hV) = X^{0, \l} / (1 \x \GLj{\l})$, which differs from $\Gr_\l(V)$ only by the $\SLi$-action on it.

Notice that the invariant random subgroup $\mu_{k, \l}$ of the group $\SLi$ can be alternatively defined as the stabilizer distribution of a random element in $\Gr_k(V) \x \Gr_\l(\hV)$, according to the measure induced from $\Haar^{k, \l}$. Therefore it is interesting to ask what the $\SLi$-invariant measures on this space are. We will show that this action is in fact uniquely ergodic.

\begin{proposition}\label{prop:grassmannian_measures}
    Let $k, \l \in \NN_0$. The only $\SLi$-invariant Borel probability measure on $\Gr_k(V) \x \Gr_\l(\hV)$ is the projection of the measure $\Haar^{k, \l}$ from $X^{k,\l}$.
\end{proposition}

The above proposition will be an immediate corollary of the following general result in the ergodic theory of group actions.

\begin{proposition}
\label{prop:unique_ergodic_quotient_amenable_uniquely_ergodic}
    Let $G$ be a locally compact amenable group acting continuously on a pair of Polish spaces $X$ and $Y$.
    Let $p : X \to Y$ be a surjective $G$-equivariant continuous map and $p_* : \cM_G(X) \to \cM_G(Y)$ the pushforward map on the spaces of $G$-invariant Borel probability measures. If the fibers of the map $p$ are compact 
    then the map $p_*$  is surjective.
\end{proposition}

\begin{proof}
    We use Zimmer's theory of amenable actions \cite{zimmer1978amenable_paper, zimmer1984ergodic_theory}. 
    Indeed \cite[Remark 6(a)]{zimmer1984ergodic_theory} states the analogous result, in the case where $X = Y \x F$ is a direct product of $Y$ with a compact fiber space $F$. 
    
    In our case, let $\nu$ be a $G$-invariant probability measure on $Y$. Consider the Borel bundle of Banach spaces over $Y$ where the fiber over every point $y \in Y$ is $C(X_y)$, the space of continuous functions on the non-empty compact fiber  $X_y = p^{-1}(y)$. Set $A_y$ to be the set of probability measures on $X_y$. This defines an affine $G$-space $A_y \subseteq C(X_y)^*$ over $Y$ in the sense of \cite[Definition 4.3.1]{zimmer1984ergodic_theory}.

    Since the group $G$ is amenable, it acts amenably on the probability space $(Y, \nu)$. Hence there exists a $\nu$-a.s.\ defined $G$-invariant section $\lambda_y$ of the dual bundle. Namely $\lambda_y \in A_y$ and every element $g \in G$ satisfies $g_* \lambda_y = \lambda_{g y}$ for $\nu$-a.e. $y \in Y$. This allows us define a measure $\mu$ on the space $X$ by 
    \[
        \mu = \int_{Y} \mu_y \diff \nu(y).
    \]
    Since the measure $\nu$ is $G$-invariant, we have for every $g \in G$ that 
    \[
        g_* \mu = \int_Y \lambda_{g y} \diff \nu(y) =  \int_Y \lambda_y \diff \nu(y) = \mu. 
    \]
    So that the measure  $\mu$ is $G$-invariant. Its projection to $Y$ is $p_*\mu = \int_Y \delta_{y} \diff \nu(y) = \nu$, as required.
\end{proof}

In particular, in the setting of \Cref{prop:unique_ergodic_quotient_amenable_uniquely_ergodic},   the unique ergodicity of the action of $G$ on $X$ implies the same for its action on $Y$.

\begin{proof}[Proof of \Cref{prop:grassmannian_measures}]
    The $\SLi$-action on $X = X^{k, \l}$ is uniquely ergodic by \Cref{thm:xkl_measures}. The fiber of the natural projection map $X^{k,\l} \to \Gr_k(V) \x \Gr_\l(\hV)$ is finite, hence compact.
    We may apply \Cref{prop:unique_ergodic_quotient_amenable_uniquely_ergodic} in this situation to deduce that $Y = \Gr_k(V) \x \Gr_\l(\hV)$ is uniquely ergodic. The desired conclusion follows.
\end{proof}

\section{Invariant random subgroups of \texorpdfstring{$\SLi$}{\SLipdf} and \texorpdfstring{$\GLi$}{\GLipdf}}
\label{section:proof_end}
Fix a prime power $q$.
We are now in a position to finish the classification of the invariant random subgroups of the group $\SLi$. The additional case of the group $\GLi$ will come as a corollary.

\subsection{Invariant random subgroups of \texorpdfstring{$\SLi$}{\SLipdf}}
Recall that $V = \prod_{n \in \NN} \FFq$ is equipped with the standard action of $\GLi$, and $\widehat{V} = V$ is equipped with the dual action of the same group.
For each $k, \l \in \NN_0$ we consider the product space $V^{k, \l}=V^k \x \hV^\l$ as a compact abelian group with $\GLi$-invariant Haar probability measure denoted $\Haar^{k,\l}$.
The invariant random subgroup $\mu_{k, \l} \in \IRS(\SLi)$ is defined to be the pushforward of the measure $\Haar^{k,\l}$ under the stabilizer map $\Stab_\SLi : V^{k, \l} \to \Sub{\SLi}$. In addition, we define $\mu_\infty = \delta_\triv$. 

In this section we calculate the Vershik trace of the invariant random subgroup $\mu_{k, \l}$ and characterize the topology of the space of all such measures.  
After that, we use the previous sections to show that any ergodic invariant random subgroups of $\SLi$ is either $\mu_{k, \l}$ or $\mu_\infty$. 
Note that $\mu_{k, \l}$ is ergodic, being the pushforward of the ergodic measure $\Haar^{k, \l}$ under an $\SLi$-equivariant map.

\subsubsection{Topology and Vershik traces of the $\mu_{k, \l}$}\label{subsubsection:sli_measures_topology}
By Skudlarek's theorem (\Cref{thm:sl_infty_characters})  the characters of the group $\SLi$ are 
\[\Ch(\SLi) = \setcomp{\chi_r}{r \in \NN_0} \union \set{\chi_\infty}\]
where for each $r \in \NN_0$ and all $g \in \SLi$
\[\chi_r(g) = q^{-r \cdot \rkt g} = q^{-r \cdot\rk(g-I)}\]
and $\chi_\infty = \delta_{I}$.

\begin{proposition}\label{prop:trace_calculation}
    For each $k,\l \in \NN_0$ the Vershik trace of the measure $\mu_{k, \l}$ is the character $\chi_{k + \l}$.
\end{proposition}

\begin{proof}
    Fix $k,\l \in \NN_0$. Let $g \in \SLi$ be any element. Write $$F_g = \setcomp{H \in \Sub{\SLi}}{g \in H}.$$ 
    Choose a point $(v^1, \ldots, v^k, \p^1, \ldots, \p^\l) \in V^{k, \l}$ uniformly at random. The value $\mu(F_g)$ is the probability that $g$ stabilizes the vectors $v^i$ and functionals $\p^j$. Note that the condition $g v^i = v^i$ is equivalent to $v^i \in \ker(g - I)$ for each $i$. This  happens with probability 
    \[q^{-\codim_V(\ker(g-I))} = q^{-\rk(g-I)}.\]
    The condition $(g^{-1})^{t} \p^j = \p^j$ is equivalent to $g^t \p^j = \p^j$ for each $j$. This  has probability \[q^{-\rk(g^t - I)} = q^{-\rk(g-I)}.\] 
    Overall, we get \begin{equation}\label{eq:mu_kl_character}
        \chi_{\mu_{k, \l}}(g) = \mu_{k, \l} (F_g) = q^{-(k + \l)\cdot  \rkt g} = \chi_{k + \l}(g)
    \end{equation}
    as required.
\end{proof}

Consider the subset $M = \setcomp{\mu_{k, \l}}{k, \l \in \NN_0} \union \set{\mu_\infty} \subseteq \IRS(\SLi)$ endowed with the weak-$*$ topology. In  \Cref{thm:main_theorem} it will be established that  $M=\Extr(\IRS(\SLi))$.

\begin{proposition}\label{prop:extreme_points_topology}
    Let $X = (\NN_0 \x \NN_0) \union \set{\infty}$ be the one-point compactification of $\NN_0 \x \NN_0$. The natural map $X \to M$ sending $(k, \l)$ to $\mu_{k, \l}$ and $\infty$ to $\mu_\infty$ is a homeomorphism.
\end{proposition}

\begin{proof}
    Denote by $h : X \to M$ the map in the statement.
    Recall that $E_{ij}$ are the elementary matrices in $\SLi$.
    One can check directly that
    \[
        \codim_V(\ker(E_{21}-I) \inter \ker(E_{23}-I)) = 2 \quad \text{and} \quad \codim_V(\ker(E_{12}) - I) \inter \ker(E_{32} -I)) = 1.
    \]
    We denote $F_g = \setcomp{H \in \Sub{\SLi}}{g \in H}$ for each element $g \in \SLi$. The computation in the proof of \Cref{prop:trace_calculation} shows that
    \begin{equation}\label{eq:mu_kl_multi_character}
        \mu_{k, \l}(F_{t_1} \inter F_{t_2}) = q^{-(2 k + \l)}.
    \end{equation}
    The two equations (\ref{eq:mu_kl_character}) and (\ref{eq:mu_kl_multi_character}) together show that the two parameters $k$ and $\l$ are determined by the measure $\mu_{k,\l}$. Hence the map $h$ is injective.
    
    To show that the map $h$ is continuous it remains to prove that $\mu_\infty = \lim_{(k, \l) \to \infty} \mu_{k, \l}$. By the Portmanteau theorem \cite[Theorem 2.1]{billingsley1999convergence} it will suffice to show that  all open subsets $U \subseteq \Sub{\SLi}$ satisfy 
    \[\liminf_{(k, \l) \to \infty} \mu_{k, \l}(U) \geq \mu_\infty(U).\]
    If $\triv \notin U$ then this inequality holds trivially. Otherwise, there exist some elements $g_1, \ldots, g_m \in \SLi \setminus \set{I}$ such that $\comp{F}_{g_1} \inter \cdots \inter \comp{F}_{g_m} \subseteq U$. In this case
    \begin{align*}        
        1 - \mu_{k, \l}(U) &\leq 
        \mu_{k, \l}\left(\comp{\left( \comp{F}_{g_1} \inter \cdots \inter \comp{F}_{g_m} \right)}\right)
        = \mu_{k, \l}\left( \bigunion_{i=1}^m F_{g_i} \right) \leq
        \\&\leq \sum_{i=1}^m \mu_{k, \l} (F_{g_i})
        = \sum_{i=1}^m q^{-(k + \l) \rkt g_i} \xrightarrow[(k, \l) \to \infty]{} 0
    \end{align*}
    as required.
    We conclude that $h$ is a continuous bijection. Since $X$ is compact and $M$ is Hausdorff $h$ is a homeomorphism.
\end{proof}

\begin{remark}
Once \Cref{thm:main_theorem} is established, we will deduce that the Vershik trace map $\IRS(\SLi) \to \Tr(\SLi)$ (\Cref{def:vershik_trace}) sends extreme points to extreme points. It is surjective and not injective, although it is finite-to-one on the extreme points.    
\end{remark}

\subsubsection{Classification of invariant random subgroups of $\SLi$}\label{subsubsection:main_thm_proof}

Let $\mu$ be a non-trivial ergodic invariant random subgroup of $\SLi$. Let $r = r(\mu)$ be the degree of $\mu$.

In \Cref{section:generalities_on_colors} we have defined for any $
\mu$-generic subgroup $H $ its associated set of colors $\Col(H)$ such that $\abs{\Col(H)} =q^r$. Here and in the following, a $\mu$-generic subgroup means one 
that satisfies the conclusions of Propositions \ref {prop:action_properties} and \ref{prop:no_fixed_points}. 

Further, we have defined a local and permutational (\Cref{def:permutational_action}) action of the group $\GLi$ on the set of color vectors $\Col(H)^\NN$. 
This action satisfies the assumptions of \Cref{thm:action_classification}. Hence there are $k, \l \in \NN_0 $ with $k + \l = r$ such that this $\GLi$-action is equivariantly isomorphic to the action of $\GLi$ on the space $V^{k, \l}$. By the ergodicity of the measure $\mu$, the quantities $k$ and $\l$ associated to the subgroup $H$ are $\mu$-a.s.\ constant.

Let $\Pi$ be the finite group of permutations $\varpi : \FFq^r \to \FFq^r$ such that the product map $\varpi^\NN : V^{k, \l} \to V^{k, \l}$ is an isomorphism of $\GLi$-actions. The group $\Pi$ is isomorphic to the direct product $\GLj{k} \x \GLj{\l}$, acting on the finite-dimensional vector space $\FFq^r = \FFq^k \x \FFq^\l$ by the product of the dual action and the standard action. 

For each $\mu$-generic subgroup $H$ denote by $\Pi_H$  the finite set of bijections $\pi : \Col(H) \to \FFq^k \times \FFq^\l = \FFq^r$ such that $\pi^\NN : \Col(H)^\NN \to V^{k, \l}$ is an isomorphism of $\GLi$-actions. It is a consequence of  \Cref{thm:action_classification}  that $P_H$  is non-empty. 
The finite group $\Pi$ admits a right action on $P_H$ by  composition. This action is simply transitive and hence $\abs{\Pi_H} = \abs{\Pi}$. In particular the cardinality of $\Pi_H$ does not depend on $H$.

We denote by $\WSub{\SLi}$  the measurable principal $\Pi$-bundle over $\Sub{\SLi}$ whose fiber over each $\mu$-generic subgroup $H $ is the set $\Pi_H$. The action of the finite group $\Pi$ on each set $\Pi_H$ extends to a measurable fiberwise action on $\WSub{\SLi}$.
The space $\WSub{\SLi}$ carries the probability measure $\hmu$, induced by the measure $\mu$ on the base and by the uniform probability measures on the fibers. 
Let $p : \WSub{\SLi} \to \Sub{\SLi}$ denote the projection. 

The action of the group $\SLi$ on the Chabauty space $\Sub{\SLi}$ can be lifted to a measurable action on the bundle $\WSub{\SLi}$ preserving the measure $\hmu$. Namely, let $g \in \SLi$ be an element and $(H,\pi) \in \WSub{\SLi}$ be a point of the bundle, where $H$ is a $\mu$-generic subgroup  and $\pi \in \Pi_H$ is a bijection. Recall from \Cref{section:generalities_on_colors} that the colors of $g H g^{-1}$ can be naturally 
identified with the colors of $H$. This allows us define the action of the element $g$ by $g \cdot (H, \pi) = (g H g^{-1}, \pi')$
where $\pi' : \Col(g H g^{-1}) \cong \Col(H) \xrightarrow{\pi} \FFq^r$ is the composition of this color identification with the bijection $\pi$.

Consider the $\hat \mu$-measurable \emph{fixed point map} $\FP: \WSub{\SLi} \to V^{k,\l}$  defined as follows.
For each point of the bundle consisting of a $\mu$-generic subgroup $H$ and a bijection $\pi \in \Pi_H$ let
\[
    \FP(H, \pi) = \pi^\NN(\col_e(H)).
\]
Recall that $\col_e(H)$ is the color vector of $H$ in the standard basis $e$.
We will see in the course of the proof of \Cref{thm:main_theorem} that the map $\FP$ selects a fixed point of $H$ in its action on $V^{k, \l}$ (hence its name).
The next proposition is a corollary of our measure classification result  \Cref{thm:xkl_measures}.
\begin{proposition}\label{prop:compatibility_of_measures}
    Consider the diagram
    \[
        \begin{tikzcd}
            (\WSub{\SLi}, \hmu) \arrow[r, "\FP"] \arrow[d, "p"]& (V^{k,\l}, \Haar^{k, \l}) \\ 
            (\Sub{\SLi}, \mu)
        \end{tikzcd}
    \]
    The two maps are $\SLi$-equivariant. We have $p_* \hmu = \mu$ and $\FP_* \hmu = \Haar^{k, \l}$.
\end{proposition}
\begin{proof}
    Both statements for the map $p$ follow directly from its definition. We study the map  $\FP$. Let $g \in \SLi$ be any element. Let $(H,\pi) \in \WSub{\SLi}$ be a point of the bundle, where $H$ is a $\mu$-generic subgroup  and $\pi \in \Pi_H$ is a bijection. 
    We have
    \begin{align*}
        g \cdot \FP(H, \pi) &= g \cdot \pi^\NN (\col_e(H))
        = \pi^\NN (\Theta_H(g) \col_e(H))
        = \pi^\NN (\col_e(g H g^{-1})) \\
        &= \FP(g \cdot (H, \pi)).
    \end{align*}
    The third equality in the above uses \Cref{prop:conjugate_action}. This shows that the map $\FP$ is $\SLi$-equivariant.
    
    The last step of the proof is to calculate the pushforward measure $\FP_* \hmu$. Notice that $\mu$-almost surely, every color of $H$ appears in the color vector $\col_e(H)$, that is $\setcomp{\col_e(H)_i}{i \in \NN} = \Col(H)$.
    In particular, $\hat \mu$-almost surely the bijection $\pi \in \Pi_H$ is such that the tuple $\pi^\NN(\col_e(H)) \in V^{k, \l}$ is  linearly independent, i.e.\ it lands in the subset $X^{k, \l} \subseteq V^{k, \l}$.
    So the probability measure $\FP_* \hmu$ is supported on $X^{k, \l}$ and $\SLi$-invariant. We obtain $\FP_* \hmu = \Haar^{k, \l}$ by \Cref{thm:xkl_measures}.
\end{proof}

We now prove our main result.

\begin{proof}[Proof of \Cref{thm:main_theorem}]
    To ease our notations, we write $G = \SLi$ throughout this proof.
    Let $\mu \in \IRS(G)$ be a non-trivial ergodic invariant random subgroup.
    We will show that the following diagram commutes almost surely.
    \[
        \begin{tikzcd}
            (\WSub{G}, \hmu) \arrow[r, "\FP"] \arrow[d, "p"] \arrow[dr] & (V^{k,\l}, \Haar^{k, \l}) \arrow[d] \\
            (\Sub{G}, \mu) \ar[r, shift right=.75ex, swap, "\id"] \ar[r, shift left=.75ex, "N_G"] & (\Sub{G}, \mu)
        \end{tikzcd}
    \]
    The objects appearing in the diagram (and in particular the indices $k,\l \in \NN_0$) have all been introduced in course of the above discussion. The two unlabeled arrows are stabilizer maps, and $N_G$ is the normalizer map.
    We need to show the following three statements a.s:
    \begin{enumerate}
        \item $N_G(H) = H$, 
        \item $N_G \mathrel{\circ} p = \Stab_G$,  and
        \item $\Stab_G \mathrel{\circ} \FP =  \Stab_G$.
    \end{enumerate}
    Let us take into account the following elementary \Cref{lemma:as_eq_functions} with respect to the space $Y = \Sub{G}$ and the collection $\cE = \setcomp{F_g}{g \in G}$,  where  $F_g = \setcomp{H \in \Sub{G}}{g \in H}$.
    Hence, the above three statements respectively are equivalent to showing for each element $g \in G$ that
    \begin{enumerate}[label=(\arabic*')]
        \item $\ssetcomp{H \in \Sub{G}}{g \in N_G(H)} = \ssetcomp{H \in \Sub{G}}{g \in H}$, 
        \item $\ssetcomp{(H, \pi) \in \WSub{G}}{g \in N_G(H)} = \ssetcomp{(H, \pi) \in \WSub{G}}{g \in \Stab_G(H, \pi)}$, and
        \item $\ssetcomp{(H, \pi) \in \WSub{G}}{g \in \Stab_G (\FP (H, \pi))} = \ssetcomp{(H, \pi) \in \WSub{G}}{g \in \Stab_G(H, \pi)}$ 
    \end{enumerate}
    up to $\mu$-null sets.
    Notice that containment in the direction $\supseteq$ clearly holds in all three statements. So it will suffice to show in each case that the two events have the same measure.
    Notice that $g$ stabilizes $(H, \pi)$ if and only if it stabilizes $H = p(H, \pi)$, i.e.\ $g \in N_G(H)$. This shows (2').
    We will show the other two statements simultaneously. Indeed, we have
    \begin{align*}
        q^{-r \cdot \rkt g } &= \Haar^{k, \l} \mleft( \ssetcomp{v \in V^{k,\l}}{g \in \Stab_G(v)} \mright) =
        \hmu \mleft( \ssetcomp{(H, \pi) \in \WSub{G}}{g \in \Stab_G(\FP(H, \pi))} \mright)  \\
        &\geq \hmu \mleft( \ssetcomp{(H, \pi) \in \WSub{G}}{g \in \Stab_G ((H, \pi))} \mright) =
        \hmu \mleft( \ssetcomp{(H, \pi) \in \WSub{G}}{g \in N_G(H)} \mright) \\ 
        &= \mu \mleft( \setcomp{H \in \Sub{G}}{g \in N_G(H)}\mright) \geq \mu \mleft( \setcomp{H \in \Sub{G}}{ g \in H} \mright) = q^{-r \cdot \rkt g}. 
    \end{align*}
    The first equality follows from \Cref{prop:trace_calculation}. The second equality on the first line and the first  on the bottom line follow from \Cref{prop:compatibility_of_measures}. The  last equality follows from \Cref{prop:num_of_colors_constant}. Equality must hold throughout. The remaining statements (1') and (3') follow and the above diagram commutes.
    To conclude the proof
    \[
       \mu_{k, \l} = (\Stab_G)_* \Haar^{k, \l} = (\Stab_G \circ \FP)_* \hmu = (\Stab_G)_* \hmu = (N_G \circ p)_* \hmu
       = (N_G)_* \mu = \mu
    \]
    as required.
\end{proof}

The following elementary measure theoretic lemma was used in the above proof of \Cref{thm:main_theorem}.

\begin{lemma}\label{lemma:as_eq_functions}
    Let $(X, \cA_X, \mu)$ be a measure space and $(Y, \cA_Y)$ be a measurable space. Let $\cE \subseteq \cA_Y$ be a countable collection of measurable subsets of $Y$ which separates points. 
    Let $f, g : X \to Y$ be two measurable functions such that 
    $ f^{-1}(E) = g^{-1}(E)$ up to a $\mu$-null set
    for all $E \in \cE$. Then $f$ and $g$ are equal $\mu$-a.s.
\end{lemma}
\begin{proof}
    The statement follows by noting that
    \begin{align*}
        \mu\mleft(\setcomp{x \in X}{f (x) \neq g(x)}\mright) &= \mu\mleft( \bigcup_{E \in \cE} \mleft( f^{-1}(E) \triangle g^{-1}(E) \mright) \mright) \leq 
         \sum_{E \in \cE} \mu \mleft( f^{-1}(E) \triangle g^{-1}(E) \mright) = 0
    \end{align*}
    where the first equality comes from the hypothesis that $\cE$ separates points.
\end{proof}
\subsection{Invariant random subgroups of \texorpdfstring{$\GLi$}{\GLipdf}}
We will show that every non-trivial ergodic invariant random subgroup of the group $\GLi$ is of the form \[\mu^K_{k, \l} = (\Stab_{\GLi}^K)_* \Haar^{k, \l}\]
where $k, \l \in \NN_0$, $K \leq \FFq^\x$ is a subgroup controlling the value of the determinant and $\Stab_{\GLi}^K : V^{k, \l} \to \Sub{\GLi}$ is the map 
\[
\Stab_{\GLi}^K(\ux) = \setcomp{A \in \GLi}{\text{$A \cdot \ux = \ux$ and $\det A \in K$}}.
\]
As the maps $\Stab_{\GLi}^K$ are $\GLi$-equivariant, the invariant random subgroups $\mu_{k, \l}^K$ are indeed ergodic.

\subsubsection{Topology and Vershik traces of the  $\mu_{k, \l}^K$}
We return to a discussion parallel to \Cref{subsubsection:sli_measures_topology} for the group $\GLi$. 
First, we state Skudlarek's classification of the characters of $\GLi$.
Let $\widehat{\FFq}^\x$ denote the Pontryagin dual group of the multiplicative group $\FFq^\x$. 
For each  character $\tau \in \widehat{\FFq^\x}$ and each $r \in \NN_0$, we define a map $\chi_r^\tau : \GLi \to \CC$ by
\[
\chi_r^\tau(g) = \tau(\det g) \cdot q^{-r \cdot \rkt g},
\]
where $\rkt g = \rk(g-I)$ as usual. In addition, we define $\chi_\infty = \delta_I$. 
\begin{theorem}[\cite{Skudlarek1976}] 
The characters of the finitary general linear group $\GLi$ are 
    \[
        \Ch(\GLi) = \setcomp{\chi_r^\tau}{r \in \NN_0, \ \tau \in \widehat{\FFq^\x}} \union \set{\chi_\infty}.
    \]
\end{theorem}

Given a subgroup $K \leq \FFq^\x$ there is a natural identification of the quotient group $\widehat{\FFq^\x / K}$ with the set of unitary  of $\FFq^\x$ which take the value $1$ on $K$. The analogous statement to \Cref{prop:trace_calculation} in the present setting is the following.

\begin{proposition}
    Fix $k, \l \in \NN_0$ and a subgroup $K \leq \FFq^\x$. The Vershik trace of the invariant random subgroup $\mu_{k, \l}^K$ is 
    \[
        \chi_{\mu_{k, \l}^K}(g) = \frac{\abs{K}}{q-1} \sum_{\tau \in \widehat{\FFq^\x / K}} \chi_{k + \l}^\tau (g) 
        = 
        \begin{cases}
            q^{-(k + \l) \cdot \rkt g} & \det g \in K \\ 
            0 & \det g \notin K
        \end{cases}
    \]
    for all $g \in G$.
\end{proposition}
\begin{proof}
    If $\det g \notin K$ then $\mu_{k, \l}^K\mleft(\setcomp{H}{g \in H}\mright) = 0$. If $\det g \in K$ then the calculation which establishes the formula on the right-hand side is the same as the calculation in the proof of \Cref{prop:trace_calculation}. The formula involving the summation over all unitary characters $ \tau \in \widehat{\FFq^\x / K}$ is equivalent (by means of standard harmonic analysis).
\end{proof}

Unlike the case of the group $\SLi$, the Vershik trace map $\tr : \IRS(\GLi) \to \Tr(\GLi)$ does not take extreme points to extreme points in general.
This result should also be compared with Thomas's calculation of the Vershik trace of the invariant random subgroups of $\Ai$ (\Cref{thm:ainfty_irs_vershik_character}).

The topology of the space $M_{\mathrm{GL}} = \setcomp{\mu_{k, \l}^K}{k, \l \in \NN_0, \ K \leq \FFq^\x} \union \set{\mu_\infty}$ is described in the following result. Once \Cref{thm:gl_inf_theorem} is established, it will follow from  it  that $M_{\mathrm{GL}}$ is  the space of ergodic invariant random subgroups of $\GLi$.

\begin{proposition}
    Let $X = (\NN_0 \x \NN_0 \x \Sub{\FFq^\x}) \union \set{\infty}$ be the one-point compactification of the space $\NN_0 \x \NN_0 \x \Sub{\FFq^\x}$. The map $X \to M_{\mathrm{GL}}$ sending $(k, \l, K)$ to $\mu_{k, \l}^K$ and $\infty$ to $\mu_\infty$ is a homeomorphism.
\end{proposition}
\begin{proof}
The proof is entirely analogous to that of \Cref{prop:extreme_points_topology}.    
\end{proof}

\subsubsection{{Classification of invariant random subgroups of $\GLi$}}

We will need a basic linear-algebraic property of the  product action on the space $V^{k,\l}$.

\begin{lemma}\label{lemma:fixed_vectors_of_stabilizer}
    Fix $k, \l \in \NN_0$. Let  $\ux = (\uv, \up) \in V^{k, \l}$. 
    If a vector $w \in V$ is stabilized by every matrix in $\Stab_\SLi(\ux)$ then $w \in \spn\set{\uv^1, \ldots, \uv^k}$. Dually, if a functional $\psi \in \hV$ is stabilized by every matrix in $\Stab_\SLi (\ux)$ then $\psi \in \spn\set{\up^1, \ldots, \up^\l}$.
\end{lemma}

\begin{proof}
    The two statements are dual to each other, and we will only show the first. 
    For each $n\in \NN$, we view the elements of $\hV_n$ as row vectors and those of $V_n$ as column vectors. So,  for each pair $\eta \in \hV_n$ and $u \in V_n$ the matrix product $\eta u \in \FFq$ corresponds to the pairing $\angles{u, \eta}$.

    We argue by contraposition. Let $w \notin \spn\set{\uv^1, \ldots, \uv^k}$ be any vector.
    Since this span is a finite set, there is an index $n > \l + 1$ large enough such that $w_{\leq n} \notin U_n$, where $U_n = \spn\set{\uv^1_{\leq n}, \ldots, \uv^k_{\leq n}} \leq V_n$.
    There exists a functional $\eta \in \hV_n$ which is zero on $U_n$ and non-zero on $w_{\leq n}$.
    Since $n > \l + 1$, there is a non-zero vector $u \in V_n$ satisfying the following collection of $\l + 1$ linear equations
    \[
        \text{$\up^j_{\leq n} u = 0$\, for all $j = 1, \ldots, \l$\, and\, $ \eta u = 0$}.
    \]

    Let $\tau = I + u \eta \in \GLn$. We may view $\tau$ as a matrix in $\GLi$ via the standard embedding  $\GLn \le \GLi$.
    Notice that
    \[
        (\tau - I)^2 = u (\eta u) \eta = 0.
    \]
    Therefore  $\tau$ is unipotent and in particular has determinant $1$ (it is even a transvection). 
    Finally,  notice that for $1 \leq i \leq k$ and $1 \leq j \leq \l$ we have 
    \[
        \tau \uv^i = \uv^i \quad \text{and} \quad \up^j \tau = \up^j.
    \]
    Hence  $\tau \in \Stab_\SLi(\ux)$. However
    \[
        \tau w = w + u (\eta w_{\leq n}) \neq w,
    \]
    which means that  $\tau$ does not stabilize $w$, as required.
\end{proof}


We extend our main result to the finitary general linear group.

\begin{proof}[Proof of \Cref{thm:gl_inf_theorem}]
Let $\mu$ be an ergodic invariant random subgroup of the group $\GLi$. The pushforward of $\mu$ under the intersection map $\Sub{\GLi} \to \Sub{\SLi}$ is an invariant random subgroup $\mu_0$ of $\SLi$. As the subgroup $\SLi$ is normal in $\GLi$, there is  a conjugation action of $\GLi$ on the Chabauty space $\Sub{\SLi}$. The measure $\mu_0$ is ergodic since the intersection map $\Sub{\GLi} \to \Sub{\SLi}$ is $\GLi$-equivariant.


    For any pair of indices $k, \l \in \NN_0$ notice that the stabilizer map $\Stab_{\SLi} : V^{k, \l} \to \Sub{\SLi}$ is $\GLi$-equivariant and the measure $\Haar^{k, \l}$ is $\GLi$-invariant.  Therefore the probability measures $\mu_{k, \l} \in \IRS(\SLi)$ are  $\GLi$-invariant and ergodic. 
    The same holds for the trivial invariant random subgroup $\mu_\infty = \delta_\triv$.
    By \Cref{thm:main_theorem} we know that $\mu_{k,\l}$ and $\mu_\infty$ are all the  $\SLi$-ergodic measures on $\Sub{\SLi}$.  Hence those are also all the  $\GLi$-ergodic measures on the same space.
     We conclude that $\mu_0$  is equal to either $\mu_{k, \l}$ for some $k, \l \in \NN_0$ or to $\delta_\triv$.

    In case $\mu_0 = \delta_{\triv}$ we claim that $\mu = \delta_\triv$ as well. Indeed, as the subgroup $\SLi$ has finite index in $\GLi$, the assumption $\mu_0 = \delta_\triv$ implies that $\mu$-a.e.\ subgroup is finite.
    Since $\GLi$ has the infinite conjugacy class property (ICC) and $\mu$ is a probability measure, this implies that $\mu = \delta_\triv$ as required. 

    We consider the remaining case, i.e. there are $k, \l \in \NN_0$ such that $\mu_0 = \mu_{k, \l}$.
Notice that the determinant map $\det : \Sub{\GLi} \to \Sub{\FFq^\x}$ is
$\GLi$-invariant and hence  a.s.\ constant by ergodicity of $\mu$. 
Let the subgroup $K \leq \FFq^\x$ be this constant value.
     We will show that $\mu = \mu_{k, \l}^K$.
     
    The fact that $\mu_0 = \mu_{k, \l}$ means that  for $\mu$-a.e.\ subgroup $H \in \Sub{\GLi}$ there is a tuple $\ux = (\uv, \up) \in V^{k, \l}$ such that
    \begin{equation}\label{eq:H_inter_sli}
        H \inter \SLi = \Stab_\SLi(\ux) = \setcomp{g \in \SLi}{g \cdot \ux = \ux}.
    \end{equation}
    As the property of regular linear independence (introduced in \Cref{def:regularly_linearly_independent})  is  $\Haar^{k, \l}$-generic, we may  assume that both $\uv$ and $\up$ are regularly linearly independent.
    We will show  that such a generic subgroup $H \in \Sub{\GLi}$ satisfies
    \[
        H = \Stab^K_\GLi (\ux) = \setcomp{A \in \GLi}{\text{$A \cdot \ux = \ux$\, and\, $\det A \in K$}}.
    \]

    We show containment in one direction. Consider a matrix $A \in H$.
    Observe that any matrix $B \in H$ having $\det B = \det A$ satisfies $B^{-1} A \in H \inter \SLi$. Hence $B^{-1} A$ fixes the point $\ux$, i.e.\
    \[
        A \cdot \ux = B \cdot \ux.
    \]
    Applying this observation to the conjugate matrices $B = g A g^{-1}$ where $g \in H \inter \SLi$ we get
    \[
        A \cdot \ux = g A g^{-1} \cdot \ux = g \cdot (A \cdot \ux).
    \]
    In other words, the tuple $A \cdot \ux$ is fixed by the subgroup $H \inter \SLi$. By \Cref{lemma:fixed_vectors_of_stabilizer} we deduce that $A \uv^i \in \spn \set{\uv^1, \ldots, \uv^k}$ and $(A^t)^{-1} \up^i \in \spn\set{\up^1, \ldots, \up^\l}$.
    
    Let $1 \leq i \leq k$ be a fixed index. Write $A \uv^i = \sum_{j=1}^k \alpha_j \uv^j$ where $\alpha_j \in \FFq$ are scalars.
   There is a sufficiently large index $ m\in \NN$ so that $A \in \SLj{m}$, and we get
    \[
        \uv^i_{> m} = (A \uv^i)_{> m} = \sum_{j=1}^k \alpha_j \uv^j_{> m}.
    \]
    As the tuple $\ux$ is regularly linearly independent, the sequence $(\uv^j_{> m})_{j = 1}^k$ is linearly independent. So $\alpha_j = \delta_{i j}$ and $A \uv^i = \uv^i$.
    The same argument applies to $\up^i$. Altogether this implies $A \cdot \ux = \ux$. Lastly, certainly $\det A \in K$ by the choice of $K$. This establishes containment in one direction.

    We show containment in the opposite direction. Let $A \in \GLi$ be a matrix which fixes the tuple $\ux$ and satisfies $\det A \in K$.
    Let $B \in H$ be any matrix with $\det B = \det A$. Such a matrix $B$  exists since $\det (H) = K$.
    By the containment in the forward direction,  the subgroup $B$ fixes the tuple $\ux$. So $B^{-1} A$ is a matrix which fixes $\ux$ and has determinant one. Hence $B^{-1} A \in H \inter \SLi$. We conclude that $A = B \cdot (B^{-1}A) \in H$, as required.

    To conclude, consider the $\mu_0$-measurable  map $\Psi : \Sub{\SLi} \to \Sub{\GLi}$ sending $\Stab_\SLi(\ux)$ to $\Stab^K_\GLi(\ux)$. The map $\Psi$ is well-defined, since given  $\ux = (\uv,\up)$ its stabilizer subgroup  $\Stab_\SLi(\ux)$ 
    determines the span of $\uv$ and of $\up$ (\Cref{lemma:fixed_vectors_of_stabilizer}).
    The above discussion shows that $H = \Psi(H \inter \SLi)$ for $\mu$-a.e.\ subgroup $H \in \Sub{\GLi}$. In other words $\mu = \Psi_* \mu_0$ by the definition of $\mu_0$. Therefore
    \[
        \mu = \left(\Psi \circ \Stab_\SLi \right)_* \Haar^{k, \l} = (\Stab^K_\GLi)_* \Haar^{k, \l} = \mu^K_{k, \l}
    \]
    as required.
\end{proof}

\section{Permutation stability of elementary enrichment}
\label{section:elementary-enrichment-stability}

In this section we give an explicit realization of the diagonal product $\bigotimes_i \PSLj{n_i}$ and prove \Cref{thm:elementary-enrichment-stability} saying this diagonal product is permutation stable. 

\subsection{Permutation stability}

We recall the definition of permutation stability. For each $n \in \NN$ let $S_n$ denote the symmetric group of degree $n$ with the bi-invariant normalized Hamming metric $d_n$ given by
\[
    d_n(\sigma, \tau) = 1 - \frac 1 n \abs{\Fix(\sigma^{-1} \tau)} \quad \forall \sigma, \tau \in S_n.
\]
Let $G$ be a finitely generated group. 
An \textit{almost-action} of $G$ is a sequence of set-theoretic maps $f_n : G \to S_n$ satisfying 
\[
    d_n (f_n (g) f_n(h),f_n(gh)) \xrightarrow[n \to \infty]{} 0 \quad \forall g, h \in G.
\]
Such an almost-action $f_n$ is said to be \textit{close to an action} if there is a sequence of group homomorphisms $\rho_n : G \to S_n$ satisfying
\[
    d_n(\rho_n(g), f_n(g)) \xrightarrow[n \to \infty]{} 0 \quad \forall g \in G.
\]

\begin{definition}
\label{def:permutation stable}
The group $G$ is called permutation stable if every almost-action of $G$ is close to an action.
\end{definition}

For finitely generated amenable groups, permutation stability can be understood in terms of properties of the space of invariant random subgroups; see \Cref{thm:blt_stability} in the introduction by Becker, Lubotzky and Thom. This is our main motivation in this section.
 
\subsection{Diagonal products}

We give an explicit realization of the diagonal product group $\bigotimes_i \PSLj{n_i}$ as a finitely generated dense subgroup of the profinite group $\prod \PSLj{n_i}$. 
We maintain notations consistent with \cite{levit_lubotzky2022} to make the parallel with the proof there as transparent as possible. In particular, these diagonal products will be denoted $G(\cP)$ here.

Let $\SLj{\ZZ}$ be the finitary special linear group over the finite field $\FFq$ with coordinates indexed by the set of integers $\ZZ$. This group is clearly isomorphic to the group $\SLi$ considered throughout this work, up to reindexing the coordinates. 

For each subset $I \subseteq \ZZ$ let $\PSLj{I}$ be the  projective linear group whose coordinates are indexed by $I$.
For each $r \in \NN$ denote $\setr{r} = \set{-r,\ldots,+r} \subset \NN$. For each $r \in \NN$ let  $\tau^+_r(a),\tau_r^-(a) \in \PSLj{\setr{r}}$ be the elementary matrices with value $a \in \FFq$ at its $(0,1)$-th or $(1,0)$-th entry respectively. Likewise, for each $r \in \NN$ let $ \sigma_r \in \PSLj{\setr{r}}$ be the permutation matrix associated to the full cycle $(-r,\ldots,r)$.

\begin{definition}

Let $\cP = (n_i)$ be a monotone sequence of odd integers with $n_1 \geq 3$. For each $i \in \NN$ write $n_i = 2r_i + 1$. Consider the profinite group
\[
    \PSL(\cP) = \prod_{i \in \NN} \PSLj{\setr{r_i}}.
\]
The group $G(\cP)$ is the  subgroup of $\PSL(\cP)$ generated by the elements $\sigma$ and $\tau^\pm(a)$ for all $a \in \FFq$  defined by setting $\tau^\pm(a)_i=\tau^\pm(a)_{r_i} \in \PSLj{\setr{r_i}}$ and $\sigma_i = \sigma_{r_i} \in \PSLj{\setr{r_i}}$ for each index $i \in \NN$.
\end{definition}

The group $G(\cP)$ has a rich subgroup structure. 
For instance, note that for each $r \in \NN$ the group $\PSL(\setr{r})$ is generated by the set of elementary matrices $\sigma_r$ and $\tau^\pm_r(a)$ for all $a \in \FFq$.  Hence the group $G(\cP)$ surjects onto each factor (and each finite product of factors) of $\PSL(\cP)$.
    For each $i \in \NN$ we denote by $L_i(\cP)$ the finite subgroup of $G(\cP)$ given by
    \[
        L_i(\cP) = \angles{\tau^+(a)^{\sigma^j}, \tau^-(a)^{\sigma^j} \; \middle| \; j \in \setr{r_i - 1}, a \in \FFq^\x }.
    \]
    The locally finite group $N(\cP)$ given by
    \[
        N(\cP) = \bigunion_{i \in \NN} L_i(\cP)
    \]
    is the commutator subgroup of $ G(\cP)$. It is normally generated in $G(\cP)$ by the elements $\tau^\pm(a)$ for all $a \in \FFq$.
    
Let $\bar \sigma$ be the outer automorphism of the group $\SLj{\ZZ}$ corresponding to conjugation by the infinitely supported full shift permutation given by $x \mapsto x+1$ for all $x \in \ZZ$. We define
$$ V = \ZZ \ltimes \SLj{\ZZ}$$
where $\ZZ \cong \left<\bar \sigma\right>$ and action on $\SLj{\ZZ}$ is via the automorphism $\bar \sigma$. There is a natural \emph{tail map} homomorphism $\tilde t : G(\cP) \to V$ defined by sending $\sigma$ to $\bar \sigma$ and $\tau^\pm(a)$ for each $a \in \FFq$ to the corresponding elementary matrix in
$\SLj{\ZZ}$. 
Denote $U(\cP) = \ker \widetilde t$. This normal subgroup of $G(\cP)$ satisfies 
\[ U(\cP) \cong \bigoplus_{i \in \NN} \PSLj{\setr{r_i}}  \le N(\cP).\]
The subgroup $N(\cP)$ is the kernel of the composition of the tail map $\widetilde t$ with the projection onto the subgroup $\ZZ$ of $V$, so that $G(\cP)/N(\cP) \cong \ZZ$.

\subsection{Invariant random subgroups of $V$}

Consider the standard action of the group $\SLz$ on the vector space $V_\ZZ = \FFq^\ZZ$ as well as the dual action on the vector space $\hV_\ZZ = V_\ZZ$ (c.f. \Cref{subsec:standard and dual}).  For each pair of indices $k, \l \in \NN_0$ we consider the vector space $V_\ZZ^{k, \l} = V_\ZZ^k \x \hV_\ZZ^\l$ endowed with the natural product action of $\SLz$ and with its Haar probability measure denoted $\Haar^{k, \l}$. In these terms, it follows from Theorem \ref{thm:main_theorem} that any non-trivial  invariant random subgroup of $\SLz$ is of the form 
$\mu_{k, \l}$, i.e. is obtained as the pushforward of $\Haar^{k, \l}$ under the stabilizer map $\Stab_\SLz : V_\ZZ^{k, \l} \to \Sub{\SLz}$.

The following lemma analyzes invariant random subgroups of $\SLz$. It corresponds to \cite[Corollary 4.5]{levit_lubotzky2022} and its argument  invokes \Cref{thm:main_theorem}. 

\begin{lemma}\label{lemma:no-cyclic-stabilizer}
    Let $\mu$ be an ergodic invariant random subgroup of $\SLz$.
    If there is some $m \in \NN$ such that $\mu$-almost every subgroup $H$ is normalized by $\bar\sigma^m g$ for  some element $g \in \SLz$ depending on the subgroup $H$ then $\mu$ is atomic.
\end{lemma}
\begin{proof}
    It suffices to verify that the invariant random subgroups $\mu_{k, \l} \in \IRS(\SLz)$ do not satisfy the aforementioned property provided that $k + \l \geq 1$. Consider a generic tuple $\ux = (\uv, \up) \in V_\ZZ^{k, \l}$ and let $H = \Stab_\SLz(\ux)$. Fix some $m \in \NN$. An element $g \in \SLz$ is such that the automorphism $\bar{\sigma}^m g$ normalizes $H$ if and only if the tuple $\bar{\sigma}^m g \cdot \ux$ is fixed by $H$. By \Cref{lemma:fixed_vectors_of_stabilizer}, this happens if and only if for all $i = 1, \ldots, k$ and $j = 1, \ldots, \l$ we have
    \[
        \bar{\sigma}^m g \uv^i \in \spn\set{\uv^1, \ldots, \uv^k} \quad \text{and} \quad \bar{\sigma}^m (g^{-1})^t \up^j \in \spn\set{\up^1, \ldots, \up^\l}.
    \]
    Given each $g \in \SLz$ and coefficient matrices $\alpha \in \Mj{k}{k}$ and $\beta \in \Mj{\l}{\l}$ we let $L(m, g, \alpha, \beta)$ be the set of tuples $\ux = (\uv, \up) \in V_\ZZ^{k, \l}$ such that for all $i = 1, \ldots, k$ and $j = 1, \ldots, \l$
    \[
        \bar{\sigma}^m g \uv^i = \sum_{a=1}^k \alpha_{ia} \uv^a \quad \text{and} \quad \bar{\sigma}^m (g^{-1})^t \up^j = \sum_{b=1}^\l \beta_{jb} \up^b.
    \]
    Clearly $L(m, g, \alpha, \beta)$ is a linear subspace of $V_\ZZ^{k, \l}$. We claim that this subspace is of infinite codimension provided $k + l \ge 1$.
    Indeed, given an element $g \in \SLz$ let $N \in \NN$ be sufficiently large such that $g \in \SLj{\setr{N-1}}$, i.e. $g$ acts as the identity on all coordinates of absolute value at least $N$. Then for any index $n \in \ZZ$ with $\abs{n} > N + m$, the equations on the $n$-th coordinate of $\ux \in V_\ZZ^{k, \l}$ to determine that $\ux$ belongs to $L(m, g, \alpha, \beta)$ are
    \[
    \uv^i_{n-m} = \sum_{a=1}^k \alpha_{ia} \uv^a_n \quad \text{and} \quad \up^j_{n-m} = \sum_{b=1}^\l \beta_{jb} \up^b_n.
    \]
    As $n$ ranges over indices of the form $n \in (m+1) \ZZ$ with $\abs{n} > N + m$ the resulting family of equations is linearly independent. Indeed, the equations for each such individual index  $n$ involve only the two coordinates $n - m$ and $n$,  and these coordinates are distinct  for different such $n$'s.
    We conclude that the solution space $L(m,g,\alpha,\beta)$ is of infinite codimension in $V_\ZZ^{k, \l}$ as claimed. 
    It follows that $\Haar^{k, \l}(L(m, g, \alpha, \beta)) = 0$.
    
    To conclude, note that by above paragraph every tuple $\ux \in V_\ZZ^{k, \l}$ whose stabilizer is normalized by $\bar{\sigma}^m g$ lies in $L(m, g, \alpha, \beta)$ for some pair of coefficient matrices $\alpha$ and $\beta$. As $\mu_{k, \l} = (\Stab_\SLz)_* \Haar^{k, \l}$, we get for every $m \in \NN$ that
    \begin{align*}
        \mu_{k, \l}&\mleft(\setcomp{H \in \Sub{\SLz}}{ \text{$\exists g \in \SLz$ such that $ \bar{\sigma}^m g$ normalizes $H$}}\mright) \leq \\
        &\leq \sum_{g \in \SLz} \sum_{\alpha \in \Mj{k}{k}} \sum_{\beta \in \Mj{\l}{\l}} \Haar^{k, \l}(L(m, g, \alpha, \beta)) = 0,
    \end{align*}
    as required.
\end{proof}





\begin{corollary}
\label{cor:rtimes-irs}
    Let $\tilde \mu$ be an ergodic invariant random subgroup of the group 
    $$V = \ZZ \ltimes \SLz \cong  \angles{\bar{\sigma}} \ltimes \SLz.$$ Then exactly one of the following two statements is true:
    \begin{enumerate}
        \item $\tilde \mu$-a.e.\ subgroup $H$ of $V$ is contained in the proper normal subgroup $\SLz$, or
        \item $\tilde \mu$ is the atomic probability measure supported on the finite-index subgroup $\left<\bar{\sigma}^m\right> \ltimes \SLz$ for some $m \in \NN$.
    \end{enumerate}
\end{corollary}

\begin{proof}
The proof of the analogous fact in the case of B.H. Neumann groups is \cite[Corollary 4.6]{levit_lubotzky2022}. The proof in the present case is exactly the same, where we make use of Lemma \ref{lemma:no-cyclic-stabilizer} in place of the analogous \cite[Corollary 4.5]{levit_lubotzky2022}.
\end{proof}

\subsection{Invariant random subgroups of $G(\cP)$}

The diagonal product $\bigotimes_i \PSLj{n_i}$ associated to the sequence $\cP = (n_i)$ and discussed in the introduction is isomorphic to the group $G(\cP)$ defined above. In view of \Cref{thm:blt_stability}, to show that the diagonal product in question $\bigotimes_i \SLj{n_i}$ is permutation stable, we are required to prove the following.

\begin{theorem}\label{thm:elementary-enrichment-irs}
All invariant random subgroups of $G(\cP)$ are co-sofic.
\end{theorem}

Recall that an invariant random subgroup $\mu$ is called co-sofic if $\mu$ is a weak-$*$ limit of invariant random subgroups supported on (i.e. giving measure one to) finite-index subgroup.

    \begin{proof}[Proof of \Cref{thm:elementary-enrichment-irs}]
    
    The proof is structurally identical to that of \cite[Theorem 1.3]{levit_lubotzky2022}. Rather than providing a full derivation, we highlight the necessary modifications to the definitions and indicate where \Cref{thm:main_theorem} is invoked. We  make use of exactly the same notations as in \cite{levit_lubotzky2022}.
    
    Fix a monotone sequence $\cP = (n_i)$ of odd integers with $n_1 \ge 3$ and consider the group $G = G(\cP)$. Denote $n_i = 2r_i+1$ for each $i \in \NN$. Let $\mu$ be an invariant random subgroup of the group $G$, which we wish to show is co-sofic. We may assume without loss of generality that $\mu$ is ergodic \cite[Corollary 2.6]{levit_lubotzky2022_metabelian}.

    Let $\tilde \mu = \tilde t_* \mu$ be the pushforward of the invariant random subgroup $\mu$ via the tail map $\tilde t$. Hence $\tilde \mu$ is an ergodic invariant random of the group $V = \ZZ \ltimes \SLz$. By \Cref{cor:rtimes-irs} there are exactly two possibilities for $\tilde \mu$.
    
    First, consider the case where $\tilde \mu$-a.e.\ subgroup $H$ is contained in the normal subgroup $\SLz \lhd V$, or equivalently, $\mu$-almost every subgroup is contained in the normal subgroup $N = N(\cP) \lhd G$.
    We define the following objects, parallel to the construction in \cite[Section 5]{levit_lubotzky2022}. For each $i \in \NN$ let 
    $$ D_i = D_i(\cP) = \ker\left(G \to \ZZ / (2r_i-1)\ZZ\right) $$
    via the quotient map $G/N \cong \ZZ$, and
    $$ G_i = G_i(\cP) = \ker\left(G \to \prod_{j =1}^i \PSLj{\setr{r_j}}\right)
    $$
    so that $D_i$ and $G_i$ are  finite-index normal subgroups of $G$. For each $i \in \NN$ denote 
    \[ F_i = F_i(\cP) = \left<\sigma^i \: : \: i\in \setr{2r_i-1}\right> \cdot L_i \]
    so that $F_i$ is a Folner sequence in $G$. Further, for a $\mu$-generic subgroup $H$ and for each $i \in \NN$ denote
    $$ K_i = K_i(H) =  (H \cap L_i)(G_i \cap D_i)$$
    so that $K_i$ is a finite-index subgroup of $G$. The subset $F_i$ is a finite union of disjoint coset transversals of the normalizer $N_G(K_i)$.
    The pairs $(K_i, F_i)$ form a \emph{Weiss approximation} for the $\mu$-generic subgroup $H$ in the sense of \cite[Definition 1.7]{levit_lubotzky2022_metabelian}.  The arguments presented there carry over to our case without modification. The existence of such a Weiss approximation in a finitely generated residually finite amenable group implies that the invariant random subgroup $\mu$ is co-sofic  \cite[Theorem 1.8]{levit_lubotzky2022_metabelian}.

    Lastly, consider the remaining case where $\tilde \mu$ is an atomic measure supported on the finite-index subgroup $\left<\bar \sigma^n\right> \ltimes \SLz$ for some $n\in \NN$.  It  follows from \cite[Proposition 2.5]{levit_lubotzky2022} that $\mu$-almost every subgroup has finite index in $G$ (while \cite[Proposition 2.5]{levit_lubotzky2022} is stated for B.H. Neumann groups, its parallel version in our setting is proved identically). Hence $\mu$ is certainly co-sofic in this case as well. This concludes the proof.
\end{proof}

\bibliographystyle{alpha}
\bibliography{references}

@article{Skudlarek1976,
  title={Die unzerlegbaren {C}haraktere einiger diskreter {G}ruppen},
  author={Skudlarek, Hans-Ludwig},
  journal={Mathematische Annalen},
  volume={223},
  number={3},
  pages={213--231},
  year={1976},
  publisher={Springer}
}

@article{borevich1984arrangement,
  title={Arrangement of subgroups in the general linear group over a commutative ring},
  author={Borevich, Zenon Ivanovich and Vavilov, Nikolay Aleksandrovich},
  journal={Trudy Matematicheskogo Instituta imeni VA Steklova},
  volume={165},
  pages={24--42},
  year={1984},
  publisher={Russian Academy of Sciences, Steklov Mathematical Institute of Russian~…}
}

@book{Kleidman_Liebeck_1990, 
  address={Cambridge},
  series={London Mathematical Society Lecture Note Series}, 
  title={The Subgroup Structure of the Finite Classical Groups}, 
  publisher={Cambridge University Press}, 
  author={Kleidman, Peter  and Liebeck, Martin }, 
  year={1990}, 
  volume={129}
}

@article{bowen2017characteristic,
  title={Characteristic random subgroups of geometric groups and free abelian groups of infinite rank},
  author={Bowen, Lewis and Grigorchuk, Rostislav and Kravchenko, Rostyslav},
  journal={Transactions of the American Mathematical Society},
  volume={369},
  number={2},
  pages={755--781},
  year={2017}
}

@article{glasner1997kazhdan,
  title={Kazhdan's property {T} and the geometry of the collection of invariant measures},
  author={Glasner, Eli and Weiss, Benjamin},
  journal={Geometric \& Functional Analysis GAFA},
  volume={7},
  number={5},
  pages={917--935},
  year={1997},
  publisher={Springer}
}

@article{Guralnick_Malle_2012,
  title={Products of conjugacy classes and fixed point spaces},
  author={Guralnick, Robert and Malle, Gunter},
  journal={Journal of the American Mathematical Society},
  volume={25},
  number={1},
  pages={77--121},
  year={2012}
}

@article{thoma1964unzerlegbaren,
  title={Die unzerlegbaren, positiv-definiten {K}lassenfunktionen der abz{\"a}hlbar unendlichen, symmetrischen {G}ruppe},
  author={Thoma, Elmar},
  journal={Mathematische Zeitschrift},
  volume={85},
  number={1},
  pages={40--61},
  year={1964},
  publisher={Springer}
}

@article{Vershik2012,
  author = {A. M. Vershik},
  journal = {Mosc. Math. J.},
  volume = {12},
  number = {1},
  pages={193--212},
  title = {Totally nonfree actions and the infinite symmetric group},
  year = {2012}
}

@article{vershik1974ergodic_thm,
  author       = {Vershik, A. M.},
  title        = {Description of invariant measures for the actions of some
                  infinite-dimensional groups},
  journal      = {Dokl. Akad. Nauk SSSR},
  year         = {1974},
  volume       = {218},
  number       = {4},
  pages        = {749--752},
  url          = {http://mi.mathnet.ru/dan38565},
}

@article{thomas2018_ainfty_irs,
  title={Invariant random subgroups of inductive limits of finite alternating groups},
  author={Thomas, Simon and Tucker-Drob, Robin},
  journal={Journal of Algebra},
  volume={503},
  pages={474--533},
  year={2018},
  publisher={Elsevier}
}

@article{thomas2022_sinfty_irs_and_characters,
  title={Characters and invariant random subgroups of the finitary symmetric group},
  author={Thomas, Simon},
  journal={Advances in Mathematics},
  volume={399},
  pages={108272},
  year={2022},
  publisher={Elsevier}
}

@article{thoma1964_choquet_simplex,
  title={{\"U}ber unit{\"a}re {D}arstellungen abz{\"a}hlbarer, diskreter {G}ruppen},
  author={Thoma, Elmar},
  journal={Mathematische Annalen},
  volume={153},
  number={2},
  pages={111--138},
  year={1964},
  publisher={Springer}
}

@article{lindenstrauss2001ergodic_theorem,
  title={Pointwise theorems for amenable groups},
  author={Lindenstrauss, Elon},
  journal={Inventiones mathematicae},
  volume={146},
  number={2},
  pages={259--295},
  year={2001},
  publisher={Springer}
}

@book{glasner2003ergodic,
  title={Ergodic Theory via Joinings},
  author={Glasner, Eli},
  series={Mathematical Surveys and Monographs},
  volume={101},
  year={2003},
  publisher={American Mathematical Society},
  address={Providence, RI}
}

@article{McLaughlin1967,
  author={McLaughlin, Jack},
  title={Some groups generated by transvections},
  journal={Archiv der Mathematik},
  year={1967},
  month={9},
  day={01},
  volume={18},
  number={4},
  pages={364-368},
  issn={1420-8938},
  doi={10.1007/BF01898827},
  url={https://doi.org/10.1007/BF01898827}
}

@article{levit_lubotzky2022_metabelian, 
  title={Infinitely presented permutation stable groups and invariant random subgroups of metabelian groups}, 
  volume={42}, 
  DOI={10.1017/etds.2021.9}, 
  number={6}, 
  journal={Ergodic Theory and Dynamical Systems}, 
  author={Levit, Arie and Lubotzky, Alexander}, 
  year={2022}, 
  pages={2028-2063}
}

@article{levit_lubotzky2022,
  author={Levit, Arie
  and Lubotzky, Alexander},
  title={Uncountably many permutation stable groups},
  journal={Israel Journal of Mathematics},
  year={2022},
  month={12},
  day={01},
  volume={251},
  number={2},
  pages={657-678},
  issn={1565-8511},
  doi={10.1007/s11856-022-2425-0},
  url={https://doi.org/10.1007/s11856-022-2425-0}
}

@article{blt2019,
author = {Oren Becker and Alexander Lubotzky and Andreas Thom},
title = {{Stability and invariant random subgroups}},
volume = {168},
journal = {Duke Mathematical Journal},
number = {12},
publisher = {Duke University Press},
pages = {2207 -- 2234},
year = {2019},
doi = {10.1215/00127094-2019-0024},
URL = {https://doi.org/10.1215/00127094-2019-0024}
}

@article{dogon2026characters,
  title={Characters of diagonal products and {H}ilbert--{S}chmidt stability},
  author={Dogon, Alon and Levit, Arie and Vigdorovich, Itamar},
  journal={Transactions of the American Mathematical Society},
  volume={379},
  number={02},
  pages={979--1021},
  year={2026}
}

@article{Hestenes1970SingerGroups,
  title={Singer Groups},
  volume={22},
  DOI={10.4153/CJM-1970-057-2},
  number={3},
  journal={Canadian Journal of Mathematics},
  author={Hestenes,
  Marshall},
  year={1970},
  pages={492-513}
}

@article{Kingman1978RandomPartitions,
    author = {Kingman, J. F. C.},
    title = {The Representation of Partition Structures},
    journal = {Journal of the London Mathematical Society},
    volume = {s2-18},
    number = {2},
    pages = {374-380},
    year = {1978},
    month = {10},
    issn = {0024-6107},
    doi = {10.1112/jlms/s2-18.2.374},
    url = {https://doi.org/10.1112/jlms/s2-18.2.374},
    eprint = {https://academic.oup.com/jlms/article-pdf/s2-18/2/374/2788610/s2-18-2-374.pdf},
}

@book{phelps2002lectures,
  title={Lectures on {C}hoquet's theorem},
  author={Phelps, Robert},
  year={2002},
  publisher={Springer}
}

@book{bekka_la_harpe2020,
  author="Bekka,  Bachir and La Harpe, Pierre de",
  title="Unitary representations of groups, duals, and characters",
  publisher="American Mathematical Society",
  year="2020",
  series="Mathematical Surveys and Monographs",
  doi={10.1090/surv/250},
}

@article{AGV2014,
  author = {Mikl{\'o}s Ab{\'e}rt and Yair Glasner and B{\'a}lint Vir{\'a}g},
  title = {{Kesten's theorem for invariant random subgroups}},
  volume = {163},
  journal = {Duke Mathematical Journal},
  number = {3},
  publisher = {Duke University Press},
  pages = {465 -- 488},
  year = {2014},
  doi = {10.1215/00127094-2410064},
  URL = {https://doi.org/10.1215/00127094-2410064}
}

@article{Gorin_Kerov_Vershik2014,
  title = {Finite traces and representations of the group of infinite matrices over a finite field},
  journal = {Advances in Mathematics},
  volume = {254},
  pages = {331-395},
  year = {2014},
  issn = {0001-8708},
  doi = {https://doi.org/10.1016/j.aim.2013.12.028},
  url = {https://www.sciencedirect.com/science/article/pii/S0001870813004775},
  author = {Vadim Gorin and Sergei Kerov and Anatoly Vershik},
}

@article{VershikKerov_U_infty_1982,
  title={Characters and factor representations of the infinite unitary group},
  author={Vershik, A. M. and Kerov, S. V.},
  journal={Dokl. Akad. Nauk SSSR},
  volume={267},
  number={2},
  pages={272--276},
  year={1982}
}

@article{Voiculescu1976,
  author={Voiculescu, Dan},
  title = {Repr\'esentations factorielles de type {$\mathrm{II}_1$} de {$U(\infty)$}},
  journal={J. Math. Pures Appl},
  volume = {55},
  pages = {1--20},
  year = {1976}
}

@article{Boyer1983,
  author = {Robert  Boyer},
  title = {Infinite traces of {AF}-algebras and characters of {$U(\infty)$}},
  journal = {Journal of Operator Theory},
  volume = {9},
  number = {2},
  pages = {205--236},
  year = {1983},
  publisher = {Theta Foundation},
}

@article{Boyer1992,
  ISSN = {03794024, 18417744},
  author = {Boyer, Robert},
  journal = {Journal of Operator Theory},
  number = {2},
  pages = {281--307},
  publisher = {Theta Foundation},
  title = {Characters and factor representations of the infinite dimensional classical groups},
  volume = {28},
  year = {1992}
}

@incollection{HiraiHirai2003,
  author = {Hirai, Takeshi and Hirai, Etsuko},
  title = {Characters for the infinite {W}eyl groups of type {$B_\infty$}/{$C_\infty$} and {$D_\infty$}, and for analogous groups},
  booktitle = {Non-Commutativity, Infinite-Dimensionality and Probability at the Crossroads},
  publisher = {World Scientific},
  year = {2003},
  pages = {296--317},
  doi = {10.1142/9789812705242_0011},
  url = {https://www.worldscientific.com/doi/abs/10.1142/9789812705242_0011}
}

@misc{gelander2015lectureinvariantrandomsubgroups,
  title={A lecture on Invariant Random Subgroups}, 
  author={Tsachik Gelander},
  year={2015},
  eprint={1503.08402},
  archivePrefix={arXiv},
  primaryClass={math.GR},
  note={arXiv:1503.08402 [math.GR]},
  url={https://arxiv.org/abs/1503.08402}, 
}

@article{kantor1980linear,
  title={Linear groups containing a {S}inger cycle},
  author={Kantor, William},
  journal={Journal of Algebra},
  volume={62},
  number={1},
  pages={232--234},
  year={1980},
  publisher={Academic Press}
}

@article{guralnick1999linear,
  title={Linear groups with orders having certain large prime divisors},
  author={Guralnick, Robert and Penttila, Tim and Praeger, Cheryl E and Saxl, Jan},
  DOI={10.1112/S0024611599001616},
  journal={Proceedings of the London Mathematical Society},
  volume={78},
  number={1},
  pages={167--214},
  year={1999},
  publisher={Cambridge University Press}
}

@incollection{billingsley1999convergence,
  author = {Billingsley, Patrick},
  title = {Weak Convergence in Metric Spaces},
  booktitle = {Convergence of Probability Measures},
  publisher = {John Wiley {\&} Sons, Ltd},
  year = {1999},
  pages = {7--79},
  isbn = {9780470316962},
  doi = {10.1002/9780470316962.ch1},
  url = {https://onlinelibrary.wiley.com/doi/abs/10.1002/9780470316962.ch1}
}

@Article{thomasbook,
author={Thomas, Simon},
title={private communication}
}

@Article{thomas2018URS,
  author={Thomas, Simon},
  title={Uniformly recurrent subgroups of simple locally finite groups},
  journal={Israel Journal of Mathematics},
  year={2018},
  month={10},
  day={01},
  volume={228},
  number={1},
  pages={321-352},
  issn={1565-8511},
  doi={10.1007/s11856-018-1766-1},
  url={https://doi.org/10.1007/s11856-018-1766-1}
}

@incollection{zimmer1984ergodic_theory,
  author = {Zimmer, Robert},
  title = {Amenability},
  booktitle = {Ergodic Theory and Semisimple Groups},
  publisher = {Birkh{\"a}user Boston},
  address = {Boston, MA},
  year = {1984},
  pages = {59--84},
  isbn = {978-1-4684-9488-4},
  doi = {10.1007/978-1-4684-9488-4_4}
}

@article{zimmer1978amenable_paper,
title = {Amenable ergodic group actions and an application to {P}oisson boundaries of random walks},
journal = {Journal of Functional Analysis},
volume = {27},
number = {3},
pages = {350-372},
year = {1978},
issn = {0022-1236},
doi = {https://doi.org/10.1016/0022-1236(78)90013-7},
author = {Robert Zimmer}
}

@article{Leinen2007,
title = {Positive definite functions of finitary isometry groups over fields of odd characteristic},
journal = {Journal of Pure and Applied Algebra},
volume = {208},
number = {3},
pages = {1003-1021},
year = {2007},
issn = {0022-4049},
author = {Felix Leinen and Orazio Puglisi},
}

\end{document}